\documentclass{amsart}

\usepackage{amsmath,amssymb, xypic, color}
\usepackage{graphicx}
\usepackage{hyperref}

\usepackage{tikz-cd}

\newcommand{\arr}{\longrightarrow}
\newcommand{\R}{\mathbb{R}}
\newcommand{\Z}{\mathbb{Z}}
\newcommand{\C}{\mathbb{C}}

\newcommand{\X}{\mathcal{X}}

\newcommand{\M}{\mathcal{M}}
\newcommand{\alb}{\mathsf{X}}
\newcommand{\xs}{\alb^*}
\newcommand{\xo}{\alb^\omega}
\newcommand{\xmo}{\alb^{-\omega}}
\newcommand{\nuke}{\mathcal{N}}

\newcommand{\img}{\mathop{\mathrm{IMG}}}

\newcommand{\G}{\Gamma}

\newcommand{\lims}[1][G]{\mathcal{J}_{#1}}
\newcommand{\limg}[1][G]{\mathcal{X}_{#1}}
\newcommand{\bim}{\mathfrak{B}}
\newcommand{\til}{\mathcal{T}}

\newcommand{\be}{\mathsf{s}}
\newcommand{\en}{\mathsf{r}}

\newcommand{\cdim}{\mathop{\mathrm{c.dim}}}

\newcommand{\Br}{\mathbf{B}}

\newcommand{\mapdown}[1]%
{\Big\downarrow\rlap{$\vcenter{\hbox{$\scriptstyle#1$}}$}}

\newtheorem{theorem}{Theorem}[section]
\newtheorem{lemma}[theorem]{Lemma}
\newtheorem{proposition}[theorem]{Proposition}
\newtheorem{corollary}[theorem]{Corollary}

\theoremstyle{definition}
\newtheorem{defi}[theorem]{Definition}
\newtheorem{example}[theorem]{Example}

\title{Polytopal substitutions and conformal dimension}
\author{Volodymyr Nekrashevych}

\begin{document}
\maketitle

\begin{abstract}
We introduce polytopal substitutions, a recursive procedure for constructing polyhedral approximations of self-similar topological spaces. For a class of substitutions called sub-hyperbolic, we show that the limit of these approximations is naturally homeomorphic to the limit space of the associated iterated monodromy group. We establish a relationship between the $\ell^p$-contraction coefficients of sub-hyperbolic substitutions and the conformal dimension of their limit spaces.  In particular, we characterize the conformal dimension in terms of asymptotic moduli and relate these moduli to $p$-energies of graphs.
\end{abstract}

\tableofcontents
\section{Introduction}

We study a recursive procedure, called \emph{polytopal substitution}, that produces polyhetopal complexes approximating self-similar spaces. This procedure was introduced (in slightly different terms) in~\cite{nek:models} under the name ``topological automaton''; see also~\cite{nek:dyngroups} and~\cite{nek:dimension}. 

The main result of our paper is a characterization of the Ahlfors regular conformal dimension of the Julia sets of subhyperbolic branched coverings in terms of $\ell^p$-contraction coefficients of the polytopal substitutions approximating them. This generalizes to higher dimensions a result of~\cite{pigrimthurston:dimension}, which treats the one-dimensional case. Our approach also applies to subhyperbolic maps and complexes of groups, rather than only to hyperbolic maps.

A polytopal substitution is a pair $\sigma=(f, \iota)$ of maps $f, \iota\colon\G_1\arr\G_0$, where $f$ maps faces of $\G_1$ isometrically onto faces of $\G_0$, while $\iota$ maps faces of $\G_1$ into faces of $\G_0$ by affine (or continuous) maps. 
We regard the faces of $\G_0$ as labels, so that $f$ assigns a label to each face of $\G_1$. 

The pair $\sigma=(f, \iota)$ is interpreted as a rule of constructing a new $\G_0$-labeled complex from a given one. Namely, let $\Delta$ be a complex equipped with a labeling $F\colon\Delta\arr\G_0$, where $F$ maps faces of $\Delta$ isometrically to faces of $\G_0$. The new $\G_0$-labeled complex $\sigma(\Delta)$
is obtained by replacing every face $A$  of $\Delta$ by the labeled complex $\iota^{-1}(F(A))\subset\G_1$. 

More formally, the complex $\sigma(\Delta)$ is obtained by completing the pullback diagram
\[\begin{array}{ccc}
\sigma(\Delta) & \stackrel{F'}{\arr} & \G_1 \\
\mapdown{\iota'} & & \mapdown{\iota} \\
\Delta &\stackrel{F}{\arr} & \G_0 
\end{array}\]
and the labeling of the faces of $\sigma(\Delta)$ is given by the composition $f\circ F'\colon\sigma(\Delta)\arr\G_0$.

Iterating the substitution produces a sequence of labeled complexes
\begin{equation}
\label{eq:deltainfty}
\Delta\stackrel{\iota'}{\longleftarrow}\sigma(\Delta)\longleftarrow\sigma^2(\Delta)\longleftarrow\cdots
\end{equation}
We define $\sigma^\infty(\Delta)$ to be the inverse limit of this sequence in the category of topological spaces.

In particular, starting with $\Delta=\G_0$, we obtain a sequence of complexes $\G_0, \G_1, \G_2, \ldots$, where $\G_{n+1}$ is obtained from $\G_n$ by the substitution procedure above. Equivalently, these complexes are obtained by completing the pullback diagrams
\[
\begin{array}{ccc}
\G_{n+1} & \stackrel{f_n}{\arr} & \G_n\\
\mapdown{\iota_n} & & \mapdown{\iota_{n-1}} \\
\G_n & \stackrel{f_{n-1}}{\arr} & \G_{n-1}
\end{array}
\]
where $f=f_0$ and $\iota=\iota_0$.

The space $\sigma^\infty(\G_0)=\G_\infty$, called the \emph{limit space} of the substitution, is the inverse limit of the topological spaces $\G_n$ with respect to the maps $\iota_n\colon\G_{n+1}\arr\G_n$. The maps $f_n\colon\G_{n+1}\arr\G_n$ induce a continuous map
\[f_\infty\colon\G_\infty\arr\G_\infty,\]
which we call the \emph{limit dynamical system} of the substitution.

The notion of polytopal substitution generalizes, or is closely related to, several constructions appearing in the literature. For example, suppose that $\G_0$ is a bouquet of loops, $\G_1$ is a bouquet of cycles, the restriction of $f$ to each cycle is orientation-preserving, and $\iota$ is an orientation-preserving homeomorphism. Then the substitution $(f,\iota)$ naturally corresponds to a word substitution $\sigma\colon E\arr E^*$, where $E$ is the set of loops of $\G_0$, $E^*$ is the free semigroup generated by the set of edges of $\G_0$, and $\sigma(e)$ is the word describing the image of the corresponding cycle $\iota^{-1}(e)$ under $f$.

Subdivision rules constitute another special case. They are substitutions $f, \iota\colon\G_1\arr\G_0$ for which $\iota$ is a homeomorphism. Such rules arise in conformal dynamics~\cite{cfp:subdivision} and in the study of self-similar tilings~\cite{goodmanstrauss:tilings,anderputn,mozes:tilings}.

In dimension one, our construction is closely related to several notions of substitutional graphs. Examples include the graph substitutions in~\cite{previte:graphsubs} and the replacement rules in~\cite{belkforrest}, both of which can be viewed as special cases of the present framework.

A construction very close to our notion of polytopal substitution, together with a discussion of its relation to biology, appears in work of M. Gromov~\cite{gromov:cell90,gromov:cell13}. Related biological substitution mechanisms include Lindenmayer systems~\cite{Lindenmayer1968I,Lindenmayer1968II}; see also~\cite{PrusinkiewiczLindenmayer1990}.

An important distinction between our notion and similar more geometric notions from fractal geometry, such as iterated function systems, is that the complexes $\sigma^n(\Delta)$ are abstract polytopla complexes rather than complexes embedded in a fixed ambient space. Consequently, the limit space \(\sigma^\infty(\Delta)\) is initially only a topological space and does not come equipped with a natural metric.

As already observed by Gromov in~\cite{gromov:cell90}, an important condition for polytopal substitutions is their \emph{hyperbolicity}. It can be formulated by requiring that $\iota\colon\G_1\arr\G_0$ be contracting on every cell. Equivalently, under substitution every cell of $\G_0$ is replaced by a uniformly larger complex. If, in addition, the complexes $\sigma^n(\Delta)$ are compact and have uniformly bounded geometry, then the complex associated with the sequence~\eqref{eq:deltainfty} is Gromov-hyperbolic, and the space $\sigma^\infty(\Delta)$ is naturally homeomorphic to its boundary. The limit space therefore carries a natural class of visual metrics, and it is natural to consider the corresponding quasisymmetric gauge.

An especially well-behaved class of substitutions is given by \emph{contracting} (or \emph{hyperbolic}) \emph{virtual endomorphisms}. These are substitutions
 $f, \iota\colon\G_1\arr\G_0$ for which $\iota$ is strictly contracting on the faces, and $f$ is a covering. It was shown in~\cite{nek:models} that if $\G_0$ is connected, then the limit dynamical system $f_\infty\colon\G_\infty\arr\G_\infty$ of a hyperbolic virtual endomorphism is uniquely determined by the associated iterated monodromy group. In particular, the limit dynamical system depends only on the actions of $f$ and $\iota$ on the fundamental groups.

More generally, we can define contracting virtual endomorphisms of complexes of groups. We  call a polytopal substitution $f, \iota\colon\G_1\arr\G_0$ \emph{subhyperbolic} if $\iota$ is contracting on the faces, and there exist structures of finite complexes of groups on $\G_0, \G_1$ for which $f$ and $\iota$ extend, respectively, to a covering and a morphism of complexes of groups. For any such structure and extension (i.e., for any \emph{uniformization} of the substitution) one can associate the iterated monodromy group. The limit dynamical system $f_\infty\colon\G_\infty\arr\G_\infty$ is then topologically conjugate to the \emph{limit dynamical system} of this iterated monodromy group.

Conversely, every \emph{contracting self-similar group} $G$ admits a subhyperbolic polytopal substitution together with a uniformization, for which $G$ is the associated iterated monodromy group~\cite{nek:models}.

These two results were proved in~\cite{nek:models} using proper cocompact group actions on complexes and equivariant maps between them, rather than complexes of groups and their morphisms. The two approaches are equivalent. Group actions are often more convenient theoretically, but they involve infinite objects. Complexes of groups, being finite, are better suited to the description of concrete examples. One of the goals of this paper is therefore to explain the techniques of finite complexes of groups in the setting of polytopal substitutions.

As a consequence of these results, if $f\colon\mathcal{J}\arr\mathcal{J}$ is a locally expanding covering map of a compact metric space, then $f$ is topologically conjugate to the limit dynamical system of a polytopal substitution. Thus a very broad class of fractal dynamical systems can be represented as limits of iterated polytopal substitutions.

It was shown in~\cite{nek:dimension} that the minimal dimension of a subhyperbolic polytopal substitution converging to an expanding covering map $f\colon\mathcal{J}\arr\mathcal{J}$ is equal to the covering dimension of $\mathcal{J}$.
In particular, if $f$ is a post-critically finite complex rational function, then its action on its Julia set is topologically conjugate to the limit dynamical system of its iterated monodromy group and can therefore be represented as the limit of a subhyperbolic polytopal substitution. If the Julia set is not the entire sphere, then it has covering dimension one, and the substitution can consequently be chosen to be one-dimensional.

The one-dimensional case, i.e., hyperbolic virtual endomorphisms of graphs, was studied by D.~Thurston and K.~Pilgrim in~\cite{dthurston:rationalmaps,dthurston:elastic,pigrimthurston:dimension}. The final section of this paper extends their results to higher dimensions and to complexes of groups. In particular, we show how polytopal substitutions can be used to study the conformal dimension of Julia sets of subhyperbolic dynamical systems and of limit spaces of contracting self-similar groups. This part of the paper is a continuation of~\cite{nek:confdim}.

We now state the results of the final section.
Let $f, \iota\colon\G_1\arr\G_0$ be a virtual endomorphism of a connected complex of groups $\G_0$.
Let $B$ be a face of $\G_1$, and let $G_B$ and $G_{f(B)}$ be the \emph{local groups} of the corresponding complexes of groups. The covering $f\colon\G_1\arr\G_0$ of the complexes of groups induces injective homomorphisms $\tilde f\colon G_B\arr G_{f(B)}$. We call the index $[G_{f(B)}:G_B]$ the \emph{local degree} of $f$ at $B$. 

Let $A$ be a face of $\G_0$. List all $f$-preimages $B_1, B_2, \ldots, B_d$ of $A$, with each preimage repeated according to the local degree of $f$ at that preimage. The number $d$ does not depend on $A$ and is equal to the degree of the covering $f$.

Consider the map $F_A\colon A\arr \prod_{i=1}^d\tilde\iota(B_i)$ defined by 
\[F_A(x)=(\iota(f_1^{-1}(x)), \iota(f_2^{-1}(x)), \ldots, \iota(f_d^{-1}(x))),\]
where $f_i$ denotes the restriction of $f$ to $B_i$.
Let $DQ_p(f, \iota)$ denote the supremum, over all faces $A$ of $\G_0$ and all points $x\in A$, of the $\ell^p$-norm of the differential of $F_A$.
We then define $\mathcal{DQ}_p(f, \iota)$ to be the infimum of $DQ_p(f, \iota')$ over all morphisms $\iota'\colon\G_1\arr\G_0$ homotopic to $\iota$.

We obtain the following characterization of the conformal dimension of the limit space.

\begin{theorem}
If $DQ_p(f, \iota)<1$, then the conformal dimension of the limit space of $(f, \iota)$ is at most $\left(\frac 1p-\frac{\log DQ_p(f, \iota)}{\log\deg f}\right)^{-1}$.

Suppose that $DQ_\infty(f, \iota)<1$. Then the asymptotic modulus of the limit space of $(f, \iota)$ is equal to
\[\theta_p=\lim_{n\to\infty}\sqrt[n]{\mathcal{DQ}_p(f^n, \iota^n)}.\]

The conformal dimension of the limit space is equal to
\[\inf\{p\colon \theta_p<1\}.\]
\end{theorem}

The asymptotic modulus $\theta_p$ of the limit space was defined and studied in~\cite{nek:confdim}. The relation between asymptotic modulus and conformal dimension appearing in the last part of the theorem is a result of that paper.

Thus polytopal substitutions provide a method for obtaining upper bounds on the conformal dimension of limit spaces of contracting self-similar groups. Such bounds are of interest not only in dynamics but also in group theory. For example, it was proved in~\cite{mattebonnekzheng:liouvills} that if the limit space of a contracting group has conformal dimension less than $2$, then the group is amenable. Consequently, exhibiting a polytopal substitution with $\mathcal{DQ}_2(f, \iota)<1$ proves amenability of its iterated monodromy group.

The asymptotic modulus $\theta_p$ can be also described in terms of \emph{$p$-energies of graphs}, following~\cite{dthurston:elastic}.
For a continuous map $\phi\colon\Delta_1\arr\Delta_2$ between polytopal complexes, denote by $E_p(\phi)$ the sum of the $p$th powers of the distances in $\Delta_2$ between the $\phi$-images of all pairs of adjacent vertices of $\Delta_1$. Let $f, \iota\colon\G_1\arr\G_0$ be a contracting virtual endomorphism of complexes of groups, and let $\widetilde E_{p, n}$ be the infimum of $E_p(\phi)$ over all morphisms $\phi\colon\G_n\arr\G_0$ homotopic to $\iota^n\colon\G_n\arr\G_0$. 

We show in Theorem~\ref{th:asymptoticDQvirtend} that 
\[\theta_p(f, \iota)=\lim_{n\to\infty} \sqrt[n]{\widetilde E_{p, n}}.\]
The relation between the graph energies $\widetilde E_{p,n}$ and conformal dimension generalizes a result of~\cite{pigrimthurston:dimension}.

The limits $\theta_p$ of the graph energies (especially in the case $p=2$) were studied from the point of view of the theory of Laplacians on fractals in~\cite{kigami:book,HPS:penergy,CaoGuQiu}.

\subsection*{Outline of the paper}

Section 2, \emph{``Polytopal substitutions,''} introduces the basic notions used throughout the paper: topological correspondences, polytopal complexes, and polytopal substitutions. We define polytopal complexes as categories of polytopes and attaching maps between them, without requiring an embedding into an ambient space. We give several examples of polytopal substitutions in this section.

Section 3, \emph{``Self-similar groups,''} introduces iterated monodromy groups of polytopal and topological virtual endomorphisms and views the basic definitions from the theory of contracting self-similar groups and their limit spaces.

Section 4, \emph{``Complexes of groups,''} develops the constructions and notation concerning complexes of groups that are needed in the sequel. We follow~\cite{bridhaefl}, with some modifications. In particular, we define abstract complexes of groups as a special class of small categories. The usual definition is recovered by choosing representatives of isomorphism classes of objects and representatives of the orbits of automorphism groups on sets of morphisms. Since these choices are not canonical, the standard formulation requires the introduction of “twisting elements” in the definitions of complexes of groups and their morphisms. The general category-theoretic formulation is more natural theoretically, whereas the concrete formulation based on chosen representatives is often more convenient for explicit examples. We therefore describe both approaches.

Section 5, \emph{``Polytopal models of self-similar groups,''} introduces subhyperbolic polytopal substitutions and reviews the results of~\cite{nek:models} and their applications to virtual endomorphisms of complexes of groups.

Finally, Section 6, \emph{``Polytopal substitutions and conformal dimension,''} characterizes the conformal dimension of the limit space in terms of the asymptotic $\ell^p$-contraction coefficients $\mathcal{DQ}_p(f, \iota)$ and graph energies $\tilde E_{p, n}$ of the corresponding polytopal substitutions.

\subsection{Acknowledgments}  The author is grateful to Kevin Pilgrim, Dylan Thurston, and Misha Hlushchanka for discussions on the topics of the paper and useful remarks. The author was supported by NSF grants DMS-2204379 and DMS-2604376.

\subsection{AI declaration} AI was used only for proofreading and editing some of the figures.

\section{Polytopal substitutions}

\subsection{Topological correspondences}

\begin{defi}
    A \emph{topological correspondence} (or a \emph{multivalued dynamical system}) is a pair of continuous maps $f, \iota\colon\M_1\arr\M_0$ between topological spaces.
\end{defi}

Topological correspondences are generalizations of topological dynamical systems, since they can be iterated. We interpret a topological correspondence as a graph with the space of vertices $\M_0$, set of edges $\M_1$, and the source and range maps $f, \iota$, respectively.

The space of length $n$ paths this graph is the subspace $\M_n\subseteq\M_1^n$ consisting of sequences $(t_1, t_2, \ldots, t_n)$ satisfying \[f(t_i)=\iota(t_{i-1})\] for all $i=2, \ldots, n$. We regard $\M_n$ as the space of \emph{orbits of length $n$} of the correspondence. (Note that we write the order of the arrows from right to left.) 

The \emph{$n$th iteration} of the correspondence is then the pair of maps $f^n, \iota^n\colon\M_n\arr\M_0$ given by
\[f^n(t_1, t_2, \ldots, t_n)=f(t_1),\qquad \iota^n(t_1, t_2, \ldots, t_n)=\iota(t_n).\]
We define the maps $f_n\colon\M_{n+1}\arr\M_n$ and $\iota_n\colon\M_{n+1}\arr\M_n$  by
\[f_n(t_1, t_2, \ldots, t_n, t_{n+1})=(t_1, \ldots, t_n),\qquad\iota_n(t_1, t_2, \ldots, t_{n+1})=(t_2, \ldots, t_{n+1}),\]
and $f_0=f, \iota_0=\iota$. Then $f^n=f_0\circ f_1\circ\cdots\circ f_{n-1}$ and $\iota^n=\iota_0\circ\iota_1\circ\cdots\iota_{n-1}$. The maps $f_n$ and $\iota_n$ are also given by the pullback diagram
\[\begin{array}{ccc}
\M_{n+1} & \stackrel{\iota_n}{\arr} & \M_n\\
\mapdown{f_n} & & \mapdown{f_{n-1}}\\
\M_n & \stackrel{\iota_{n-1}}{\arr} & \M_{n-1}\end{array}\]

\begin{example}
If $f\colon\M_1\arr\M$ is a continuous map, where $\M_1$ is a subset of $\M$, then we can consider it as a correspondence by setting $\iota\colon\M_1\arr\M$ to be the identity embedding. Then iterations of the partial map $f$ coincide with the iterations of the correspondence, so that $\M_n$ is the domain of $f^n$, and $f_n\colon\M_{n+1}\arr\M_n$ is the restriction of $f$ to $\M_{n+1}$. This correspondence is a virtual endomorphism if $f$ is a covering map. In this case, we call $f$ a \emph{partial self-covering}.
\end{example}

The \emph{limit space} $\M_\infty$ of a correspondence $f, \iota\colon\M_1\arr\M$ is the inverse limit of the spaces $\M_n$ with respect to the maps $\iota_n$. It is homeomorphic to the subspace of the Cartesian product $\M_1^{-\omega}$ consisting of sequences $(\ldots, t_2, t_1)$ such that $f(t_n)=\iota(t_{n+1})$ for all $n$. In other words, it is the space of all \emph{forward orbits} of the correspondence. The limit of the maps $f_n$ is the map $f_\infty(\ldots, t_2, t_1)=(\ldots, t_3, t_2)$.

\begin{example}
    In the setting of partial self-maps $f\colon\M_1\arr\M$ for $\M_1\arr\M$, the limit space $\M_\infty$ is naturally identified with the intersection of the domains $\M_n$ of the iterations $f^n$ of the partial self-covering. 
\end{example}

\subsection{Polytopal complexes}

Our definition of polytopal complexes is based on~\cite[Chapter~III.C]{bridhaefl}.

If $\mathcal{C}$ is a small category, then, borrowing notation from theory of groupoids, we denote by $\mathcal{C}^{(0)}$ and $\mathcal{C}^{(1)}$ its sets of objects and morphisms, respectively. If $\varphi\in\mathcal{C}^{(1)}$, then we denote by $\be(\varphi)$ and $\en(\varphi)$ the source and range of $\varphi$, respectively, so that $\varphi$ is a morphism from $\be(\varphi)$ to $\en(\varphi)$. We compose morphisms using left notation: a composition $\varphi_1\circ\varphi_2$ is defined if and only if $\en(\varphi_2)=\be(\varphi_1)$. We denote the set of sequences $(\varphi_1, \varphi_2, \ldots, \varphi_n)$ of morphisms such that $\varphi_1\circ\varphi_2\circ\cdots\circ\varphi_n$ is defined by $\mathcal{C}^{(n)}$.

We often consider small categories up to equivalence, i.e., replacing a category $\mathcal{C}$ by a full sub-category $\mathcal{C}_1$ intersecting every isomorphism class of $\mathcal{C}$.

\begin{defi}
    A \emph{scwol} (\emph{small category without loops}) is a small category in which all endomorphisms are identity automorphisms.
\end{defi}

For example, a small category in which for every pair of objects $\alpha, \beta$, the set of morphisms from $\alpha$ to $\beta$ has at most one element, i.e., the category associated with a partially ordered set, is a scwol. In this sense scwols are generalizations of posets.

\begin{defi}
A \emph{polytopal complex} $\Delta$ of dimension $n$ consists of the following data.

\begin{enumerate}
\item
For each $0\leq k\leq n$, a set $\mathcal{P}_k$ of $k$-dimensional
convex polytopes in $\R^k$.

\item
For every $A\in\mathcal{P}_k$ and every $l$-dimensional face
$B\subseteq A$, where $l<k$, an isometry from $B$ onto an element
$B'\in\mathcal{P}_l$. We regard this isometry as a partially defined
map $\varphi_{A,B}\colon A\dashrightarrow B'$
with domain $B$, and call it an \emph{attaching map}. In addition,
for every $A\in\mathcal{P}_k$, the identity map $1_A\colon A\longrightarrow A$ is an attaching map.
\end{enumerate}

We require the collection of attaching maps to be closed under
composition. Thus, with objects $\bigcup_{k=0}^n\mathcal{P}_k$ 
and morphisms given by the attaching maps, it forms a category, which
we call the \emph{scwol} of $\Delta$ and denote by
$\mathcal{C}(\Delta)$.

The \emph{support} $|\Delta|$ of $\Delta$ is the quotient of the
disjoint union $\bigsqcup_{k=0}^n\bigsqcup_{A\in\mathcal{P}_k}A$
by the identifications induced by the attaching maps.
\end{defi}

Elements of the sets $\mathcal{P}_i$ are called \emph{faces} of the complex. To avoid ambiguity, a face of a polytope $A\in\mathcal{P}_i$ is called in this case a \emph{subface} of $A$.

The composition of the morphisms $\varphi_{A, B}\colon A\dashrightarrow B'$ and $\varphi_{B', C}\colon B'\dashrightarrow C'$ is the morphism $\varphi_{A, \varphi_{A, B}^{-1}(C)}$. In particular, sequences of composable attaching maps of length $n$  are in a bijection with flags $A\supseteq A_1\supseteq A_2\supseteq \cdots\supseteq A_n$ of subfaces of faces $A\in\mathcal{P}_k$ of $\Delta$. The set of all such flags is in a natural way the set of faces of the \emph{barycentric subdivision} $\Br(\Delta)$ of $|\Delta|$. 

The support of a connected polytopal complex is naturally a length metric space with the length of curves induced by the Euclidean structure on the faces.

\begin{example}
A 1-dimensional polytopal complex is a graph (with loops and multiple edges allowed) together with a choice of lengths of its edges.
\end{example}

\begin{example}
    Let $\Delta$ be the classical representation of the torus as the quotient of the square $[0, 1]\times [0,1]\subset\R^2$ with identified opposite sides. It is a polytopal complex with one two-dimensional face (the square), two edges (representing the sides $\{0\}\times [0, 1]$ and $[0, 1]\times\{0\}$), and one vertex. Hence, there are four objects in the associated scwol $\mathcal{C}(\Delta)$.
    
There are twelve non-identity morhpisms in $\mathcal{C}(\Delta)$: two attaching maps from the two-dimensional face to each of the edges, two attaching maps from each edge to the vertex, and four attaching maps from the two-dimensional face to the vertex.
\end{example}

\begin{example}
Every simplicial complex can be turned into a polytopal complex by identifying each $k$-dimensional simplex with the standard simplex \[\left\{(x_0, \ldots, x_k)\colon 0\le x_i\le 1, \sum_{i=0}^kx_i=1\right\}\subset\R^{k+1}\] with the natural attaching maps. Note that in this case the associated scwol is the category associated with the partially ordered set of simplices, i.e., the set of morphisms between any two objects has at most one element (which is not true for general polytopal complexes, see the previous example).
\end{example}

Not every scwol is the scwol of a polytopal complex. For example, a scwol with the set $\{A, B\}$ of two objects and more than two morphisms from $A$ to $B$ is not associated with a polytopal complex. 

Conversely, every scwol $\mathcal{C}$ has a natural \emph{geometric realization}, which is a simplicial complex. Its vertices are objects of $\mathcal{C}$, and a simplex corresponds to a composable sequence $(\varphi_n, \ldots, \varphi_2, \varphi_1)\in\mathcal{C}^{(n)}$ of non-identity morphisms.

If $\mathcal{C}=\mathcal{C}(\Delta)$ is the scwol of a polytopal complex $\Delta$, then its geometric realization is naturally isomorphic to the \emph{barycentric subdivision} $\Br(\Delta)$ of $\Delta$. In particular, its support is naturally homeomorphic to $|\Delta|$.

The objects of $\mathcal{C}(\Delta)$ are identified with the barycenters of the corresponding faces, the attaching map from a face $A$ to its subface $B$ is identified with the edge of $\Br(\Delta)$ from the barycenter of $A$ to the barycenter of $B$ inside $A$. The \emph{canonical orientation} of an edge of $\Br(\Delta)$ is towards the cell of lower dimension, so that it agrees with the direction of the corresponding morphisms of $\mathcal{C}(\Delta)$. We will often identify the scwol associated with a complex $\Delta$ with the one-skeleton of its barycentric subdivision $\Br(\Delta)$. 

\begin{defi}
\label{def:morphismofpcomplexes}
A \emph{morphism} $f\colon\G_1\arr\G_0$ between two polytopal complexes is given by a functor of the associated scwols $\tilde f\colon\mathcal{C}(\G_1)\arr\mathcal{C}(\G_0)$, and a collection of continuous maps $f\colon A\arr \tilde f(A)$ between the corresponding affine polyhedra, agreeing with the attaching maps in the following sense. If $\varphi\colon A\dashrightarrow B$ is an attaching map in $\G_1$, and $\tilde f(\varphi)\colon\tilde f(A)\dashrightarrow\tilde f(B)$ is its image under the functor $\tilde f$, then we have the equality $f\circ \varphi=\tilde f(\varphi)\circ f$ of the corresponding continuous maps on the domain of $\varphi$. We say that the morphism $f$ is \emph{defined over the functor $\tilde f$}.

We say that the morphism is \emph{affine} or \emph{smooth} if the maps $f\colon A\arr \tilde f(A)$ are affine or smooth, respectively.
\end{defi}

A morphism $f\colon\G_1\arr\G_0$ naturally induces a continuous map $|\G_1|\arr |\G_0|$, which we will also denote by $f$. If we identify $\mathcal{C}(\G_i)$ with the barycentric subdivisions $\Br(\G_i)$, then the functor $\tilde f$ will induce a map $\tilde f\colon |\Br(\G_1)|\arr |\Br(\G_0)|$ homotopic to $f$.

We stress that the functor $\tilde f$ is a part of the definition of the morphism. One can have two morphisms with different functors $\tilde f_1, \tilde f_2$, but with the same induced map on the supports of the complexes.

\subsection{Polytopal substitutions}
\label{ss:polsubstitutions}

\begin{defi}
\label{def:labeling}
Let $\G$ be a polytopal complex. A \emph{labeling} of a polytopal complex $\Delta$ by $\G$ is a morphism $f\colon\Delta\arr\G$ such that for every face $A$ of $\Delta$ the map $f\colon A\arr \tilde f(A)$ is a surjective isometry.
\end{defi}

Note that the condition of Definition~\ref{def:labeling}  implies that for every face $A$ of $\Delta$ the map
\begin{equation}
    \label{eq:branchedcovercondtion}
\tilde f\colon \be^{-1}(A)\arr\be^{-1}(f(A))
\end{equation}
is a bijection.

\begin{defi}
\label{def:polysubstitution}
A \emph{polytopal substitution} is a pair of morphisms $f, \iota\colon\G_1\arr\G_0$ of polytopal complexes such that $f$ is a surjective labeling of $\G_1$ by $\G_0$. 
\end{defi}

Let $f, \iota\colon\G_1\arr\G_0$ be a polytopal substitution, and let $\Delta$ by a polytopal complex with a labeling $F\colon\Delta\arr\G_0$. Then the \emph{result of application} of the substitution $\sigma=(f, \iota)$ to the labeled complex $\Delta$ is the $\G_0$-labeled complex $\sigma(\Delta)$ obtained by closing the following pullback diagram of polytopal complexes:
\[\begin{array}{ccc}
\sigma(\Delta) & \arr & \Delta\\
\mapdown{F'} & & \mapdown{f}\\
\G_1  & \stackrel{\iota}{\arr} & \G_0\end{array}\]
where the labeling is given by the composition $f\circ F'\colon\sigma(\Delta)\arr\G_0$.

Let us describe the complex $\sigma(\Delta)$ in more explicit terms. The faces of $\sigma(\Delta)$, i.e., the objects of the associated scwol $\mathcal{C}(\sigma(\Delta))$ are in a bijection with the set of pairs of faces $(A_1, A_2)\in\mathcal{C}^{(0)}(\G_1)\times\mathcal{C}^{(0)}(\Delta)$ of $\G_1$ and $\Delta$, such that $\tilde\iota(A_1)=\tilde f(A_2)$. The polytope associated with the pair $(A_1, A_2)$ is isometric to $\tilde f(A_2)$. The set of morphisms of $\mathcal{C}(\sigma(\Delta))$ is also equal to the set of pairs $(\varphi_1, \varphi_2)\in\mathcal{C}^{(1)}(\G_1)\times\mathcal{C}^{(1)}(\Delta)$ such that $\tilde\iota(\varphi_1)=\tilde f(\varphi_2)$, with coordinate-wise composition.

In other words, the labeled complex $\sigma(\Delta)$ is obtained from $\Delta$ by replacing each face $A$ of $\Delta$ by the subcomplex $\sigma(A)=(\tilde\iota)^{-1}(\tilde f(A))$ of $\G_1$ equal to the preimage under $\iota$ of the label of $A$. 
Since $\sigma(\Delta)$ is also $\G_0$-labeled, we can apply the substitution $\sigma=(f, \iota)$ iteratively, and get the complexes $\sigma^n(\Delta)$. 

In particular, we can iteratively apply the substitution $\sigma=(f, \iota)$ to $\G_0$ with the identity labeling, and obtain complexes $\G_n=\sigma^n(\G_0)$. Note that $\sigma(\G_0)$ is naturally isomorphic to $\G_1$. 

The objects and morphisms of $\mathcal{C}(\G_n)$, are identified with sequences $(\alpha_1, \alpha_2, \ldots, \alpha_n)$ of objects and morphisms, respectively, of $\mathcal{C}(\G_1)$ such that $\tilde f(\alpha_i)=\tilde\iota(\alpha_{i-1})$ for all $i=2, \ldots, n$. Compare this with the definition of the spaces $\M_n$ of the length $n$ orbits of a topological correspondence.

We denote, same as for topological correspondences,
\[\tilde f_n(\alpha_1, \alpha_2, \ldots, \alpha_n, \alpha_{n+1})=(\alpha_1, \ldots, \alpha_n),\qquad \tilde\iota_n(\alpha_1, \alpha_2, \ldots, \alpha_n, \alpha_{n+1})=(\alpha_2, \ldots, \alpha_{n+1})\]
and 
\[\tilde{f^n}(\alpha_1, \alpha_2, \ldots, \alpha_n)=\tilde f(\alpha_1),\qquad \tilde{\iota^n}(\alpha_1, \alpha_2, \ldots, \alpha_n)=\tilde\iota(\alpha_n).\]
A face $(A_1, A_2, \ldots, A_n)$ is isometric to the face $\tilde f(A_1)$ of $\G_0$.

It follows directly from the definitions that the supports $|\G_n|$ and the maps $f_n, \iota_n\colon|\G_{n+1}|\arr|\G_n|$ will be naturally homeomorphic (respectively, conjugate) to the spaces and maps obtained by iterating the topological correspondence $f, \iota\colon |\G_1|\arr |\G_0|$.

Let $(A_1, A_2, \ldots, A_n)\in\mathcal{C}(\G_n)^{(0)}$, and let $\varphi_1\in\mathcal{C}(\G_1)^{(1)}$ be such that $\be(\varphi_1)=A_1$. Then $\tilde f(A_2)=\tilde\iota(A_1)$, so $\be(\iota(\varphi_1))=\iota(A_1)=\tilde f(A_2)$. By Definition~\ref{def:labeling}, there is a unique $\varphi_2\in\mathcal{C}(\G_1)^{(1)}$ such that $\be(\varphi_2)=A_2$ and $\tilde f(\varphi_2)=\tilde\iota(\varphi_1)$. Continuing this argument, we conclude that for every $A=(A_1, A_2, \ldots, A_n)\in\mathcal{C}(\G_n)^{(1)}$ and every $\varphi_1\in\be^{-1}(A_1)$ there is a unique $\varphi=(\varphi_1, \varphi_2, \ldots, \varphi_n)\in\mathcal{C}(\G_n)^{(1)}$ such that  $\be(\varphi)=A$. Thus, every $\varphi=(\varphi_1, \varphi_2, \ldots, \varphi_n)$ is uniquely determined by its source $(\be(\varphi_1), \be(\varphi_2), \ldots, \be(\varphi_n))$ and $\varphi_1$.

The process of finding $(\varphi_1, \varphi_2, \ldots, \varphi_n)$ from $(A_1, A_2, \ldots, A_n)$ can be interpreted as the operation of an automaton (a \emph{transducer}) in the following way. 

It is the automaton with the set of internal states $\mathcal{C}(\G_0)^{(1)}$ accepting and giving on output strings of elements of $\mathcal{C}(\G_1)^{(0)}$. 
Automaton being in the state $\varphi\in\mathcal{C}(\G_0)^{(1)}$ accepts only letters $A\in\mathcal{C}(\G_1)^{(0)}$ such that $\iota(A)=\be(\varphi)$. If it accepts it, then its next state  is the unique element $\varphi_1\in\mathcal{C}(\G_1)^{(1)}$ such that $\be(\varphi_1)=A$ and $\tilde\iota(\varphi_1)=\varphi$. The automaton outputs the letter $\en(\varphi_1)$, and its next state will be $\iota(\varphi_1)$.

The automaton being in the state $\varphi_1$ accepts a string $(A_1, A_2, \ldots, A_n)$ if and only if $(A_1, A_2, \ldots, A_n)\in\mathcal{C}(\G_n)^{(0)}$ and $\iota(A_1)=\be(\varphi_1)$. The sequence of the successive states of the automaton, when it is processing the string, is the unique sequence $(\varphi_1, \varphi_2, \ldots, \varphi_n)\in\mathcal{C}(\G_n)^{(1)}$, and its output is the string $(\en(\varphi_1), \en(\varphi_2), \ldots, \en(\varphi_n))\in\mathcal{C}(\G_n)^{(0)}$. 

We denote the output by
\[S_{\varphi_1}(A_1, A_2, \ldots, A_n)=(\en(\varphi_1), \en(\varphi_2), \ldots, \en(\varphi_n)).\]

It follows that the polytopal complex $\G_n$ can be described in the following way.  For every $(A_1, A_2, \ldots, A_n)\in\mathcal{C}(\G_n)^{(0)}$ we take a copy of the polytope isometric to $A_1$. If $\varphi_1\colon A_1\dashrightarrow B_1$ is an attaching map in $\G_1$ (it is an isometry of a subface $A_1'$ of $A_1$ with $B_1$), then we attach the subface $A_1'$ of the copy $(A_1, A_2, \ldots, A_n)$ of $A_1$ to the copy $S_{\varphi_1}(A_1, A_2, \ldots, A_n)$ of $B_1$ by $\varphi_1$. By functoriality, this will be a well-defined polytopal complex.

\begin{example}
    Let $\G_0$ be the graph consisting of one vertex $v$ and one loop $e$. The corresponding scwol has two objects (corresponding to the loop and the vertex) and two non-identity morphisms $\varphi_0, \varphi_1$ from $e$ to $v$. Let $\G_1$ be the graph with two vertices $v_0, v_1$ and two edges $e_0$, $e_1$, connecting them. We have four attaching maps $\varphi_{i, j}$ from $e_i$ to $v_j$. We assume that all edges have length 1.
    
    Let $f\colon\G_1\arr\G_0$ be the covering map acting by isometries of the edges and on the scwol by the functor \[\tilde f(\varphi_{0, 1})=\varphi_1, \tilde f(\varphi_{1, 0})=\varphi_1, \tilde f(\varphi_{0, 0})=\varphi_0, \tilde f(\varphi_{1, 1})=\varphi_0.\] 
    
Let $\iota\colon\G_1\arr\G_0$ be the homeomorphism of the corresponding circles over the functor $\tilde\iota$ that acts by
\[\tilde\iota(v_1)=\tilde\iota(e_0)=\tilde\iota(\varphi_{0, 1})=\tilde\iota(\varphi_{1, 1})=v,\]
where we identify objects with their identity automorphisms, and
\[\tilde\iota(v_0)=v,\quad\tilde\iota(\varphi_{0, 0})=\varphi_0,\quad\tilde\iota(\varphi_{1, 0})=\varphi_1.\]

The corresponding substitution replaces an edge by a path of two edges, so the complex $\G_n$ is a cycle of $2^n$ edges.
\end{example}

We will be less formal in the subsequent examples.

\begin{example}
\label{ex:wordsubstitution}
Let $\G_0$ be the graph with one vertex and a set $E$ of loops. For every loop $e\in E$ choose a path $\sigma(e)$ in $\G_0$. Let $\G_1$ be the graph consisting of $|E|$ cycles $W_e$, labeled by the loops $e\in E$ with one common vertex $v_0$ (only $v_0$ is common for any two cycles $W_e$). Let $\iota\colon\G_1\arr\G_0$ be the homeomorphism mapping homeomorphically every cycle $W_e$ to the loop $e$. Let $f\colon\G_1\arr\G_0$ map each cycle $W_e$ to the path $\sigma(e)$. If we orient the loops of $\G_0$ and assume that the paths $\sigma(e)$ follow the orientation, then the polytopal substitution $f, \iota\colon\G_1\arr\G_0$ models a homomorphism $\sigma\colon E^*\arr E^*$ of the free monoid. 

If $\Delta$ is a graph labeled by $\G_0$, i.e., a directed graph whose edges are labeled by elements of $E$, then $\sigma(\Delta)$ is obtained by replacing each edge $\alpha$ of $\Delta$ by a chain of edges such that the labels on the chain spell $\sigma(e)$, where $e$ is the label of $\alpha$.

In particular, the graph $\G_n$ is a bouquet of loops mapped by $f^n$ to the paths corresponding to the $n$th iteration $\sigma^n(e)$ of the \emph{word substitution} $\sigma\colon E^*\arr E^*$.
\end{example}

\begin{figure}
    \centering
    \includegraphics{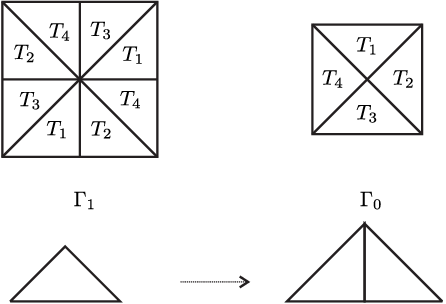}
    \caption{A double covering of the torus}
    \label{fig:torus}
\end{figure}

\begin{example}
\label{ex:doubletorusmap}
Consider the torus obtained by identifying the opposite sides of $[0, 1]\times [0, 1]$, and let us additionally subdivide it into four triangles by its diagonals. Let $\G_0$ be the obtained polytopal complex. Let $\G_1$ be the torus obtained from $\G_0$ by subdividing the triangles of $\G_0$ by the bisectors of their right angles, see Figure~\ref{fig:torus}. Let $\iota\colon\G_1\arr\G_0$ be the natural homeomorphism (the identity map, if we consider $\G_1$ to be a subdivision of $\G_0$). Let us multiply all lengths in $\G_1$ by $\sqrt{2}$, and let $f\colon\G_1\arr\G_0$ be the covering map acting isometrically on the faces of $\G_1$ and mapping triangles according to the labels in Figure~\ref{fig:torus} by orientation-preserving isometries. Then $\G_n$ is obtained from $\G_0$ by $n$ successive subdivisions of triangles by bisectors of their right angles.
\end{example}

In all previous examples, the map $\iota\colon|\G_1|\arr|\G_0|$ was a homeomorphism. Such polytopal substitutions are called \emph{subdivision rules}. We can identify then $|\G_1|$ with $|\G_0|$ by $\iota$, thus transforming $f\colon|\G_1|\arr |\G_0|$ into a dynamical system $f\colon|\G_0|\arr|\G_0|$ on a topological space. The complex $\G_{n+1}$ is then a subdivision of the complex $\G_0$. This can be used to define dynamical systems and to study their properties. See~\cite{cfp:subdivision} for applications of this approach to complex dynamics and hyperbolic geometry. 

Other examples of subdivision rules are self-similar tilings. For instance, the subdivision rule producing Penrose tilings can be formulated as a polytopal substitution as follows.

\begin{example}
\label{ex:penrose}
Consider the polytopal complex $\G_0$ obtained by gluing two triangles (one acute with angles $72^\circ, 72^\circ, 36^\circ$ and one obtuse with angles $36^\circ, 36^\circ, 108^\circ$) along their sides to form a trapezoid shown on the right-hand side of Figure~\ref{fig:match}, and then identifying the lateral sides of the trapezoid by the isometry preserving the orientation shown in the figure. We will get then a cylindrical complex shown on the bottom of Figure~\ref{fig:match}. Let $\G_1$ be the subdivision of $\G_0$ shown on the left-hand side of the figure. 

\begin{figure}
\centering
\includegraphics{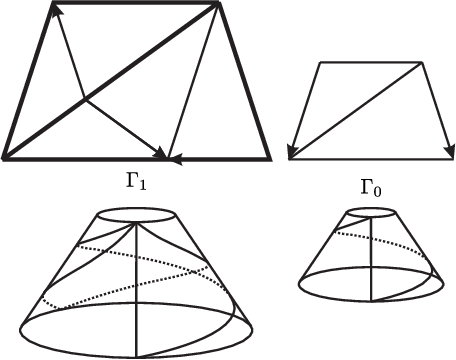}
\caption{Penrose tiling substitution}
\label{fig:match}
\end{figure}

As before, $\iota$ is the homeomorphism (similarity). We have highlighted the preimages of the edges of $\G_0$ in $\G_1$.
The map $f$ acts on the faces by isometries preserving the orientation of the oriented side of the triangles.
The iteration of the substitution produces patches of Penrose tilings. The complex $\G_5$, for example, is shown in Figure~\ref{fig:penrose5}.

\begin{figure}
\centering
\includegraphics{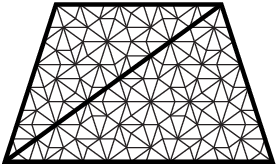}
\caption{The complex $\G_5$ of the Penrose tiling subdivision}
\label{fig:penrose5}
\end{figure}

If we apply the described substitution to an arbitrary Penrose tiling $\mathcal{P}$, seen as a $\G_0$-labeled complex, then we get a subdivision $\sigma(\mathcal{P})$ of the tiling $\mathcal{P}$. The tiling $\mathcal{P}$ is the \emph{inflaction} of the tiling $\sigma(\mathcal{P})$.
\end{example}

Let us consider some examples that are not subdivision rules, i.e.,  when $\iota$ is not a homeomorphism of the underlying topological spaces.

\begin{figure}
\centering
\includegraphics{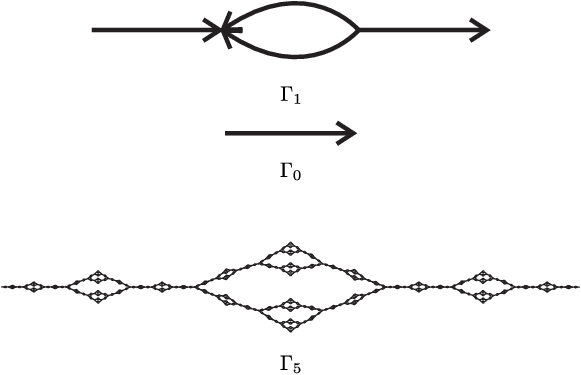}
\caption{A ``bubble'' substitution}
\label{fig:gromov}
\end{figure}

\begin{example}
\label{ex:gromovsubst}
The following example is from~\cite{gromov:cell90}, adapted to our approach. Let $\G_0$ and $\G_1$ be as in Figure~\ref{fig:gromov}. The map $f\colon\G_1\arr\G_0$ is the orientation-preserving map isometric on the edges. The map $\iota$ collapses the two central edges to a point, and maps the other two arrows by orientation preserving maps to halves of the segment $\G_0$. The graphs $\G_n$ are obtained by iteratively replacing the middle third of all edges by a double edge.
\end{example}

\begin{example}
\label{ex:sierptriangle}
Consider the graph $\G_0$ shown in Figure~\ref{fig:sierptriangle}, and the substitution shown below it. Vertices labeled by $A, B, C$ are mapped by the substitution to one vertex with the same label. Vertices labeled by $O$ are mapped to a complex consisting of three disjoint vertices labeled by $A, B, C$. The substitution on edges is shown in the figure.

\begin{figure}
    \includegraphics{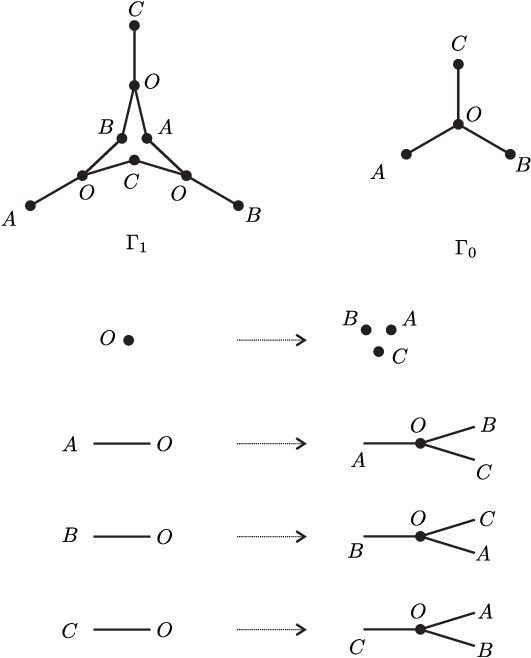}
    \caption{Sierpi\'nski gasket substitution}
    \label{fig:sierptriangle}
\end{figure}

\begin{figure}
\includegraphics{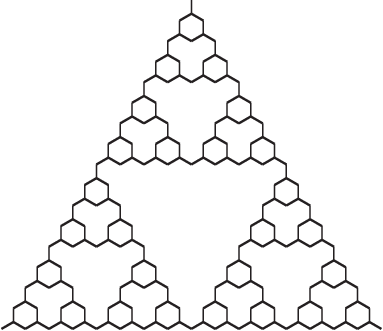}
\caption{Sierpi\'nski gasket approximation}
\label{fig:striangle}
\end{figure}

The result of applying this substitution to $\G_0$ once is the graph $\G_1$ shown in Figure~\ref{fig:sierptriangle}. The map $f\colon\G_1\arr\G_0$ maps the vertices of $\G_1$ to the vertices of $\G_0$ labeled by the same letters. It maps an edge with endpoints $X, Y$ of $\G_1$ isometrically to the edge labeled by $X, Y$ of $\G_0$. The map $\iota$ is inverse to the substitution, so it maps smaller tripods of $\G_1$ attached to $A$, $B$, $C$ to the edges $AO, BO, CO$, by mapping the central vertices labeled by $A, B, C$ to $O$ and the vertices labeled by $O$ to the middles of the edges $AO, BO, CO$ of $\G_0$. See Figure~\ref{fig:striangle} for the graph $\G_4$ obtained by iterating this substitution.
\end{example}

\begin{example}
\label{ex:dendrites}
The following substitution appears in~\cite{nek:dendrites}. It is associated with the map 
\begin{equation}
\label{eq:dendrites}
f(z_1, z_2)=\left(\left(1-\frac{2z_1}{z_2}\right)^2, \left(1-\frac{2}{z_2}\right)^2\right).
\end{equation}

Let $\G_0$ be the simplicial complex obtained by gluing three tetrahedra $A_1ABC$, $B_1ABC$, $C_1ABC$ along their common face $ABC$. 

Let $\G_1$ be obtained from a square pyramid $CA_1A_2A_4A_3$ by attaching two tetrahedra $C_1'CA_1A_4$ and $C_2'CA_2A_3$ along diagonal cross-sections of the pyramid. Subdivide the square pyramid by viewing it as a cone over the barycentric subdivision of its base. Each attached tetrahedron is subdivided into two simplices by dividing the corresponding diagonal cross-section into two triangles along the height $BC$ of the pyramid. See Figure~\ref{fig:tetrahedral}.

\begin{figure}
\centering
\includegraphics{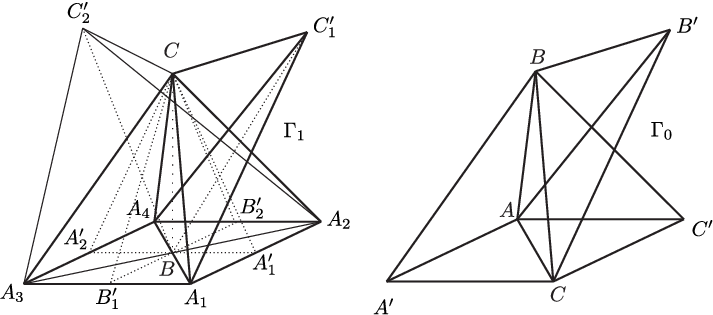}
\caption{``Hubbard complex'' for a 2D complex map}
\label{fig:tetrahedral}
\end{figure}

Define the map $f\colon\G_1\arr\G_0$ to be the simplicial map defined by the rules $f(X_i)=X$ and $f(X_i')=X'$. The map $\iota\colon\G_1\arr\G_0$ is affine on the simplices of $\G_1$, maps the 
vertices as follows:
\[\iota\colon A_3\mapsto A',\quad A_1\mapsto C, \quad A_4\mapsto A, \quad A_2\mapsto C',\quad C_1'\mapsto B', \quad C_2'\mapsto B.\]
In particular, it maps $B$, $B_1'$, $B_2'$, $A_1'$, $A_2'$ to the middles of the segments $AC$, $A'C$, $AC'$, $CC'$, $AA'$, respectively.
\end{example}

\section{Self-similar groups}
\subsection{Polytopal virtual endomorphisms and iterated monodromy groups}

Let $f, \iota\colon\G_1\arr\G_0$ be a polytopal substitution. We are sometimes using $f$ instead of $\tilde f$, since $\tilde f(A)=f(A)$, so there is no ambiguity.

As we have seen in Subsection~\ref{ss:polsubstitutions}, every morphism $\varphi\in\mathcal{C}_0^{(1)}$ defines a map
\[S_\varphi\colon f^{-n}(\be(\varphi))\arr f^{-n}(\en(\varphi)).\]

For every face $A$ of $\G_0$ the disjoint union $T_A=\bigsqcup_{n\ge 0}f^{-n}(A)$ of fibers of the map $f^n\colon\G_n\arr\G_0$ has a natural structure of a rooted tree in which $(A, A_2, \ldots, A_n)\in T_n$ is connected to $(A, A_2, \ldots, A_{n-1})$. Here $f^{-0}(A)=\{A\}$ is the root of $T_A$. 

For every $\varphi\in\mathcal{C}(\G_0)^{(1)}$, the map $S_\varphi\colon T_{\be(\varphi)}\arr T_{\en(\varphi)}$ is a morphism of trees, i.e., it maps the root to the root, preserves the levels and adjacency of the vertices.

\begin{defi}
    A topological correspondence $f, \iota\colon\M_1\arr\M_0$ such that $f$ is a finite degree covering map is called a \emph{virtual endomorphism} of $\M_0$.
    
    A polytopal substitution $f, \iota\colon\G_1\arr\G_0$ is called a \emph{(polytopal) virtual endomorphism} (of $\G_0$) if the morphism $f$ is a covering of finite degree.
\end{defi}

Here the degree of the covering is the number of $f$-preimages.
We assume throughout this paper that $\G_0$ is connected, so that the degree is constant.

Note that if $f, \iota\colon\G_1\arr\G_0$ is a polytopal substitution, then $f$ is a covering if and only if for every face $A$ of $\G_1$ the map 
\[\tilde f\colon\en^{-1}(A)\arr\en^{-1}(\tilde f(A))\]
is a bijection. Recall that we already have that $\tilde f\colon\be^{-1}(A)\arr\be^{-1}(\tilde f(A))$ is a bijection, according to the definition of a polytopal substitution.

It is easy to prove then by induction that if the substitution $f, \iota\colon\G_1\arr\G_0$ is a virtual endomorphism, then the maps $f^n\colon\G_n\arr\G_0$ are coverings and the maps $S_\varphi\colon f^{-n}(\be(\varphi))\arr f^{-n}(\en(\varphi))$ are bijections.

Consequently, the maps $S_\varphi\colon T_{\be(\varphi)}\arr T_{\en(\varphi)}$ are isomorphisms of rooted trees.

Let $e$ be an oriented edge of the barycentric subdivision $\Br(\G_0)$ and let $\varphi\in\mathcal{C}(\G_0)^{(1)}$ be the corresponding morphism of the scwol. If the orientation of $e$ is canonical (towards lower dimension), then we set $S_e=S_\varphi$. Otherwise, we set $S_e=S_\varphi^{-1}$. 

 If $\gamma$ is a path in the one-skeleton of the barycentric subdivision $\Br(\G_0)$ starting in a vertex $A$ and ending in a vertex $B$, then we get the associated composition of the maps $S_e$ along the path. It is a bijection from $f^{-n}(A)$ to $f^{-n}(B)$, coinciding with the classical \emph{monodromy action} of the fundamental groupoid of $\Br(\G_0)$ on the fibers of the covering map $f^n\colon\G_n\arr\G_0$.

  In particular, we get an action of the fundamental group $\pi_1(\G_0, A)$ on the tree $T_A$ for every face $A$ of $\G_0$ (i.e, for every vertex of $\Br(\G_0)$), which we call the \emph{iterated monodromy action}. The quotient of $\pi_1(\G_0, A)$ by the kernel of the iterated monodromy action is called the \emph{iterated monodromy group} of $(f, \iota)$ and is denoted $\img(f, \iota)$. It does not depend, up to conjugacy, of the actions on rooted trees on the choice of the basepoint $A$, if $\G_0$ is connected.
 
 More generally, if $f, \iota\colon\M_1\arr\M_0$ is a topological correspondence, where $\M_0$ is path-connected, then the fundamental group $\pi_1(\M_0, t)$ acts naturally by the monodromy action on the tree of preimages $T_t=\bigsqcup_{n\ge 0}f^{-n}(t)$ of the basepoint under the covering maps $f^n\colon\M_n\arr\M_0$. In the remainder of this sub-section we give a review of the main notions of the theory of iterated monodromy groups of topological virtual endomorphisms, keeping in mind that virtual endomorphisms of polytopal complexes are a special case of topological virtual endomorphisms. For more details, see~\cite{nek:models} and~\cite{nek:dyngroups}.

Let $f, \iota\colon\M_1\arr\M_0$ be a topological virtual endomorphism, where $\M_0$ is path-connected. Choose a basepoint $t\in\M_0$, and consider the corresponding iterated monodromy group $\img(f, \iota)$.

Let $t_1\in f^{-1}(t)$. Then the map
\[(t_2, t_3, \ldots, t_n)\mapsto (t_1, t_2, \ldots, t_n)\]
from $f^{-(n-1)}(\iota(t_1))$ to $f^{-n}(t)$ gives a natural isomorphic identification $T_{\tilde\iota(t_1)}\arr T_t$ of rooted trees. 

If $\gamma$ is a path in $\M_0$ starting in $t$ and ending in $\iota(t_1)$, then the map $S_\gamma\colon T_t\arr T_{\iota(t_1)}$ is an isomorphism of rooted trees, and so the above identification interprets it as an isomorphic embedding $S_\gamma\colon T_t\arr T_t$ of $T_t$ to its sub-tree. Let us denote the set of such isomorphisms by $\bim$. It is easy to check that the set $\bim$ is invariant under pre- and post-compositions with elements of $\img(f, \iota)$. We get two commuting actions (one left and one right) of the iterated monodromy group on $\bim$. We call $\bim$ the \emph{biset} of the iterated monodromy group. The (right) action by pre-compositions is free (i.e., $x\cdot g=x$ implies that $g$ is the identity) and has $\deg f$ orbits labeled by the elements of $f^{-1}(t)$.

\begin{defi}
A \emph{self-similar group} is a pair $(G, \bim)$, where $G$ is a group, and $\bim$ is a $G$-biset, i.e., a set with commuting left and right actions of $G$ on it, such that the right action of $G$ is free. We assume throughout the paper that the number of right $G$-orbits is finite. 

We say that two self-similar groups $(G_1, \bim_1)$ and $(G_2, \bim_2)$ are \emph{equivalent} if there exists an isomorphism $\phi\colon G_1\arr G_2$ and a bijection $F\colon\bim_1\arr\bim_2$ such that $F(g_1\cdot x\cdot g_2)=\phi(g_1)\cdot F(x)\cdot \phi(g_2)$ for all $x\in\bim$ and $g_1, g_2\in G$.
\end{defi}

We always consider the iterated monodromy group $\img(f, \iota)$ as the self-similar group with the biset $\bim$, described above.

It is also convenient to encode the iterated monodromy action of the fundamental group $\pi_1(\M_0)$ in a self-similar group, which is done in the following way.

\begin{defi}
Denote by $\bim_{f, \iota, t}=\bim_{f, \iota}$ the set of pairs $(z, [\ell])$ where $z\in f^{-1}(t)$, and $[\ell]$ is the homotopy class of a path in $\M_0$ starting in $t$ and ending in $\iota(z)$. Define the left and the right actions of $\pi_1(\M_0, t)$ on $\bim_{f, \iota}$ in the following way. For $[\gamma]\in\pi_1(\M_0, t)$, we set
\[(z, [\ell])\cdot [\gamma]=(z, [\ell\gamma]),\qquad [\gamma]\cdot (z, [\ell])=(S_\gamma(z), [\iota(\gamma_z)\ell]),\]
where $\gamma_z$ is the lift of $\gamma$ by $f$ to a path starting in $z$, and $S_\gamma(z)$ is the end of $\gamma_z$. We compose paths as functions: in a concatenation $\beta\alpha$, the path $\alpha$ is passed before $\beta$.
\end{defi}

The following proposition is straightforward.

\begin{proposition}
    Let $K\lhd \pi_1(\M)$ be the kernel of the iterated monodromy action. Let $\bim_{f, \iota}/K$ be the set of the right $K$-orbits on $\bim$. Then $(\pi_1(\M)/K, \bim/K)$  is equivalent to the iterated monodromy group $\img(f, \iota)$.
\end{proposition}

If $\bim_1, \bim_2$ are two $G$-bisets, then $\bim_1\otimes\bim_2$ is defined as the quotient of $\bim_1\times\bim_2$ by the equivalence relation $x_1\otimes x_2=x_1\cdot g\otimes g^{-1}\cdot x_2$, for $g\in G$ and $x_i\in\bim_i$. It is a $G$-biset with respect to the actions $g_1\cdot (x_1\otimes x_2)\cdot g_2=(g_1\cdot x_1)\otimes (x_2\cdot g_2)$. The bisets $(\bim_1\otimes\bim_2)\otimes\bim_3$ and $\bim_1\otimes (\bim_2\otimes\bim_3)$ are isomorphic, so we have well defined $G$-bisets $\bim^{\otimes n}$ for all $n\ge 1$. We also define $\bim^{\otimes 0}=G$ with the natural left and right actions (by multiplication).

We thus obtain the biset $\bim^*=\bigsqcup_{n\ge 0}\bim^{\otimes n}$. It is also a semigroup with respect to the operation $\otimes$. The biset structure agrees with multiplication by the elements of the subgroup $G=\bim^{\otimes 0}$ of the semigroup $\bim^*$.

In the case of the iterated monodromy group $(\img(f, \iota), \bim)$, the semigroup $\bim^*$ is naturally isomorphic to the semigroup of isomorphisms $T_t\hookrightarrow T_t$ of the tree of preimages $T_t$ with the trees of preimages $T_z$ for $z\in T_t$. The multiplication $\otimes$ in $\bim^*$ corresponds to compositions of the isomorphisms.

In the case of the self-similar group $(\pi_1(\M_0, t), \bim_{f, \iota})$, the bisets $\bim^{\otimes n}_{f, \iota}$ are naturally isomorphic to the bisets $\bim_{f^n, \iota^n}$. For $(z_1, [\ell_1])\in\bim_{f^{n_1}, \iota^{n_1}}$ and $(z_2, [\ell_2])\in\bim_{f^{n_2}, \iota_{n_2}}$, the product $(z_1, [\ell_1])\otimes (z_2, [\ell_2])$ is the element \[(z_3, [\iota^{n_2}(\ell_{1, z_2})\ell_2])\in\bim_{f^{n_1+n_2}, \iota^{n_1+n_2}},\] where $\ell_{1, z_2}$ is the lift of $\ell_1$ by $f^{n_2}$ to a path starting in $z_2$, and $z_3$ is the end of this lift.

A subset $\alb\subset\bim$ is a call a \emph{basis} of the biset if it intersects every orbit exactly once, i.e., if every element of $\bim$ is uniquely written as $x\cdot g$ for $x\in\alb$ and $g\in G$ (uniqueness follows from the fact that the right action is free).

Let $\alb$ be a basis of $\bim$. Then for every $x\in\alb$ and $g\in G$ there exist unique $y\in\alb$ and $h\in G$ such that $g\cdot x=y\cdot h$. The obtained map $x\mapsto y$ is a permutation of $\alb$ (corresponding to the left action of $G$ on $\bim/G$). We denote, therefore, $y=g(x)$. We also denote $h=g|_x$, and call it \emph{section} of $g$ in $x$.

If $\alb$ is a basis of $\bim$, then $\alb^{\otimes n}$ is a basis of $\bim^{\otimes n}$. We denote $x_1\otimes x_2\otimes \cdots\otimes x_n\in\alb^{\otimes n}$ just by $x_1x_2\ldots x_n$, and identify $\alb^{\otimes n}$ with the direct power $\alb^n$. We denote $\xs=\bigcup_{n\ge 0}\alb^n$, where $\alb^0$ consists of a single \emph{empty word} $\varnothing$.

Same as for $\alb$, we get an action of $G$ on $\alb^n$ and the corresponding notion of the section, so that for every $v\in\alb^n$ and $g\in G$ we have
\[g\cdot v=g(v)\cdot g|_v.\]
We distinguish the action $v\mapsto g\cdot v$ of $G$ on $\bim^{\otimes n}$ and the action $v\mapsto g(v)$ of $G$ on $\alb^{\otimes n}$.

For small bases $\alb$, it is convenient to encode the action of an element $g\in G$ on $\alb$ and the sections $g|_x$ for $x\in\alb$ into one map
\[g\mapsto \sigma(g|_x)_{x\in\alb},\]
where $\sigma$ is the permutation induced by $g$ on $\alb$. We interpret $\sigma(g|_x)_{x\in\alb}$ as an element of the semi-direct product $\mathrm{Sym}_{\alb}\ltimes G^\alb$ of the symmetric group on $\alb$ with the direct power $G^\alb$ (the \emph{permutational wreath product} of $\mathrm{Sym}_\alb$ with $G$). We get a homomorphism $G\arr\mathrm{Sym}_\alb\ltimes G^\alb$, which we call the \emph{wreath recursion}. 

Let us show some examples of the computation of iterated monodromy groups of polytopal substitutions.

\begin{example}
\label{ex:basilicasubst}
Consider the substitution shown in Figure~\ref{fig:basilicam}. The arrows on the edges are shown only to specify the actions of the maps $f, \iota$ as orientation-preserving. The morphism $f$ is the label and orientation preserving isometry on the edges. The morphism $\tilde\iota$ maps the vertex $1$, the loop labeled by $b$ adjacent to it, and both edges labeled by $a$ to the edge $b$ of $\G_0$. It maps the rest of the graph (i.e., the loop at $0$ labeled by $b$) to the edge $a$ of $\G_0$. It maps the vertex $0$ of $\G_1$ to the unique vertex of $\G_0$. The maps $\iota$ map the vertex $1$ and the loop labeled by $b$ adjacent to it to the middle of the loop $b$ of $\G_0$, and maps the arrows labeled by $a$ in $\G_0$ to two halves of the loop $b$. The loop labeled by $b$ adjacent to the vertex $0$ is mapped by $\iota$ homeomorphically to the loop $b$.

The corresponding substitution $\sigma$ mapping edges of $\G_0$ to their preimages in $\G_1$ is shown on the right-hand side of the figure. Figure~\ref{fig:basilicagraph2} shows the corresponding graph $\G_6$.

\begin{figure}
    \centering
    \includegraphics{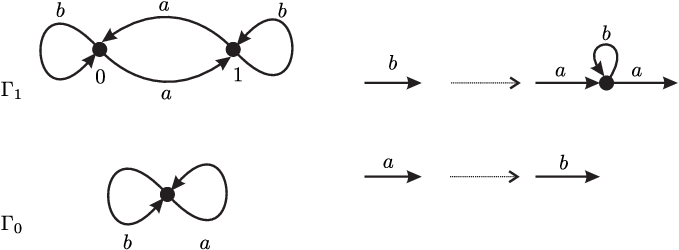}
    \caption{Basilica substitution}
    \label{fig:basilicam}
\end{figure}

Let us compute the iterated monodromy action of $\pi_1(\G_0)$ (and thus the iterated monodromy group as its faithful quotient).
We take the loops $a$ and $b$ as the generators of the fundamental group. Choose the vertex of $\G_0$ as the basepoint $t$. We interpret $\tilde\iota$ as the map between the barycentric subdivisions (identifying the scwols with them). Then $\tilde\iota(f^{-1}(t))$ consists of $t$ and the midpoint of the loop $b$. Let us choose the connecting paths $\ell_0, \ell_1$ to be the trivial path at $t$ for the former and the path from $t$ to the midpoint of $t$ going against the orientation of the arrow in Figure~\ref{fig:basilicam} (i.e., the lower half).

The loop $a$ is lifted by $f$ to the two edges of $\G_1$ connecting the vertices $0$ and $1$. Their images under $\iota$ are halves of the loop $b$ of $\G_0$. As we assume that $f$ and $\tilde\iota$ preserve the orientation, we get that, in the biset $\bim_{f, \iota}$, we have 
\[a\cdot(0, [\ell_0])=(1, [\ell_1])\cdot 1,\qquad a\cdot(1, [\ell_1])=(0, [\ell_0])\cdot b).\]
We also have
\[b\cdot (0, [\ell_0])=(0, [\ell_0])\cdot 1,\qquad b\cdot (1, [\ell_1])=(1, [\ell_1])\cdot a.\]

\begin{figure}
    \centering
    \includegraphics{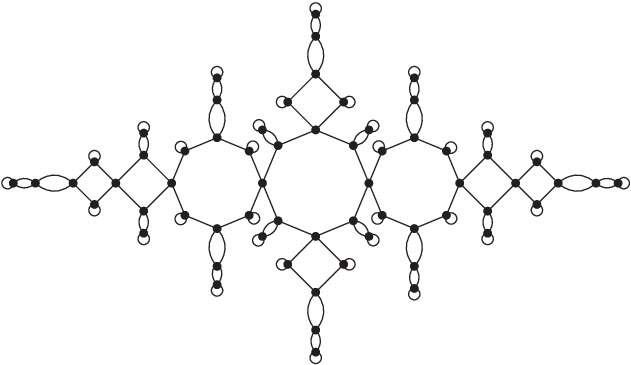}
    \caption{The graph $\G_6$ approximating the Julia set of $z^2-1$}
    \label{fig:basilicagraph2}
\end{figure}

If we use the symbols $0=(0, [\ell_0])$, $1=(1, [\ell_1])$, we write this as the wreath recursion
\[a\mapsto (01)(1, b),\qquad b\mapsto  (1, a).\]

The group generated by this wreath recursion is the iterated monodromy group of the complex polynomial $z^2-1$ (i.e., of its restriction to the Julia set), and is called the \emph{basilica group}, see~\cite{zukgrigorchuk:3st} and its computation in~\cite{nek:book}.
\end{example}

\begin{example}
Let us compute the iterated monodromy group of the double-covering of the torus, described in Example~\ref{ex:doubletorusmap}.
Let us take the common vertex of the right angles of the triangles $T_i$ of $\G_0$ as the basepoint $t$. Then $f^{-1}(t)$ consists of the midpoints of the sides of the square shown on the right-hand side of Figure~\ref{fig:torus}. They represent two points of the torus. Their images under $\tilde\iota$ are the midpoints of the hypotenuses of the triangles $T_i$ in $\G_0$. Let us connect them to the basepoint by the paths $\ell_0, \ell_1$ equal to the bisectors in the triangles $T_3$ and $T_4$, respectively.

\begin{figure}
    \centering
    \includegraphics{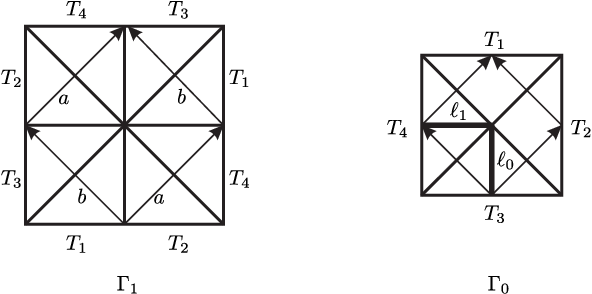}
    \caption{Iterated monodromy group of a torus map}
    \label{fig:imgtorus}
\end{figure}

Let us take, as generators of the fundamental group $\pi_1(\G_0, t)$, the simple horizontal loop $a$ and the simple vertical loop $b$. We orient them left-to-right and bottom-to-top, respectively. The curves $f^{-1}(a), f^{-1}(b)$ and their $\iota$-images are shown as arrows in Figure~\ref{fig:imgtorus}. The lifts of $a$ and $b$ are marked by the respective letters in $\G_1$. It follows that the wreath recursion defining the iterated monodromy group is
\[a=(01)(a, b),\qquad b=(01)(1, a^{-1}b).\]
\end{example}

\begin{example}
Consider the virtual endomorphism shown on the top part of Figure~\ref{fig:thincarpet}. The labels show which vertices and edges of $\G_1$ are mapped by $f$ to which vertices and edges of $\G_0$. Namely, the covering $f\colon\G_1\arr\G_0$ is the unique orientation and label preserving map from $\G_1$ to $\G_0$. 

\begin{figure}
    \centering
    \includegraphics{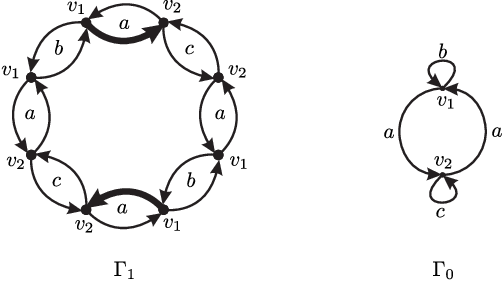}
    \caption{A substitution modeling the Julia set of $z^2-\frac{1}{16z^2}$}
    \label{fig:thincarpet}
\end{figure}

The map $\iota\colon\G_1\arr\G_0$ maps the thicker arrows on the top and on the bottom of $\G_1$ to the vertices $v_1$ and $v_2$ of $\G_0$, respectively. Arrows labeled by $a$ adjacent to them are mapped to the loops $b$ and $c$ of $\G_0$, respectively. The left and the right chains of three cycles of length 2 are mapped by $\iota$ to the arrows of $\G_1$ labeled by $a$: the left one to the arrow from $v_1$ to $v_2$, and the right one to the arrow from $v_2$ to $v_1$. The corresponding substitution of vertices and edges is shown in Figure~\ref{fig:carpetsub}. Figure~\ref{fig:thincarpetgr} shows the graph $\G_4$.

\begin{figure}
    \centering
    \includegraphics{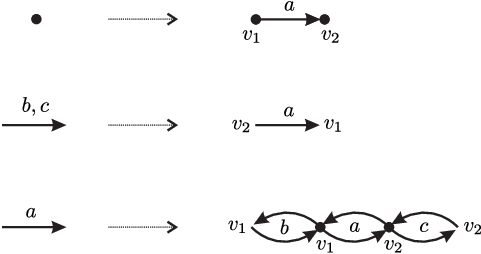}
    \caption{The substitution}
    \label{fig:carpetsub}
\end{figure}

\begin{figure}
    \centering
    \includegraphics{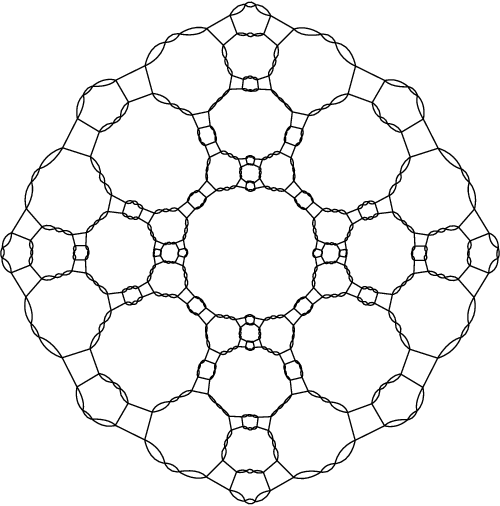}
    \caption{An approximation of the Julia set of $z^2-\frac{1}{16z^2}$}
    \label{fig:thincarpetgr}
\end{figure}

Let us compute $\img{(f, \iota)}$. Choose the basepoint $t=v_1$ in $\G_0$. It has four preimages in $\G_1$ (vertices labeled by $v_1$). They are mapped by $\iota$ to $v_1, v_2$, and some points in the interiors of the edges of $\G_0$ labeled by $a$. Let us connect $v_1$ to the points of $\iota(f^{-1}(v_1))$ by paths in $\G_0$ inside the arrows labeled by $a$ in the direction agreeing with the orientation shown on the figure (i.e., from the left side). The generators $\alpha, \beta, \gamma$ of the fundamental group $\pi_1(\G_0, v_1)$ are as follows. The loop $\alpha$ is the loop formed by two arrows of $\G_0$ labeled by $a$, but in the \emph{opposite} direction to the one shown on the figure (i.e., clockwise). The generator $\beta$ is the loop labeled by $b$, and the generator $\gamma$ is the arrow labeled by $a$ (starting in $v_1$ and ending in $v_2$) followed by the loop $c$, and then followed by the inverse of the arrow labeled by $a$. We will get the wreath recursion

\begin{align*}
\alpha &= (\beta, 1, \gamma, 1),\\
\beta &= (14)(23),\\
\gamma &= (12)(34)(\alpha, \alpha^{-1}, 1, 1).
\end{align*}
It coincides with the wreath recursion for the iterated monodromy group of $z^2-\frac{1}{16z^2}$ (see~\cite{mattebonnekzheng:liouvills}).
\end{example}

\begin{example}
\label{ex:torusbasilica}
Let $f, \iota\colon\G_1\arr\G_0$ be the virtual endomorphism constructed in Example~\ref{ex:doubletorusmap}. Denote by $O_0$ and $O_1$ the vertices of the acute angles of the triangular faces of $\G_1$, where $O_0$ is the one mapped by $\iota$ to the vertex of acute angles of the faces of $\G_0$, and $O_1$ is the one mapped to the vertex of the right angles.  Let $\Delta_0$ be obtained by taking two copies $\G_{0, 0}$ and $\G_{0, 1}$ of $\G_0$ and pasting them together by identifying their vertices of acute angles. Let $\Delta_1$ be obtained by taking a copy of $\G_1$, and two copies of $\G_0$, and attaching the vertices of acute angles of the copies of $\G_0$ to the vertices $O_0$ and $O_1$ of $\G_1$, see Figure~\ref{fig:torusbasilica}. The covering $f_1\colon\Delta_1\arr\Delta_0$ maps $\G_1$ to $\G_{0, 0}$ by $f$, and maps the copies of $\G_0$ in $\Delta_1$ to $\G_{0, 1}$ isomorphically, so that we get a degree 2 covering $f_1\colon\Delta_1\arr\Delta_0$. The morphism $\iota_1\colon\Delta_1\arr\Delta_0$ maps the copy of $\G_0$ attached to $O_1$ to the vertex of $\G_{0, 1}$ corresponding to $\iota(O_1)$ (i.e., to the vertex of the right angles of $\G_{0, 1}$), maps the other copy of $\G_0$ to $\G_{0, 0}$ isomorphically, and the copy of $\G_1$ to $\G_{0, 1}$ by $\iota$.

    \begin{figure}
        \includegraphics{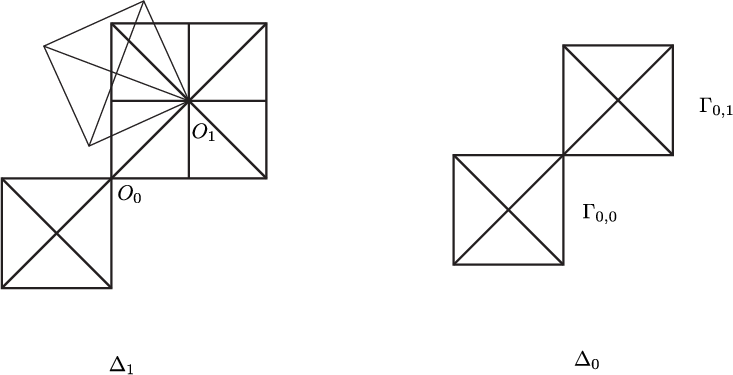}
        \caption{Basilica operation applied to a double-covering of the torus}
        \label{fig:torusbasilica}
    \end{figure}

    The iterated monodromy group of $f_1, \iota_1\colon\Delta_1\arr\Delta_0$ is  generated by the wreath recursion
    \begin{alignat*}{2}
    a_0 &= (01)(a_1, b_1), &\qquad b_0 &= (01)(1, a_1^{-1}b_1),\\
    a_1 &= (1, a_0), &\qquad b_1&= (1, b_0).
    \end{alignat*}

This group is the result of application of the ``Basilica operation'' (defined and studied in~\cite{petschickrajeev})  to the iterated monodromy group of the torus self-covering from Example~\ref{ex:doubletorusmap}.
\end{example}

\subsection{Contracting groups and their limit spaces}

\begin{defi}
We say that a virtual endomorphism $f, \iota\colon\M_1\arr\M_0$ of compact metric spaces is \emph{contracting} (or \emph{hyperbolic}) if $f$ is a local isometry and $\iota$ is contracting, i.e., if there exist constants $\delta>0$ and $L>1$ such that
\[d(f(\xi), f(\zeta))=d(\xi, \zeta)\]
for all $\xi, \zeta\in\M_1$ such that $d(\xi, \zeta)<\delta$, and \[d(\iota(\xi), \iota(\zeta))\le L^{-1}d(\xi, \zeta)\]
for all $\xi, \zeta\in\M_1$.
\end{defi}

In particular, a virtual endomorphism $f, \iota\colon\G_1\arr\G_0$ of finite polytopal complexes, then it is contracting if and only if the maps $\iota\colon A\arr\tilde{\iota^n}(A)$ are contracting for all faces $A$ of $\G_1$.

\begin{defi}
    Let $(G, \bim)$ be a self-similar group, and let $\alb\subset\bim$ be a basis. We say that $(G, \bim)$ is \emph{contracting} if there exists a finite set $\nuke\subset G$ such that for every $g\in G$ and there exists $n$ such that $g|_v\in\nuke$ for all $v\in\xs$ of length $\ge n$.

    The smallest set $\nuke$ with this property is called the \emph{nucleus} of the action.
\end{defi}

It is shown in~\cite[Corollary~2.11.7]{nek:book} that the property of being contracting does not depend on the choice of the basis (though the nucleus does depend on it).

For a proof of the following proposition, see~\cite[Theorem~5.10]{nek:models} and~\cite[Theorem~4.5.36]{nek:dyngroups}.  It is essentially the same as the proof of~\cite[Theorem~5.5.3]{nek:book}.

\begin{proposition}
Suppose that $f, \iota\colon\M_1\arr\M_0$ is a virtual endomorphism such that $f^n, \iota^n\colon\M_n\arr\M_0$ is contracting for some $n\ge 1$. Then the self-similar groups $(\pi_1(\M_0), \bim_{f, \iota})$ and the iterated monodromy group $(\img{(f, \iota)}, \bim)$ are contracting.
\end{proposition}

Let $(G, \bim)$ be a contracting self-similar group, and let $\alb\subset\bim$ be a basis. Consider the space $\xmo$ of left-infinite sequences of elements of $\alb$ with the direct product topology. We say that $\ldots x_2x_1, \ldots y_2y_1\in\xmo$ are \emph{asymptotically equivalent} if there exists a finite set $N\subset G$ and a sequence $g_k\in N$ such that 
\[g_n(x_n\cdots x_2x_1)=y_n\cdots y_2y_1\]
for every $n$.

We say that two sequences $\ldots x_2x_1\cdot g, \ldots y_2y_1\cdot h\in\xmo\times G$ are \emph{asymptotically equivalent} if there exists a finite set $N\subset G$ and a sequence $g_n\in N$ such that 
\[g_n\cdot x_n\ldots x_2x_1\cdot g=y_n\ldots y_2y_1\cdot h\]
in $\bim^{\otimes n}$ for every $n$.

The quotients of the spaces $\xmo$ and $\xmo\times G$ by the asymptotic equivalence relations are called the \emph{limit space $\lims$} and \emph{limit $G$-space} $\limg$, respectively.

The action of $G$ on $\xmo\times G$ by right multiplication induces a right action of $G$ on $\limg$ by homeomorphisms. This action does not depend on the choice of the basis $\alb$, up to topological conjugacy.

The shift $\ldots x_2x_1\mapsto \ldots x_3x_2$ induces a well-defined continuous map  $f\colon\lims\arr\lims$. We call it the \emph{limit dynamical system} of $(G, \bim)$. It does not depend, up to topological conjugacy on the choice of the basis $\alb$.

The following is proved in~\cite[Proposition~4.5.23]{nek:dyngroups}.

\begin{theorem}
    The limit dynamical system is a covering map if and only if the group $(G, \bim)$ is \emph{regular}, i.e., if for every $g\in G$ and every $x_1x_2\ldots \in\xo$ there exists $n$ such that $g(x_1x_2\ldots x_n)\ne x_1x_2\ldots x_n$ or $g|_{x_1x_2\ldots x_n}=1$.

    The iterated monodromy group of a locally expanding covering map is a regular self-similar group.
\end{theorem}

The next theorem is a partial case of~\cite[Theorem~5.10]{nek:models}, see also~\cite[Theorem~4.5.36]{nek:dyngroups}.

\begin{theorem}
\label{th:mainimgregular}
Let $f, \iota\colon\M_1\arr\M_0$ be a hyperbolic virtual endomorphism. Then the limit dynamical system $f_\infty\colon\M_\infty\arr\M_\infty$ is topologically conjugate to the limit dynamical system of the iterated monodromy group $\img{(f, \iota)}$.

In particular, if $f\colon\M\arr\M$ is a locally expanding covering, then it is topologically conjugate to the limit dynamical system of $\img{f}$.
\end{theorem}

It was shown in~\cite[Theorem~6.7]{nek:models} that for every regular contracting group $(G, \bim)$ there exists a polytopal hyperbolic virtual endomorphism $f, \iota\colon\G_1\arr\G_0$ such that $\img{(f, \iota)}$ is equivalent to $(G, \bim^{\otimes n})$ for some $n$. We will describe this in Subsection~\ref{ss:gequivariantmaps} below. In particular, this gives a method of approximating arbitrary spaces $\M$ admitting a locally expanding covering map $f\colon\M\arr\M$ by substitutional polytopal complexes.

\section{Complexes of groups}

We would like to extend Theorem~\ref{th:mainimgregular} to all contracting groups (not only the regular ones), so that we get approximations of their limit spaces by substitutional polytopal complexes. This will also help us to include more examples of dynamical systems (for example, subhyperbolic complex rational functions).

\subsection{The definition}

We give here an overview of the theory of complexes of groups. We follow~\cite{bridhaefl} with some modifications.

\begin{defi}
\label{def:complexofgroups}
An \emph{(abstract) complex of groups} is a small category $\mathcal{CG}$ satisfying the following conditions:
\begin{enumerate}
    \item Every endomorphism $e\in\mathcal{CG}^{(1)}$ is an automorphism.
    \item For any two objects $\alpha, \beta\in\mathcal{CG}^{(0)}$ the automorphism groups $G_\alpha=\mathop{\mathrm{Aut}}(\alpha)$ and $G_\beta=\mathop{\mathrm{Aut}}(\beta)$ act freely (by compositions) on the set of morphisms from $\alpha$ to $\beta$, and the $G_\alpha$-orbits are contained in the $G_\beta$-orbits.
\end{enumerate}
\end{defi}

We consider small categories up to equivalence of categories, i.e., passing to a full sub-category spanned by a set of objects intersecting every isomorphism class. 

\begin{defi}
    \label{def:polyhedralcomplexofgroupsgeneral}
    A \emph{polytopal complex of groups} $\Delta$ consists of an abstract complex of groups $\mathcal{CG}(\Delta)$, an assignment of convex Euclidean polytopes to the objects of $\mathcal{CG}(\Delta)$ and a functorial assignment to every morphism $\varphi\in\mathcal{CG}(\Delta)^{(1)}$ of an \emph{attaching map} $\varphi$, which is an isometry from a subface of $\be(\varphi)$ to the face $\en(\varphi)$. The attaching maps associated with automorphisms act by identity isometries. 
\end{defi}

Let $\mathcal{CG}$ be a complex of groups. By condition (2) of Definition~\ref{def:complexofgroups}, for every morphism $\varphi\in\mathcal{CG}^{(1)}$ and every element $g\in G_{\be(\varphi)}$ there exists an element $h\in G_{\en(\varphi)}$ such that $\varphi\cdot g=h\cdot\varphi$. Since the actions of the automorphism groups on the sets of morphisms are free, the element $h$ is unique. The map $g\mapsto h$ is an injective homomorphism of groups, which we will denote $\Psi_\varphi\colon G_{\be(\varphi)}\arr G_{\en(\varphi)}$.

Denote by $\mathcal{C}$ the category with the same set of objects as $\mathcal{CG}$ and with morphisms equal to the orbits $G_{\en(\varphi)}\cdot\varphi=G_{\en(\varphi)}\cdot\varphi\cdot G_{\be(\varphi)}$ of morphisms $\varphi\in\mathcal{CG}^{(1)}$ under the action of the automorphism groups of their sources and ranges. For every composable pair $(\varphi_1, \varphi_2)\in\mathcal{CG}^{(2)}$ and any automorphisms $g_i\in G_{\en(\varphi_i)}$, we have $g_1\cdot\varphi_1\cdot g_2\cdot \varphi_2=g_1\Psi_{\varphi_1}(g_2)\cdot e_1e_2$. Consequently, the orbit of $\varphi_1\varphi_2$ depends only on the orbits of $\varphi_1$ and $\varphi_2$, hence we get a well-defined category. The category $\mathcal{C}$ is a scwol. 

If $\Delta$ is a polytopal complex of groups, then the attaching maps associated with the morphisms of $\mathcal{CG}(\Delta)$ are invariant under the actions of the automorphism groups. Consequently, the polytopal complex defined by the complex of groups is the same as the one defined by the associated scwol. 

Natural examples of complexes of groups are the complexes associated with actions of groups on polytopal complexes.
A polytopal complex $\Delta$ together with an action of a group $G$ on it by automorphisms is called a (polytopal) \emph{$G$-complex}. We always assume that every element of $G$ that leaves a face invariant acts on it by identity. In particular, this condition ensures that the quotient $\Delta/G$ is also a polytopal complex in the natural way.

We say that the $G$-complex is \emph{proper}, resp.\ \emph{cocompact}, if the corresponding action on the support is proper, resp.\ cocompact. An action on a locally finite complex is proper if and only if the stabilizer of every face is finite. An action on a finite-dimensional polytopal complex is cocompact if and only if the set of orbits of the action on faces is finite.

\begin{defi}
\label{def:complexofaction}
Let $\Delta$ be a right $G$-complex. The complex of groups $\Delta/G$ is defined by the category with the set of objects equal to the set of faces of $\Delta$. A morphism from $A$ to $B$ is a pair $(g, \varphi)$, where $g\in G$ and $\varphi\colon A\dashrightarrow B\cdot g$ is an attaching map from $A$ to $B\cdot g$ (i.e., an element of the scwol of $\Delta$). We include pairs $(g, Id_A)$, which will be isomorphisms of the category $\mathcal{CG}(\Delta/G)$. 

If $\varphi_1\colon A\dashrightarrow B\cdot g_1$, and $\varphi_2\colon B\dashrightarrow C\cdot g_2$ are attaching maps in $\Delta$, then $\varphi_2\cdot g_1\colon B\cdot g_1\dashrightarrow C\cdot g_2g_1$ is also an attaching map, so we can define composition of elements of $\mathcal{CG}(\Delta/G)$ by
\[(g_2, \varphi_2)\circ(g_1, \varphi_1)=(g_2g_1, (\varphi_2\cdot g_1)\circ \varphi_1),\]
where the composition on the right-hand side of the equality is in the scwol of $\Delta$, where $G$ acts naturally on the set of attaching maps. It is checked directly that we get in this way a complex of groups.
\end{defi}

If $\Delta$ is a $G$-complex, then we denote by $\Delta/G$ the associated polytopal complex of groups, and by $|\Delta/G|$ its underlying complex (i.e., complex of orbits of the action of $G$ on $\Delta$).

\begin{example}
Let $\Delta$ be a $G$-complex. Suppose that $\mathcal{T}$ is a set of faces of $\Delta$ intersecting every $G$-orbit exactly once.
Then the category $\mathcal{CG}(\Delta/G)$ is equivalent to its full subcategory with the set of objects $\mathcal{T}$. If the $G$-complex $\Delta$ is proper and cocompact, then this subcategory is finite. The scwol associated with this sub-category is isomorphic to the scwol of the quotient polytopal complex $|\Delta/G|$. 
\end{example}

Let $\mathcal{C}$ be the scwol of a complex of groups $\mathcal{CG}$. Suppose that every isomorphism class of $\mathcal{CG}$ is a singleton (like in the above case of the restriction to a $G$-orbit transversal $\mathcal{T}$).
We label every object $\alpha\in\mathcal{C}^{(0)}$ by its automorphism group $G_\alpha$ in $\mathcal{CG}$.
Choose for every morphism $\varphi\in\mathcal{C}^{(1)}$ a representative $\overline\varphi\in\mathcal{CG}^{(1)}$ of the corresponding orbit. It will define the associated monomorphism $\Psi_\varphi=\Psi_{\overline\varphi}\colon G_{\be(\varphi)}\arr G_{\en(\varphi)}$.

For every composable pair $(\varphi_1, \varphi_2)\in\mathcal{G}^{(2)}$ there is a unique  \emph{twisting element} $t_{\varphi_1, \varphi_2}\in G_{\en(\varphi_1)}$ such that
\[\overline{\varphi_1}\cdot\overline{\varphi_2}=t_{\varphi_1, \varphi_2}\cdot\overline{\varphi_1\varphi_2}.\]

It is checked directly that the homomorphism $\Psi_\varphi$ and the elements $t_{\varphi_1, \varphi_2}$ satisfy the identity
\begin{equation}
    \label{eq:comgroup1}
\Psi_{\varphi_1}\circ\Psi_{\varphi_2}(g)=t_{\varphi_1, \varphi_2}\cdot\Psi_{\varphi_1\varphi_2}(g)\cdot t_{\varphi_1, \varphi_2}^{-1}
\end{equation}
for all $(\varphi_1, \varphi_2)\in\mathcal{C}^{(2)}$ and every element $g\in G_{\en(\varphi_2)}$, and the identity
\begin{equation}
    \label{eq:compgroup2}
    \Psi_{\varphi_1}(t_{\varphi_2, \varphi_3})t_{\varphi_1, \varphi_2\varphi_3}=t_{\varphi_1, \varphi_2}t_{\varphi_1\varphi_2, \varphi_3}
\end{equation}
for every composable triple $(\varphi_1, \varphi_2, \varphi_3)\in\mathcal{C}^{(3)}$.

The category $\mathcal{CG}$ is then reconstructed from the above data in the following way. Its set of objects is equal to the set of objects of $\mathcal{C}$. The set of morphisms from $\alpha$ to $\beta$ is $\mathop{Mor}_{\mathcal{CG}}(\alpha, \beta)=G_\beta\times\mathop{Mor}_{\mathcal{G}}(\alpha, \beta)$. The morphisms are composed by the rule
\begin{equation}
    \label{eq:compositionrule}
(g_1, \varphi_1)(g_2, \varphi_2)=(g_1\Psi_{\varphi_1}(g_2)t_{\varphi_1, \varphi_2}, \varphi_1\varphi_2).
\end{equation}

If we make a different choice $r_\varphi\cdot\overline\varphi$ of the representatives of orbits of morphisms, for some $r_\varphi\in G_{\en(\varphi)}$, then the homomorphisms $\Psi_\varphi$ will be changed to 
\[\Psi_\varphi'(g)=r_\varphi\cdot \Psi_\varphi(g)\cdot r_\varphi^{-1},\]
and the twisting element $t_{\varphi_1, \varphi_2}$ will be changed to \[t_{\varphi_1, \varphi_2}'=r_{\varphi_1}\Psi_{\varphi_1}(r_{\varphi_2})t_{\varphi_1, \varphi_2}r_{\varphi_1\varphi_2}^{-1}.\]

This leads to the following version of the definition of complexes of groups, as in~\cite{bridhaefl}. We are using here an interpretation of elements of the scwol as directed edges of the barycentric subdivision $\Br(\Delta)$. We direct the edges of $\Br(\Delta)$ towards faces of lower dimension and call them \emph{arrows}.

\begin{defi}
    \label{def:complexofgroupstriple}
    Let $\Delta$ be a polytopal complex.
A \emph{polytopal complex of groups} with the \emph{underlying complex} $\Delta$ is defined by the following data:
\begin{enumerate}
    \item a group (called \emph{local group}) $G_A$ for every face $A$ of $\Gamma$, i.e., every vertex of $\Br(\Delta)$;
    \item an injective homomorphism $\Psi_\varphi\colon G_A\arr G_B$ for every arrow $\varphi$ of $\Br(\Delta)$ from $A$ to $B$;
    \item an element $t_{\varphi_1, \varphi_2}\in G_C$ for every path $(\varphi_1, \varphi_2)$ of two arrows of $\Br(\Delta)$;
\end{enumerate}
such that they satisfy the identities~\eqref{eq:comgroup1} and~\eqref{eq:compgroup2}.
\end{defi}

\begin{figure}
    \centering
    \includegraphics{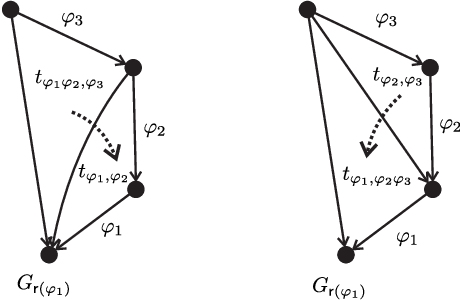}
    \caption{Twisting elements}
    \label{fig:grcomplexid}
\end{figure}

Geometrically, the twisting elements are associated with the triangles of $\Br(\Delta)$ and are elements of the local group of the vertex of the triangle that is the barycenter of a face of the lowest dimension among them. The identity~\eqref{eq:compgroup2} describes connection between the twisting elements associated with faces of a tetrahedron, see Figure~\ref{fig:grcomplexid}.

\begin{example}
    If the complex $\Delta$ in Definition~\ref{def:complexofgroupstriple} is one-dimensional, i.e., if it is a graph, then there are no directed paths of length 2 in $\Br(\Delta)$, so there are no twisting elements $t_{\varphi_1, \varphi_2}$. Therefore, a \emph{graph of groups} is given by a graph $\Delta$, local \emph{vertex groups} and \emph{edge groups} associated to the vertices and edges of $\Delta$, and monomorphisms from the edge groups to the vertex groups of the vertices they connect. For example, if $\Delta$ is one vertex $v$ and a loop, then a graph of groups is defined by groups $G_e$, $G_v$, and two embeddings $G_e\arr G_v$. Theory of graphs of groups and their fundamental groups is developed in~\cite{serr:trees}.
\end{example}

\begin{example}
Consider the group $G$ of affine transformations of $\R^2$ of the form $\vec x\mapsto \pm\vec x+\vec a$ for $\vec a\in\Z^2$. Let us tile $\R^2$ by the images of the squares $A_1=[0, 1/2]\times [0, 1/2]$ and $A_2=[0, 1/2]\times [1, 1/2]$. Their union is a fundamental domain of the action. Let $\Delta$ be the 2-dimensional polytopal complex of the tiling.

Let us describe the associated complex of groups $\Delta/G$. We choose the following transversal $\mathcal{T}$. Two-dimensional cells are $A_1$ and $A_2$. One-dimensional and zero-dimensional cells are the sides and vertices of the square $A_1$. The sides and vertices of $A_2$ that are not common with $A_1$ are not chosen, as they are images of the sides of $A_1$ under the action of $G$.

The complex $|\Delta/G|$ is a ``pillow'' obtained by gluing the squares $A_1$ and $A_2$ along their boundaries: the side $[(0, 1), (1/2, 1)]$ is glued to $[(0, 0), (1/2, 0)]$ by the map $\vec x\mapsto \vec x-(0, 1)$, the side $([1/2, 1), (1/2, 1/2)]$ to $[(1/2, 0), (1/2, 1/2)]$ by $\vec x\mapsto -\vec x+(1, 1)$, and the side $[(0, 1/2), (0, 1)]$ to $[(0, 1/2), (0, 0)]$ by $\vec x\mapsto -\vec x+(0, 1)$.

The set of morphisms of the scwol $\mathcal{G}(\Delta/G)$ is described by the barycentric subdivision of the pillow. The local groups are trivial except for the ones of the vertices, which are of order 2. The twisting elements $t_{\varphi_1, \varphi_2}$ label triangles of the barycentric subdivision. They depend on the choice of representatives of the elements of $\mathcal{G}(\Delta/G)$. If we choose the identification of the vertices $(0, 1)$ and $(1/2, 1)$ with $(0, 0)$ and $(1/2, 0)$ to be the same as the identification of the corresponding sides of the squares $A_2$ and $A_1$ (i.e., the map $\vec x-(0, 1)$), then the twisting elements will be trivial for all triangles except for the triangles of $A_2$ adjacent to the sides $[(0, 0), (0, 1/2)]$ and $[(1/2, 1/2), (1/2, 1)]$, where they are equal to the generator of the local group of the corresponding vertex of $\Delta/G$.
\end{example}

\subsection{Morphisms between complexes of groups}

\begin{defi}
\label{def:morphismsasfunctors}
Let $\Delta_i$, for $i=1, 2$, be complexes of groups (in the formalism of Definition~\ref{def:polyhedralcomplexofgroupsgeneral}). A \emph{morphism of complexes of groups} $f\colon\Delta_1\arr\Delta_2$ is given by 
\begin{enumerate}
    \item a functor of categories $\tilde f\colon\mathcal{CG}(\Delta_1)\arr\mathcal{CG}(\Delta_2)$.
    \item a collection of continuous maps between polytopes $f\colon A\arr \tilde f(A)$ agreeing (in the sense of Definition~\ref{def:morphismofpcomplexes}) with the functor $\tilde f\colon\mathcal{CG}(\Delta_1)\arr\mathcal{CG}(\Delta_2)$.
\end{enumerate}
\end{defi}

Note that since the attaching maps are invariant under the composition with automorphisms of sources and ranges, one has to check condition (2) only for the functor induced on the associated scwols $\mathcal{C}(\Delta_i)$.

Let us describe how to define functors between complexes of groups in the setting of Definition~\ref{def:complexofgroupstriple}. Again, the construction will be less elegant than Definition~\ref{def:morphismsasfunctors}, but it will be useful for describing concrete examples.

Let $f\colon\mathcal{CG}_1\arr\mathcal{CG}_2$ be a functor between abstract complexes of groups. Let $\mathcal{C}_i$ be the associated scwols, and let us choose the representatives $\overline\varphi\in\mathcal{CG}_i$ of morphisms $\varphi\in\mathcal{C}_i^{(1)}$, as before. 

The functor $f$ defines homomorphisms $f\colon G_\alpha\arr G_{f(\alpha)}$ of the local groups.
We have $f(\overline\varphi)=q_\varphi\overline{f(\varphi)}$ for some element $q_\varphi\in G_{\en(f(\overline\varphi))}$.
The elements $q_\varphi$ satisfy the identities
\begin{equation}
    \label{eq:functor1}
    f(\Psi_\varphi(g))=q_\varphi\Psi_{f(\varphi)}(f(g))q_\varphi^{-1}
\end{equation}
and
\begin{equation}
    \label{eq:functor2}
    f(t_{\varphi_1, \varphi_2})q_{\varphi_1\varphi_2}=q_{\varphi_1}\Psi_{f(\varphi_1)}(q_{\varphi_2})t_{f(\varphi_1), f(\varphi_2)},
\end{equation}
where $\Psi_{f(\varphi)}$ is the identity map, if $f(\varphi)$ is an identity automorphism.

Conversely, a functor $f\colon\mathcal{C}_1\arr\mathcal{C}_2$ of the scwols, a collection of elements $q_\varphi\in G_{\en(f(\varphi))}$, and a collection of group homomorphisms $f\colon G_\alpha\arr G_{f(\alpha)}$ defines a functor $\mathcal{CG}_1\arr\mathcal{CG}_2$ by
\[(g, \varphi)\mapsto (f(g)q_\varphi, f(\varphi)).\]

This leads to the following definition of morphisms of complexes of groups.

\begin{defi}
\label{def:morphismofcomplexestr}
Let $\Delta_i$, for $i=1, 2$, be polytopal complexes of groups defined by the data of Definition~\ref{def:complexofgroupstriple}. A \emph{morphism} $f\colon\Delta_1\arr\Delta_2$ is defined by the following data:
\begin{enumerate}
\item A morphism $f\colon\Delta_1\arr\Delta_2$ of polytopal complexes over a functor $\tilde f\colon\mathcal{C}(\Delta_1)\arr\mathcal{C}(\Delta_2)$ of the associated scwols;
\item homomorphisms $\tilde f\colon G_A\arr G_{\tilde{f}(A)}$ for every face $A$ of $\Delta_1$;
\item elements $q_\varphi\in G_{\tilde f(\en(\varphi))}$ for every $\varphi\in\mathcal{C}(\Delta_1)$;
\end{enumerate}
satisfying the identities~\eqref{eq:functor1} and~\eqref{eq:functor2}.
\end{defi}

A special class of morphisms between complexes of groups are \emph{coverings}.  

\begin{defi}
\label{def:covering}
Let $\mathcal{CG}$ be an abstract complex of groups. Let $P\colon X\arr\mathcal{CG}^{(0)}$ be a surjective map, where $X$ is a set. 
An \emph{action} of $\mathcal{CG}$ on $X$ \emph{over $P$} is a map $(\varphi, x)\mapsto \varphi(x)$ from the set of elements $(\varphi, x)$ of $\mathcal{CG}^{(1)}\times X$ such that $P(x)=\be(\varphi)$ to $X$ satisfying the following conditions:
\begin{enumerate}
\item For every $\varphi\in\mathcal{CG}^{(1)}$ the map 
\[x\mapsto\varphi(x)\] is a bijection from $P^{-1}(\be(\varphi))$ to $P^{-1}(\en(\varphi))$.
\item $\varphi_2(\varphi_1(x))=(\varphi_2\varphi_1)(x)$ 
for all $\varphi_1, \varphi_2$ and $x$ such that $P(x)=\be(\varphi_2)$ and $\en(\varphi_2)=\be(\varphi_1)$.
\end{enumerate}
\end{defi}

\begin{defi}
Consider an action of a complex of groups $\mathcal{CG}$ on a set $X$ over a map  $P\colon X\arr\mathcal{CG}^{(0)}$.
    The associated \emph{covering complex} is the category $\mathcal{CG}_X$ with the set of objects $X$ and morphisms $(\varphi, x)$, for $P(x)=\be(\varphi)$. A morphism $(\varphi, x)$ has domain $x$ and range $\varphi(x)$. They are composed according to the rule
\[(\varphi_2, \varphi_1(x))\cdot(\varphi_1, x)=(\varphi_2\varphi_1, x).\]

The map $(\varphi, x)\mapsto\varphi:\mathcal{CG}_X^{(1)}\arr\mathcal{CG}^{(1)}$ is a functor defining the associated \emph{covering map}, which we will also denote by $P$.
\end{defi}

Let us describe when a functor between polytopal complexes of groups (in the sense of Definition~\ref{def:morphismofcomplexestr}) defines a covering, i.e., when we can compose it with an equivalence of complexes of groups to get a covering in the sense of the above definition. Again, it will be much less elegant than Definition~\ref{def:covering}.

Let $f\colon\Delta_1\arr\Delta_2$ be a morphism given by the data in Definition~\ref{def:morphismofcomplexestr}. Let $\varphi_1\in\mathcal{CG}(\Delta_1)^{(1)}$, and denote $\varphi_2=\tilde f(\varphi_1)$. For every $g\in G_{\be(\varphi_1)}$ we have, by~\eqref{eq:functor1}:
\[q_{\varphi_1}\Psi_{\varphi_2}(\tilde f(g))q_{\varphi_1}^{-1}=\tilde f(\Psi_{\varphi_1}(g)),\]
hence the homomorphism $x\mapsto q_{\varphi_1}\Psi_{\varphi_2}(x)q_{\varphi_1}^{-1}$ from $G_{\be(\varphi_2)}$ to $G_{\en (\varphi_2)}$ maps $\tilde f(G_{\be(\varphi_1)})$ to $\tilde f(G_{\en(\varphi_2)})$.

Consequently, the homomorphism $x\mapsto q_{\varphi_1}\Psi_{\varphi_2}(x)q_{\varphi_1}^{-1}$ induces a well defined map $G_{\be(\varphi_2)}/\tilde f(G_{\be(\varphi_1)})\arr G_{\en(\varphi_2)}/\tilde f(G_{\en(\varphi_1)})$ between the sets of cosets. Consequently, the same is true for the map \[x\mapsto \Psi_{\varphi_2}(x)q_{\varphi_1}^{-1}.\]

\begin{defi}
\label{def:polyhedralcovering}
Let $f\colon\Delta_1\arr\Delta_2$ be a morphism of complexes of groups given by the data in Definition~\ref{def:morphismofcomplexestr}. It is a \emph{covering} if the following conditions are satisfied:
\begin{enumerate}
\item The map $|\Delta_1|\arr|\Delta_2|$ induced by $f$ is surjective, and the map $f\colon A\arr\tilde f(A)$ is a surjective isometry for every face $A$ of $\Delta_1$.
\item The homomorphisms $\tilde f\colon G_A\arr G_{\tilde f(A)}$ are injective.
\item For every attaching map $\varphi\colon A\dashrightarrow B$ in $\Delta_2$, and for every $B_1\in f^{-1}(B)$, the map
\[\bigsqcup_{\tilde f(\varphi_1)=\varphi, \en(\varphi_1)=B_1} G_{\tilde f(\be(\varphi_1))}/\tilde f(G_{\be(\varphi_1)})\arr G_B/\tilde f(G_{B_1})\]
equal to the union of the maps induced by $x\mapsto\Psi_\varphi(x)q_{\varphi_1}^{-1}$ is a bijection.
\end{enumerate}
\end{defi}

Since $f\colon A\arr\tilde f(A)$ is surjective, we will usually write $f(A)$ instead of $\tilde f(A)$, when $f$ is a covering.

\begin{proposition}
\label{pr:coveringofpolyhedralcomplexes}
Suppose that a morphism $f\colon\Delta_1\arr\Delta_2$ satisfies the conditions of Definition~\ref{def:polyhedralcovering}. Then there exists an action of $\mathcal{CG}(\Delta_2)$ on a set $X$ and an equivalence of categories $\mathcal{CG}(\Delta_1)\hookrightarrow\mathcal{CG}(\Delta_2)_X$ such that the composition of the equivalence with the covering $P\colon\mathcal{CG}(\Delta_2)_X\arr\mathcal{CG}(\Delta_2)$ is equal to $\tilde f$.
\end{proposition}

\begin{proof}
Let us define $X$ as the set $\bigsqcup_{A\in \Delta_1}G_{f(A)}/\tilde f(G_A)$ of all cosets 
$h\tilde f(G_A)$ in the groups $G_{f(A)}$ for all faces $A$ of $\Delta_1$. The map $P\colon X\arr\mathcal{CG}(\Delta_2)^{(0)}$ is given by $P(hG_{f(A)})=f(A)$.

Let $\varphi\in\mathcal{G}(\Delta_2)^{(1)}$ and $A=\be(\varphi), B=\en(\varphi)$. Let $A_1\in\mathcal{CG}(\Delta_1)^{(0)}$ be such that $f(A_1)=A$. By Condition (1) of Definition~\ref{def:polyhedralcovering}, there is a unique element $\varphi_1\in\mathcal{G}(\Delta_1)^{(1)}$ such that $\be(\varphi_1)=A_1$ and $\tilde f(\varphi_1)=\varphi$. Consider a coset $h\tilde f(G_{A_1})\in X$, and let $(g, \varphi)\in\mathcal{CG}^{(1)}(\Delta_2)^{(1)}$. Then we set
\begin{equation}
    \label{eq:coveringaction}
(g, \varphi)(h\tilde f(G_{\be(\varphi_1)}))=g\Psi_\varphi(h)q_{\varphi_1}^{-1}\tilde f(G_{\en(\varphi_1)}).
\end{equation}

By Condition (3) of Definition~\ref{def:polyhedralcovering}, the map $(g, \varphi)\colon P^{-1}(\be(\varphi))\arr P^{-1}(\en(\varphi))$ is a bijection. 
Condition (2) of Definition~\ref{def:covering} is checked directly. This implies that we get an action. 

The covering complex of groups $\mathcal{CG}(\Delta_2)_X$ has $X$ as the set of objects. The morphisms are triples $(g, \varphi, h\tilde f(G_A))$, where $(g, \varphi)\in\mathcal{CG}(\Delta_2)^{(1)}$, and $f(A)=\be(\varphi)$, and $h\in G_{\tilde f(A)}$. 
They are composed according to the rule
\[(g_1, \varphi_1, h_1\tilde f(G_{A_1}))\circ (g_2, \varphi_2, h_2\tilde f(G_{A_2}))=
((g_1, \varphi_1)\circ (g_2, \varphi_2), h_2\tilde f(G_{A_2})),\]
where $(g_1, \varphi_1)\circ (g_2, \varphi_2)$ is the composition in $\mathcal{CG}(\Delta_1)$ given by~\eqref{eq:compositionrule}, and the composition above is defined if and only if $(g_2, \varphi_2)(h_2\tilde f(G_{A_2}))=h_1\tilde f(G_{A_1})$, using~\eqref{eq:coveringaction}.

Let us construct an isomorphic embedding $M\colon\mathcal{CG}(\Delta_1)\arr\mathcal{CG}(\Delta_2)_X$. We define it on the objects by setting $M(A)$ to be equal to the coset $1\tilde f(G_A)$ in $G_{f(A)}$, and on morphisms by
\[M(g, \varphi)=(\tilde f(g)q_\varphi, \tilde f(\varphi), 1\tilde f(G_{\be(\varphi)})).\]

We have $\be(M(g, \varphi))=1\tilde f(G_{\be(\varphi)})=M(\be(g, \varphi))$ and 
\begin{multline*}
    \en(M(g, \varphi))=(\tilde f(g)q_\varphi, \tilde f(\varphi))(1\tilde f(G_{\be(\varphi)}))=\\
    \tilde f(g)q_\varphi\Psi_{\tilde f(\varphi)}(1)q_\varphi \tilde f(G_{\en(\varphi)})=1\tilde f(G_{\en(\varphi)})=M(\en(g, \varphi)).
\end{multline*}

Let us check that $M$ is a functor. Let $(g_1, \varphi_1), (g_2, \varphi_2)\in\mathcal{CG}(\Delta_2)^{(1)}$ be composable morphisms. Then
\begin{multline*}
M(g_1, \varphi_1)M(g_2, \varphi_2)=\\ (\tilde f(g_1)q_{\varphi_1}, \tilde f(\varphi_1), 1\tilde f(G_{\be(\varphi_1)})(\tilde f(g_2)q_{\varphi_2}, \tilde f(\varphi_2), 1\tilde f(G_{\be(\varphi_2)})=\\
((\tilde f(g_1)q_{\varphi_1}, \tilde f(\varphi_1))(\tilde f(g_1)q_{\varphi_1}, \tilde f(\varphi_1)), 1\tilde f(G_{\be(\varphi_2)})=\\
(\tilde f(g_1)q_{\varphi_1}\Psi_{\tilde f(\varphi_1)}(\tilde f(g_2))q_{\varphi_1}^{-1}q_{\varphi_1}\Psi_{\tilde f(\varphi_1)}(q_{\varphi_2})t_{\tilde f(\varphi_1), \tilde f(\varphi_2)}, \tilde f(\varphi_1\varphi_2), 1\tilde f(G_{\be(\varphi_2)})=\\
(\tilde f(g_1\Psi_{\varphi_1}(g_2)t_{\varphi_1, \varphi_2})q_{\varphi_1\varphi_2}, \tilde f(\varphi_1\varphi_2), 1\tilde f(G_{\be(\varphi_2)})=\\
M(g_1\Psi_{\varphi_1}(g_2)t_{\varphi_1, \varphi_2}, \varphi_1\varphi_2)=M((g_1, \varphi_1)(g_2, \varphi_2)).
\end{multline*}

Consider a morphism $(g, \varphi, h\tilde f(G_{A_1}))\in\mathcal{CG}(\Delta_1)_X^{(1)}$ from $M(A_1)=1\tilde f(G_{A_1})$ to $M(A_2)=1\tilde f(G_{A_2})$. Then $1\tilde f(G_{A_2})=\en(g, \varphi, 1f\tilde (G_{A_1}))=g\Psi_{\tilde f(\varphi)}(1)q_\varphi^{-1}\tilde f(G_{A_2})$, so $gq_\varphi^{-1}\in \tilde f(G_{A_2})$, i.e., $g=\tilde f(g_1)q_\varphi$ for some $g_1\in G_{A_2}$.
Therefore, the range of the functor $M\colon\mathcal{CG}(\Delta_1)^{(1)}\arr\mathcal{CG}(\Delta_2)_X^{(1)}$ is a full subcategory.
\end{proof}

\begin{example}
\label{ex:coveringstar}
Let $\G_1$ be a star consisting of $k$ edges $e_1, e_2, \ldots, e_k$ attached to one vertex $O_1$. Let us denote the other vertex (the leaf) of an edge $e_i$ by $A_i$. We assume that all the vertices $A_1, A_2, \ldots, A_k$ are pairwise different. Let $\G_0$ be the graph consisting of one edge $e$ connecting two vertices $A$ and $O$. Let $f\colon\G_1\arr\G_0$ be given by $f(e_i)=e$, $f(O_1)=O$, $f(A_i)=A$.

Let us now introduce structures of graphs of groups on $\G_1$ and $\G_0$ by choosing local groups $G_1=G_{O_1}, G_0=G_O$ at the vertices $O_1$ and $O$ of $\G_1$ and $\G_0$, respectively. All the other local groups are trivial. The morphism $f\colon\G_1\arr\G_0$ is then given by a homomorphism $\tilde f\colon G_1\arr G_0$ and twisting elements $q_{e_i}\in G_0$ associated with the attaching maps from $e_i$ to $O_1$.

Let us check for which choices of $\tilde f\colon G_1\arr G_0$ and $q_{e_i}$ the conditions of Definition~\ref{def:polyhedralcovering} are satisfied. Condition (1) is satisfied for any choice. Condition (2) requires $\tilde f\colon G_1\arr G_0$ to be injective. Condition (3) has to be checked only for the attaching map from $e$ to $O$ in $\G_0$. It requires that the map from the formal disjoint union $U$ of $k$ single-element sets associated with edges $e_i$ to $G_0/\tilde f(G_1)$ mapping the singleton associated with $e_i$ to the coset $q_{e_i}^{-1}\tilde f(G_1)$ is a bijection $U\arr G_0/\tilde f(G_1)$. This is equivalent to the condition that $\{q_{e_i}^{-1}\}_{i=1, \ldots, k}$ is a  left coset transversal of $G_0$ modulo $\tilde f(G_1)$. Equivalently, $\{q_{e_i}\}_{i=1, \ldots, k}$ is a right coset transversal.

Thus, the morphism $f$ is a covering of complexes of groups if and only if the local group $G_1$ of $O_1$ is a subgroup of index $k$ of the local group $G_0$ of $O$, and the twisting elements form a right coset transversal of $G_0$ modulo $G_1$.
\end{example}

\subsection{Fundamental groups}
\label{ss:fundamentalgroups}

Let $\mathcal{CG}$ be a complex of groups (in the sense of Definition~\ref{def:complexofgroups}). For $\varphi\in\mathcal{CG}^{(1)}$, we consider its \emph{formal inverse} $\varphi^{-1}$ with domain and range equal to the range and domain, respectively, of $\varphi$. The notion of a composable sequence of morphisms and their formal inverses are defined in the same way as for morphisms, with the notion of domain and range defined above. (We do not compose composable sequences, though.)

\begin{defi}
\label{def:pathsinG}
Let $A\in\mathcal{CG}^{(0)}$.  A \emph{path} in $\mathcal{CG}$ is a composable sequence $(\varphi_1, \varphi_2, \ldots, \varphi_n)$ of elements of $\mathcal{CG}^{(1)}$  and their formal inverses. Its \emph{beginning} is the domain of $\varphi_n$. Its \emph{end} is the range of $\varphi_1$. We define $(\varphi^{-1})^{-1}=\varphi$.

Two paths are \emph{homotopic} if one can transform one into the other using the following operations and their inverses:
\begin{enumerate}
\item removing a consecutive pair of the form $(\varphi^{-1}, \varphi)$ or $(\varphi, \varphi^{-1})$;
\item replacing a pair $(\varphi_i, \varphi_{i+1})$ by $\varphi_i\varphi_{i+1}$ or a pair $(\varphi_i^{-1}, \varphi_{i+1}^{-1})$ by $(\varphi_{i+1}\varphi_i)^{-1}$.
\end{enumerate}

The \emph{fundamental group} $\pi_1(\mathcal{CG}, A)$ of $\mathcal{CG}$ is the set of homotopy classes of paths starting and ending in $A$ with the group operation induced by the concatenation operation on sequences and the inverse $(\varphi_1, \varphi_2, \ldots, \varphi_n)^{-1}=(\varphi_n^{-1}, \ldots, \varphi_2^{-1}, \varphi_1^{-1})$.
\end{defi}

\begin{example}
    If all local groups $G_A$ of the complex of groups $\mathcal{CG}$ are trivial (i.e., $\mathcal{CG}$ is a scwol), then the fundamental group $\pi_1(\mathcal{CG}, A)$ is naturally isomorphic to the fundamental group of the geometric realization of $\mathcal{CG}$.
\end{example}

As always, the corresponding translation into the setting of Definition~\ref{def:complexofgroupstriple} is a little more complicated. Consider a complex of groups as in Definition~\ref{def:complexofgroupstriple} with the underlying complex $\Delta$, local groups $G_A$, homomorphisms $\Psi_\varphi$, and twisting elements $t_{\varphi_1, \varphi_2}$. Then, paths are represented by sequences (formal products)
\[g_1\cdot\varphi_1\cdot g_2\cdot \varphi_2\cdots\varphi_n\cdot g_{n+1},\]
where $g_i\in G_{A_i}$, and $\varphi_i$ is an oriented edge in $\Br(\Delta)$ from $A_{i+1}$ to $A_i$ (the orientation does not have to be the canonical one). We concatenate paths as functions: from right to left.

Two paths are homotopic if one can be transformed into the other using the following operations.
(Here $\varphi, \varphi_1, \varphi_2$ are edges of $\Br(\Delta)$ with the \emph{canonical} orientation, i.e., going in the direction towards faces of lower dimension.) 
\begin{alignat*}{2}
    \varphi\cdot h&=\Psi_\varphi(h)\cdot\varphi,&\qquad h\cdot\varphi^{-1}&=\varphi^{-1}\cdot\Psi_\varphi(h),\\
    \varphi_1\cdot  1\cdot \varphi_2&=t_{\varphi_1, \varphi_2}\cdot(\varphi_1\varphi_2), &\qquad \varphi_2^{-1}\cdot 1\cdot \varphi_1^{-1}&=(\varphi_1\varphi_2)^{-1}\cdot t_{\varphi_1, \varphi_2}^{-1},\\
\varphi\cdot 1\cdot \varphi^{-1}&=1,&\qquad \varphi^{-1}\cdot 1\cdot \varphi&=1.
\end{alignat*}
Here $\varphi_1\varphi_2$ represents the composition of the elements of the scwol, i.e., third side of the triangle of $\Br(\Delta)$ associated with the directed path $(\varphi_1, \varphi_2)$.

\begin{example}
\label{ex:uncovofcomplex}
Let $\mathcal{CG}$ be a connected complex of groups. Let $X$ be the set of homotopy classes $[\gamma]$ of paths $\gamma$ starting in a fixed point $A\in\mathcal{CG}^{(0)}$, and denote by $P([\gamma])$ the end of $\gamma$. Then \[\varphi((\varphi_1, \varphi_2, \ldots, \varphi_n))=(\varphi, \varphi_1, \varphi_2, \ldots, \varphi_n)\]
is an action of $\mathcal{CG}$ on $X$. The covering $P\colon\mathcal{CG}_X\arr\mathcal{CG}$ defined by this action is called the \emph{universal covering} of $\mathcal{CG}$. The complex $\mathcal{CG}$ is developable if and only if $\mathcal{CG}_X$ is a scwol, i.e., if all its local groups are trivial. A connected polytopal complex of groups $\G$ is called \emph{developable} if and only if it is isomorphic to the complex of groups $\Delta/G$ defined by an action of a group $G$ on a polytopal complex $\Delta$, see~\cite[Proposition III.C.3.9]{bridhaefl}. 
\end{example}

Every functor $f\colon\mathcal{CG}_1\arr\mathcal{CG}_2$ transforms paths into paths by
\[f(\varphi_1, \varphi_2, \ldots, \varphi_n)=(f(\varphi_1), f(\varphi_2), \ldots, f(\varphi_n)),\]
where we define $f(\varphi^{-1})=f(\varphi)^{-1}$. The action preserves the homotopy classes, so it defines a homomorphism $f_*\colon\pi_1(\mathcal{CG}_1, A)\arr\pi_1(\mathcal{CG}_2, f(A))$ of fundamental groups.

If a morphism $f\colon\mathcal{CG}(\Delta_1)\arr\mathcal{CG}(\Delta_2)$ is given by the data in Definition~\ref{def:morphismofcomplexestr}, then it acts on paths by the rule
\[f(g_1\cdot \varphi_1\cdot g_2\cdot \varphi_2\cdots \varphi_n\cdot g_{n+1})=f(g_1)\cdot f'(\varphi_1)\cdot f(g_2)\cdot f'(\varphi_2)\cdots f'(\varphi_n)\cdot f(g_{n+1}),\]
where $f'(\varphi)=q_\varphi f(\varphi)$ and $f'(\varphi^{-1})=f(\varphi)^{-1}q_\varphi^{-1}$ for canonically directed $\varphi$.

We can \emph{lift} paths by coverings. Let us describe the process of lifting for both definitions of coverings.

Consider the covering $P\colon\mathcal{CG}_X\arr\mathcal{CG}$ defined by an action of a complex of groups $\mathcal{CG}$ on a set $X$. Let $\gamma=(\varphi_1, \varphi_2, \ldots, \varphi_n)$ be a path, and let $x\in X$ be such that $P(x)=\be(\varphi_n)$. Then the lift of $\gamma$ starting in $x$ is the unique path
\[((\varphi_1, x_1), (\varphi_2, x_2), \ldots, (\varphi_n, x_n))\]
mapped to $\gamma$ by $P$ such that $x_n=x$. It is determined by the condition that for every $i=n, \ldots, 2$, we have $x_{i-1}=\phi_i(x_i)$ if $\phi_i\in\mathcal{CG}^{(1)}$ and $\phi_i^{-1}(x_{i-1})=x_i$ if $\phi_i^{-1}\in\mathcal{CG}^{(1)}$. 

Suppose now that $f\colon\Delta_1\arr\Delta$ is a covering in the sense of Definition~\ref{def:polyhedralcovering}. 
In this case, the endpoints of the lifts are, by definition, the elements of the set $X$ constructed in the proof of Proposition~\ref{pr:coveringofpolyhedralcomplexes}, i.e., cosets $h\tilde f(G_A)$ in $G_{f(A)}$ for faces $A$ of $\Delta$. Let us choose, for every face $A$ of $\Delta_1$ a coset representative system $R_A$ of $G_{f(A)}$ modulo $\tilde f(G_A)$. This will give us a choice of isomorphisms between objects of $\mathcal{CG}(\Delta)_X$ and the images of the objects of $\mathcal{CG}(\Delta)$ under the equivalence constructed in the proof of Proposition~\ref{pr:coveringofpolyhedralcomplexes}.

An \emph{$f$-preimage} of a face $A$ of $\Delta$ is, by definition, a pair $(A_1, r\tilde f(G_{A_1}))$, where $A_1\in f^{-1}(A)$. We write it as a pair $(A_1, r)$ for the coset representative $r\in R_{A_1}$.

We will define the lifts (with their endpoints in $X$) of paths consisting of a single element $g\in G_A$, $\varphi\in\mathcal{CG}(\Delta)^{(1)}$ or $\varphi^{-1}$. The general case will follow inductively from the condition that a lift of $\gamma_1\gamma_2$ to a path starting in $x$ is the concatenation $\gamma_1'\gamma_2'$, where $\gamma_2'$ is the lift of $\gamma_2$ starting in $x$, and $\gamma_1'$ is the lift of $\gamma_1$ starting in the end of $\gamma_2'$.

Let $g\in G_A$ for a face $A$ of $\Delta$, and let $(A_1, r)$ be an $f$-preimage of $A$. Then there exists a unique representative $r_1\in R_{A_1}$ such that $gr\in r_1\tilde f(G_{A_1})$. Then the lift of $g$ to a path starting in $(A_1, r)$ 
is the path \[(\tilde f)^{-1}(r_1^{-1}gr),\] and its end is $(A_1, r_1)$. We are using here the fact that $\tilde f$ defines monomorphisms of local groups.

Let now $\varphi$ be an arrow of $\Br(\Delta)$ in the canonical orientation, i.e., going to a face $B$ of lower dimension and identified with an element of $\mathcal{CG}(\Delta)^{(1)}$. Let $(A_1, r)$ be an $f$-preimage of $A$.  Since $f$ is an isometry from the interior of $A_1$ to the interior of $A$, there is a unique edge $\varphi_1$ of $\Br(\Delta_1)$ inside $A_1$ that is mapped to $\varphi$ by $f$. Let $B_1$ be the end of $\varphi_1$. Then $\Psi_\varphi(r)q_{\varphi_1}^{-1}\in G_B$, so there exists a unique element $r_1\in R_{B_1}$ such that $r_1^{-1}\Psi_\varphi(r)q_{\varphi_1}^{-1}\in \tilde f(G_{B_1})$. The end of the lift of the path $(\varphi)$ is then $(B_1, r_1)$, and the lift is the path
\[(\tilde f)^{-1}(r_1^{-1}\Psi_\varphi(r)q_{\varphi_1}^{-1})\cdot \varphi_1.\]

Consider now a path consisting of a single edge $\varphi^{-1}$ in the orientation opposite to the canonical one, i.e., going from a face $A$ to a face $B$ of higher dimension. Then $\varphi\in\mathcal{G}(\Delta)^{(1)}$. Let $(A_1, r)$ be an $f$-preimage of $A$. By Definition~\ref{def:polyhedralcovering}, there exist a unique element $\varphi_1\in\mathcal{CG}(\Delta_1)^{(1)}$ such that $\en(\varphi_1)=A_1$ and $\tilde f(\varphi_1)=\varphi^{-1}$, and a unique coset representative $r_1\in R_{\be(\varphi_1)}$ such that $\Psi_\varphi(r_1)q_{\varphi_1}^{-1}\in r\tilde f(G_{A_1})$. We then define the lift as
\[\varphi_1^{-1}\cdot (\tilde f)^{-1}(q_{\varphi_1}\Psi_\varphi(r_1^{-1})r).\]

\subsection{Homotopy}

Suppose that $f_0, f_1\colon\G_1\arr\G_0$ are two morphisms of polytopal complexes over the same functor $\tilde f\colon\mathcal{C}(\G_1)\arr\mathcal{C}(\G_2)$. Then for every face $A$ of $\G_1$ the images $f_0(A), f_1(A)$ are contained in the same face $\tilde f(A)$. Since faces are convex, we can define a homotopy $F(\cdot, t)\colon A\times[0, 1]\arr\tilde f(A)$ from $f_0|_A$ to $f_1|_A$ by setting $F(\xi, t)$ to be the point $(1-t)f_0(\xi)+tf_1(\xi)$, and, by functoriality of $\tilde f$, these homotopies will agree with the attaching maps, so that we get a homotopy (in the category of continuous maps between the supports) from $f_0$ to $f_1$.

Suppose now that $f\colon\G_1\arr\G_0$ is an arbitrary morphism of polytopal complexes over a functor $\tilde f\colon\mathcal{C}(\G_1)\arr\mathcal{C}(\G_0)$. For a face $A$ of $\G_1$, denote by $\tilde f_1(A)$ the subface of $\tilde f(A)$ of smallest dimension containing the set $f(A)$. Then $\tilde f_1\colon\mathcal{C}(\G_1)\arr\mathcal{C}(\G_0)$ is a functor, which together with the maps $f\colon A\arr\tilde f_1(A)$ defines a morphism of polytopal complexes producing the same induced map $|\G_1|\arr|\G_0|$ as $f$. Denote by $\delta_A$ the attaching map from the face $\tilde f(A)$ to its subface $\tilde f_1(A)$. Then we have $\be(\delta_A)=\tilde f(A)$, $\en(\delta_A)=\tilde f_1(A)$ and for every attaching map $\varphi\in\mathcal{C}(\G_1)^{(1)}$ we have $\delta_{\en(\varphi)}\tilde f(\varphi)=\tilde f_1(\varphi)\delta_{\be(\varphi)}$.

\begin{defi}
We say that two morphisms of complexes of groups $\Delta_1\arr\Delta_0$ are \emph{homotopic} if one can transform one into the other using the following operations and their inverses.
    \begin{enumerate}
        \item \emph{Homotopy relative to the functor.} Changing the maps $f\colon A\arr\tilde f(A)$ of faces without changing the functor $\tilde f$.
        \item \emph{Natural transformation.} We say that a morphism $f_1\colon\Delta_1\arr\Delta_0$ is a natural transformation of $f_0\colon\Delta_1\arr\Delta_0$ if there exists a collection of elements $\delta_A\in\mathcal{CG}(\Delta_2)^{(1)}$, for all faces $A$ of $\Delta_1$, such that the following conditions hold
         \begin{itemize}
             \item $\be(\delta_A)=\tilde f_0(A)$ and $\en(\delta_A)=\tilde f_1(A)$;
             \item $f_0(A)$ is contained in the domain of the attaching map $\delta_A$ and $f_1=\delta_A\circ f_0$;
             \item for all $\varphi\in\mathcal{CG}(\Delta_1)^{(1)}$ we have
        \[\delta_{\en(\varphi)}\tilde f_0(\varphi)=\tilde f_1(\varphi)\delta_{\be(\varphi)}.\]
        \end{itemize}
    \end{enumerate}
\end{defi}

We have the following characterization of homotopy, closer to the original topological notion.

\begin{proposition}
    Morphisms of complexes of groups $f_0, f_1\colon\Delta_1\arr\Delta_2$ are homotopic if and only if there exists a morphism $F\colon\Delta_1\times I\arr\Delta_2$, where $I$ is a finite chain of intervals (seen as a polytopal complex), such that restrictions of $F$ to the endpoints of $I$ are the morphisms $f_0, f_1$.
\end{proposition}

Here, the direct product of polytopal complexes of groups is defined in the straightforward way (taking direct products of faces and the corresponding categories).

The proposition is proved by splitting $F$ into homotopies in the interiors of the sub-intervals of $I$ and the attaching maps from the sub-intervals to their endpoints.

It is checked directly that if morphisms $f_1, f_2\colon\Delta_1\arr\Delta_0$ are homotopic, then the induced homomorphisms $f_{1, *}, f_{2, *}\colon\pi_1(\Delta_1)\arr\pi_1(\Delta_2)$ coincide up to an inner automorphism of the fundamental group of $\Delta_2$ (corresponding to the path connecting the $f_1$-image of the basepoint of $\Delta_1$ to its $f_2$-image defined by the homotopy).

\section{Polytopal models of self-similar groups}

\subsection{Subhyperbolic polytopal substitutions}

\begin{defi}
A \emph{virtual endomorphism of complexes of groups} is a pair of morphisms $f, \iota\colon\Delta_1\arr\Delta_0$ of complexes of groups, such that $f$ is a covering of finite degree. We require (unless the opposite is explicitly stated) that the homomorphisms of the local groups induced by $\tilde\iota$ are surjective.
\end{defi}

\begin{example}
\label{ex:virtend}
Let $\G_1, \G_0$ be single vertices with local groups $G_1, G_0$. Any morphism of complexes of groups $\G_1\arr\G_0$ is given just by a group homomorphism $G_1\arr G_0$. Therefore, a virtual endomorphism $f, \iota\colon G_1\arr G_0$ is given by a monomorphism $\tilde f\colon G_1\arr G_0$ such that $\tilde f(G_1)$ is a subgroup of finite index in $G_0$, and a homomorphism $\tilde\iota\colon G_1\arr G_0$. In other words, it is a \emph{virtual endomorphism} $\tilde\iota\circ (\tilde f)^{-1}$ of the group $G_0$.
\end{example}

Let $f, \iota\colon\G_1\arr\G_0$ be a virtual endomorphism of complexes of groups such that $\G_0$ is connected. Choose a face $A$ of $\G_0$ (for example, a vertex). Then the $\pi_1(\G_0, A)$-biset $\bim_{f, \iota}$ and the iterated monodromy group $\img{(f, \iota)}$ are defined in the same way as for topological virtual endomorphisms. We lift $\mathcal{CG}(\G_0)$-paths by $f$ and then map them back by $\tilde\iota$ using the constructions from Subsection~\ref{ss:fundamentalgroups}. 

For every local group $G_B$ of  $\mathcal{CG}(\G_0)$, we have a homomorphism $G_B\arr\pi_1(\G_0, A)$ defined up to an inner automorphism of the fundamental group (depending on a path connecting $B$ to the basepoint $A$). Composing this homomorphism with the epimorphism $\pi_1(\G_0, A)\arr\img{(f, \iota)}$, we get the \emph{iterated monodromy representation} of the local group.

\begin{defi}
Let $f, \iota\colon\G_1\arr\G_0$ be a virtual endomorphism of a complex of groups. We say that the local groups of $\G_0$ are \emph{faithfully represented in $\img{(f, \iota)}$} if the kernel of the iterated monodromy representation of every local group is trivial.
\end{defi}

Note that if the local groups of $\G_0$ are faithfully represented in the iterated monodromy group, then they are also faithfully represented in the fundamental group. The latter implies that the complex of groups $\G_0$ is developable, see Example~\ref{ex:uncovofcomplex}.

\begin{example}
\label{ex:dihedralimg}
Let $\G_0$ be the graph of groups consisting of two vertices $A$ and $B$ with local groups $\langle a\rangle$ and $\langle b\rangle$ of order 2 and one edge $e$ with trivial local group. We denote by $\varphi_A, \varphi_B$ the corresponding morphisms (attaching maps) from the edge to the vertices $A$ and $B$, respectively.

\begin{figure}
    \centering 
    \includegraphics{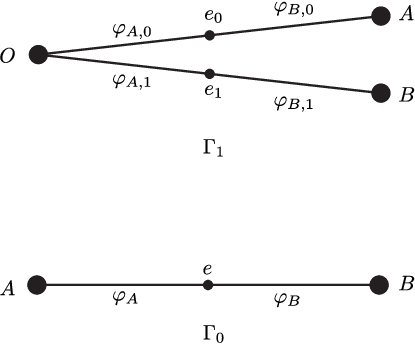}
    \caption{Tent map}
    \label{fig:tentmap}
\end{figure}

Let $\G_1$ be its barycentric subdivision with the same local groups. We denote by $O$ the new vertex, corresponding to the barycenter of the edge of $\G_0$, and by $A$ and $B$ the old vertices. We label the edges of $\G_1$ by $e_0, e_1$, and the attaching maps by $\varphi_{A, 0}, \varphi_{A, 1}$, $\varphi_{B, 0}$, $\varphi_{B, 1}$, as it is shown in Figure~\ref{fig:tentmap}.

Let $\iota\colon\G_1\arr\G_0$ be the natural homeomorphism mapping the vertices $A$ and $B$ of $\G_1$ to the vertices $A$ and $B$ of $\G_0$, and the edges $e_0$ and $e_1$ of $\G_1$ to halves of the edge of $\G_0$. The corresponding functor of scwols can be given by 
\[\tilde\iota\colon \varphi_{A, 0}\mapsto\varphi_A, \quad\varphi_{A, 1}\mapsto\varphi_B\quad \varphi_{B, 0}\mapsto e, \quad\varphi_{B, 1}\mapsto e, \quad e_i\mapsto e.\]
The homomorphisms $G_A\arr G_A$ and $G_B\arr G_B$ induced on the local groups are assumed to be identities, the twisting elements $q_\varphi$ are trivial.

Consider the covering $f\colon\G_1\arr\G_0$ mapping central vertex of $\G_1$ to $A$, as it is shown in Figure~\ref{fig:tentmap}. It is defined by the functor
\[\tilde f\colon\varphi_{A, i}\mapsto\varphi_A,\quad \varphi_{B, i}\mapsto\varphi_B,\quad A\mapsto B,\quad B\mapsto B,\quad O\mapsto A, \quad e_i\mapsto e.\]
The homomorphisms $\tilde f\colon G_A\arr G_A$ and $\tilde f\colon G_B\arr G_A$ are isomorphisms. According to Example~\ref{ex:coveringstar}, we need the set $\{q_{\varphi_{A_0}}, q_{\varphi_{A_1}}\}$ to be equal to $G_A$, in order for $f$ to be a covering of graphs of groups. The other twisting elements should be trivial. Let us choose
\[q_{\varphi_{A, 0}}=a,\qquad q_{\varphi_{A_1}}=1.\]
It is checked directly that the conditions of Definition~\ref{def:polyhedralcovering} for $f$ are satisfied.

Let us compute the iterated monodromy group of the virtual endomorphism $f, \iota\colon\G_1\arr\G_0$. The fundamental group $\pi_1(\G_0, e)$ is generated by the loops $\varphi_A^{-1}\cdot a\cdot\varphi_A$ and $\varphi_B^{-1}\cdot b\cdot\varphi_B$.

Since preimages $e_0$ and $e_1$ are mapped by $\tilde\iota$ back to $e$, we may choose the trivial connecting paths $\ell_0$ and $\ell_1$.

Since $q_{\varphi_{B, i}}$ are trivial and $\tilde f\colon G_A\arr G_B$, $\tilde f\colon G_B\arr G_B$ are isomorphisms, the lifts of $\varphi_B^{-1}\cdot b\cdot\varphi_B$ by $f$ are the paths $\varphi_{B, 0}^{-1}\cdot a\cdot\varphi_{B, 0}$ and $\varphi_{B, 1}^{-1}\cdot b\cdot\varphi_{B, 1}$. Since $q_{\varphi_{A, 0}}=a$ and $q_{\varphi_{A, 1}}=1$, we have
\[f( \varphi_{A, 1}^{-1}\varphi_{A, 0})=f(\varphi_{A, 0}^{-1}\varphi_{A_1})=\varphi_A^{-1}\cdot a\cdot\varphi_A.\]
Applying $\iota$ to the lifts of the generators of $\pi_1(\G_0, e)$, we get the wreath recursion
\[a=(01),\qquad b=(a, b)\]
for $\img(f, \iota)$.
\end{example}

\begin{defi}
Let $f, \iota\colon\G_1\arr\G_0$ be a polytopal substitution. Its \emph{uniformization} is a structure of complexes of groups on $\G_i$ and an extension of $f$ and $\iota$ to morphisms of complexes of groups such that $f$ is a covering.
\end{defi}

\begin{example}
\label{ex:sierpvirtend}
We can uniformize the substitution producing the Sierpi\'nski triangle from Example~\ref{ex:sierptriangle} in a way similar to the previous example. We set the local groups $G_A, G_B, G_C$ to be groups of order two with generators $a, b, c$, respectively. The local group $G_O$ is trivial. Similarly, the local groups of $\G_1$ are trivial, except for the vertices of degree 1, which have groups of order 2 as local groups. The covering $f\colon\G_1\arr\G_0$ is defined in the same way as in the previous example, where $q_\varphi$ are trivial except for one morphism per each vertex of $\G_1$ at which $f$ is not a local homeomorphism. The iterated monodromy group of this virtual endomorphism is given by the wreath recursion
\[a=(01)(1, 1, a),\qquad b=(02)(1, b, 1),\qquad c=(12)(c, 1, 1).\]
This is the \emph{Hanoi tower group}, see~\cite{grisunik:hanoi}.
\end{example}

\begin{example}
Example~\ref{ex:gromovsubst} can be uniformized by endowing the segment $\G_0$ with the structure of a graph of group with two vertex groups $\Z/3\Z$ and trivial edge group. The iterated monodromy group will be
\[a=(234)(a, 1, 1, 1),\qquad b=(321)(1, 1, 1, b),\]
where $a$ and $b$ are the generators of the vertex groups.

The elements $ab$ and $a^{-1}b^{-1}$ satisfy
\[ab=(142)(1, a, 1, b), a^{-1}b^{-1}=(143)(1, 1, a^{-1}, b^{-1}).\]
If we conjugate the right-hand side by $(1, a^{-1}, a, 1)$, which correspond to a change of the basis of the biset, we get
\[ab=(142)(1, 1, 1, ab),\qquad a^{-1}b^{-1}=(143)(1, 1, 1, a^{-1}b^{-1}).\]
The group given by the recursion
\[x=(142)(1, 1, 1, x),\qquad y=(143)(1, 1, 1, y)\]
is equivalent to the iterated monodromy group of $F(z)=\frac{-z^4+6z^2+3}{8z}$. See the Julia set of $F(z)$ in Figure~\ref{fig:JuliaGr}. The dashed circles mark fixed points $i\sqrt{3}$ and $-i\sqrt{3}$ of $F$.

\begin{figure}
\centering
\includegraphics{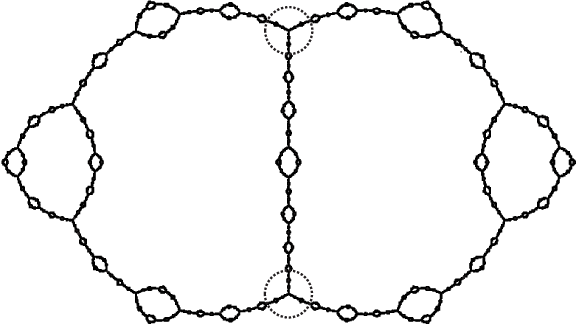}
\caption{The Julia set of $\frac{-z^4+6z^2+3}{8z}$}
\label{fig:JuliaGr}
\end{figure}

The subgroup $\langle ab, a^{-1}b^{-1}\rangle$ is normal of index 3, which implies that the inclusion induces a degree 3 branched covering of the limit space $\G_\infty$ of the substitution by the Julia set of $F(z)$. This covering corresponds to taking quotient of the Julia set by the action of a cyclic group of order 3. One can show that this cyclic group can be realized as the group of order 3 of rotations of the Riemann sphere around the points $i\sqrt{3}$ and $-i\sqrt{3}$. The limit orbispace of $\langle a, b\rangle$ is the orbispace of the action of this group of rotations on the Julia set of $F$. 
\end{example}

\begin{example} Let us describe all uniformizations and the corresponding iterated monodromy groups for the \emph{tent map}, i.e., the polytopal substitution $f, \iota\colon\G_1\arr\G_0$ described in Example~\ref{ex:dihedralimg}.  We restrict ourselves only to uniformizations in which the local groups are faithfully represented in the iterated monodromy group.

We want $f$ to be a degree two covering of graphs of groups. This implies that the edge groups of $\G_0$ and $\G_1$ are trivial, if they are faithfully represented in the iterated monodromy group. Consequently, the local group $G_A$ of $\G_0$ has to be of order 2, since its preimage under $f$ is mapped by $\iota$ to the edge of $\G_0$. Let us denote the non-trivial element of $G_A$ by $a$. Since $f$ is a covering, the functor $\tilde f$ induces isomorphisms of the vertex groups of $\G_1$ with $G_B$. Let us identify them with $G_B$ by these isomorphisms. We get the morphisms $\psi_0\colon G_B\arr G_A$ and $\psi_1\colon G_B\arr G_B$ induced by $\tilde\iota$. Repeating computation of the iterated monodromy group given in Example~\ref{ex:dihedralimg}, we see that it is generated by the wreath recursion
\[a=(01),\qquad g=(\psi_0(g), \psi_1(g)),\]
for $g\in G_B$. It follows from the recursion that $G_B$ acts on the tree $\{0, 1\}^*$ by commuting automorphisms of order two. Therefore, if $G_B$ is finite, then its image in the iterated monodromy group is an elementary abelian 2-group $(\Z/2\Z)^n$ for some $n$. Since we assume that $\tilde\iota$ is surjective on the local groups, $\psi_1$ is an automorphism, and we get that the iterated monodromy groups of such virtual endomorphisms of graphs of groups coincide with the class of groups acting on binary trees described in~\cite{sunic:hausdorff}.

If we do not require that $\tilde\iota$ induces surjective homomorphisms of groups, then the elements of $G_B$ belonging to the kernel of an iteration of $\psi_1$ will act by a finitary automorphism of $\{0, 1\}^*$ (i.e., by an automorphism $g$ having all sections $g|_v$ trivial for $v\in\alb^m$ for some $m$). In particular, the nucleus of the iterated monodromy group will be the same as the nucleus of the iterated monodromy group in which $G_B$ is replaced by the intersection of the ranges of $\psi_1^n$.
\end{example}

\begin{defi}
We say that a polytopal substitution $f, \iota\colon\Delta_1\arr\Delta_0$ is \emph{subhyperbolic} if the maps $\iota\colon A\arr\tilde\iota(A)$ are uniformly contracting for every face $A$ of $\Delta_1$, and there exists a uniformization of $(f, \iota)$ with finite local groups.
\end{defi}

This definition can be seen as a polytopal version of the notion of a subhyperbolic rational function, see~\cite{milnor}.

The next theorem follows from~\cite[Theorem~5.10]{nek:models}, see also~\cite[Theorem~4.5.36]{nek:dyngroups}.

\begin{theorem}
\label{th:limitofcontractingvirtend}
    Let $f, \iota\colon\G_1\arr\G_0$ be a subhyperbolic polytopal substitution. Choose its uniformization with finite local groups. Then the iterated monodromy group $\img{(f, \iota)}$ is contracting and the limit dynamical system $f_\infty\colon\G_\infty\arr\G_\infty$ of the topological correspondence $f, \iota\colon |\G_1|\arr|\G_0|$ is topologically conjugate to the limit dynamical system of the iterated monodromy group. In particular, $\G_\infty$ and $\lims[\img{(f, \iota)}]$ are homeomorphic.
\end{theorem}

\begin{example}
The limit space of the substitution from Example~\ref{ex:sierptriangle} is the Sierpi\'nski triangle, which is homeomorphic to the limit space of the Hanoi tower group, see~\cite{grisunik:hanoi}.
\end{example}

\begin{figure}
\centering
\includegraphics{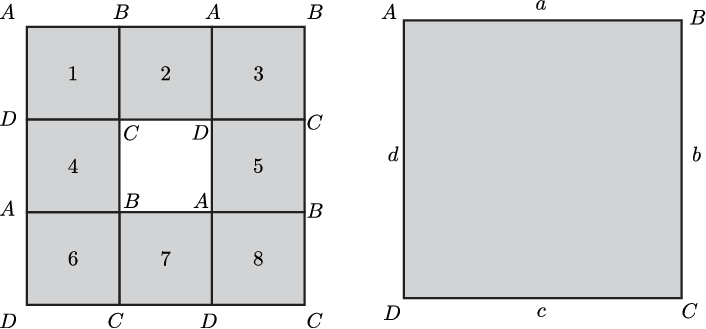}
\caption{Sierpi\'nski carpet substitution}
\label{fig:sierpcarpet}
\end{figure}

\begin{example}
\label{ex:sierpsquare}
Let $\G_0$ be the square $ABCD$ on the right-hand side of Figure~\ref{fig:sierpcarpet}. Let $\G_1$ be the complex obtained from a $3\times 3$ square after removing the central square, shown on the right-hand side. Consider the substitution $f, \iota\colon\G_1\arr\G_0$, where $f$ maps the squares of $\G_1$ isometrically onto $\G_0$ by the label-preserving map (i.e., we are folding $\G_1$ eight times onto $\G_0$). Let $\iota$ be the similarity with coefficient $\frac 13$ from $\G_1$ to a subset of $\G_0$. Iterations of the substitution $(f, \iota)$ will produce the classical approximations  of the Sierpi\'nski carpet.

Let us set the local groups of the edges of the square $\G_0$ to be groups of order 2, and let $a, b, c, d$ be the non-trivial elements of the groups $G_{AB}, G_{BC}, G_{CD}, G_{DA}$, respectively. The local groups of vertices will be generated by the local groups of the adjacent edges, and hence dihedral. Let us see what orders of these dihedral groups should be in order to have a uniformization with faithfully represented local groups. We choose the basepoint and the connecting paths in the interior of the square $\G_0$. Then it follows directly from the description of the substitution, that the elements $a, b, c, d$ of $\img(f, \iota)$ satisfy the recursion
\begin{align*}
a&=(46)(58)(a, a, a, 1, 1, 1, 1, 1),\\
b&=(12)(67)(1, 1, b, 1, b, 1, 1, b),\\
c&=(14)(35)(1, 1, 1, 1, 1, c, c, c),\\
d&=(23)(78)(d, 1, 1, d, 1, d, 1, 1),
\end{align*}
where the digits $1, \ldots, 8$ enumerate the preimages of the basepoint as it is shown in Figure~\ref{fig:sierpcarpet}. 

It is checked directly that the elements $ab$, $bc$, $cd$, $da$ of the iterated monodromy group have order 6. Consequently, the vertex groups should be dihedral groups $D_6$ of order 12, if we want them to be faithfully represented in the iterated monodromy group and be generated by the adjacent edge groups. This group was introduced in~\cite{mattebonnekzheng:liouvills}. It is the first example of an infinite amenable group with property FW, see~\cite{mattebonnekzheng:FW}.
\end{example}

\begin{example}
Consider the substitution from Example~\ref{ex:dendrites}.
The corresponding uniformization can be found in a similar way as in Example~\ref{ex:sierpsquare}. We set the local groups of the two-dimensional faces of $\G_0$ to be groups of order two, except for the face $ABC$, which is assumed to be trivial. 

Let $g_{XYZ}$ denote the non-trivial element of the local group of the face $XYZ$. Let us denote
\begin{alignat*}{3}
g_{A'BC}&=\alpha, &\qquad g_{B'BC}&=a, &\qquad g_{C'BC}&=a', \\
g_{A'AB} &=c , &\qquad 
g_{B'AB} &=c', &\qquad 
g_{C'AB} &=\gamma \\
g_{A'AC} &=b, &\qquad 
g_{B'AC} &=\beta, &\qquad 
g_{C'AC} &=b'
\end{alignat*}

 Then the wreath recursion for the corresponding iterated monodromy group is as follows
\begin{alignat*}{3}
\alpha &=(12)(34), &\qquad a &=(13)(24), &\qquad a' &=(14)(23), \\
c &=(b', b', b, b) , &\qquad 
c' &=(b, b', b, b'), &\qquad 
\gamma &=(\beta, 1, 1, \beta) \\
b &=(a', a', c, c), &\qquad 
\beta &=(\alpha, \gamma, \alpha, \gamma), &\qquad 
b' &=(c', 1, 1, c')
\end{alignat*}

It is checked directly that, in the iterated monodromy group we have 
\[a'=a\alpha,\quad b'=b\beta,\quad c'=c\gamma.\]
We choose the edge and vertex groups to be generated by the local groups of the adjacent two-dimensional faces. See more about this iterated monodromy group in the paper~\cite{nek:dendrites}. In particular, it is shown in~\cite[Proposition~5.5]{nek:dendrites} that some iteration $\iota^n\colon\G_n\arr\G_0$ is contracting on the faces. It is also shown that $\langle\alpha, \beta, \gamma, ac, cb\rangle$ is a self-similar subgroup of index 2 equivalent to the iterated monodromy group of the function $f$ given in~\eqref{ex:dendrites}.
\end{example}

\subsection{Examples of non-subhyperbolic substitutions}

Let us show some examples of polytopal substitutions that are not subhyperbolic.

\begin{example}
If a polytopal substitution $f, \iota\colon\Delta_1\arr\Delta_0$ is such that, for some face $A$, the map $\tilde f\colon\en^{-1}(A)\arr\en^{-1}(f(A))$ is not surjective, then $(f, \iota)$ is not subhyperbolic. In particular, the substitutions described in Example~\ref{ex:wordsubstitution} are not subhyperbolic, if $|E|\ge 2$ and the length of some word $\sigma(e)$ is more than 1. Similarly, the Penrose tiling subdivision rule from Example~\ref{ex:penrose} is not subhyperbolic.
\end{example}

Another necessary condition for a substitution to be subhyperbolic is bounded ramification of its iterations.

\begin{example}
Let $\G_0$ be a segment, and let $\G_1$ be the square, shown on the top of Figure~\ref{fig:diamond}. We choose $f\colon\G_1\arr\G_0$ to be the label-preserving simplicial map. 

\begin{figure}
\centering
\includegraphics{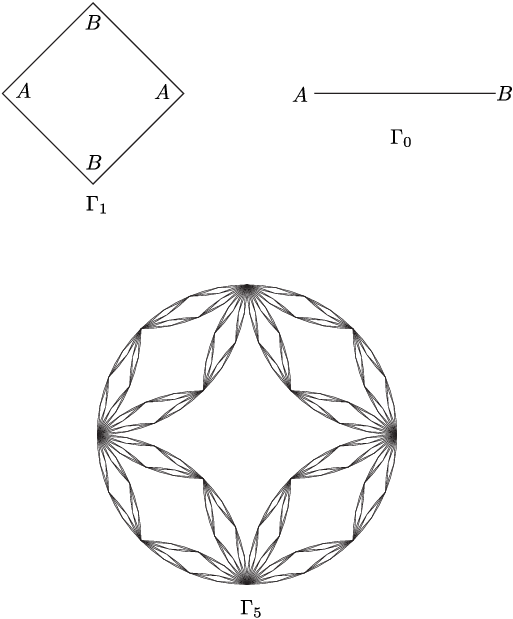}
\caption{Diamond substitution}
\label{fig:diamond}
\end{figure}

We define $\iota\colon\G_1\arr\G_0$ as the map mapping the vertices of $\G_1$ labeled by $A$ to the endpoints $A$ and $B$ of the segment  $\G_0$, and the vertices of $\G_1$ labeled by $B$ to the middle of the segment. The map $\iota$ contracts the lengths of the edges by a factor of $1/2$. The graphs $\G_n$ where used in~\cite{NewmanRabinovich:diamond,JohsonSchechtman:diamond} in their study of Banach spaces. Before that, they appeared in the physics literature~\cite{BerkerOstlund1979}.

 In the bottom of Figure~\ref{fig:diamond} the graph $\G_5$ is shown. It is easy to see that the local degree of the covering $f^n\colon\G_n\arr\G_0$ at the endpoints of $\G_0$ goes to infinity. Therefore, this substitution is not subhyperbolic. 
\end{example}

On the other hand, there are uniformizations of the substitution from the last example, if we drop the condition that the local groups are finite. 

Let $G_A$ be the local group of the vertex $A$ of $\G_0$. Then the local group of the preimage of $A$ in $\G_1$ is a subgroup of index $2$ in $G_A$, and $\tilde\iota$ defines a homomorphism from it to $G_A$, i.e., we get a virtual endomorphism $\phi\colon G_A\dashrightarrow G_A$. This virtual endomorphism defines a structure of a self-similar group on $G_A$ (over a binary alphabet). 
We also have to define a virtual homomorphism $\rho\colon G_A\dashrightarrow G_B$ defined on a subgroup of index 2 of $G_A$, induced by $\tilde\iota$. 

For instance, we can choose $G_B=\langle b\rangle\cong\Z/2\Z$, $G_A=\langle a\rangle\cong\Z$, choose $\phi\colon G_A\dashrightarrow G_A$ to be $n\mapsto n/2$ defined on $2\Z$, and choose $\rho\colon G_A\dashrightarrow G_B$ to be the surjective homomorphism from $2\Z$ to $G_B$. Then $G_A$ is the binary adding machine, and the iterated monodromy group of the uniformization will be the group generated by
\begin{equation}
\label{eq:diamondimg}
a=(12)(34)(a, 1, b, 1),\qquad b=(14)(23).
\end{equation}

However, the limit space of this group is different from $\G_\infty$, i.e., Theorem~\ref{th:limitofcontractingvirtend} does not hold for this uniformization. In fact, different choices of the self-similar group $G_A$ may produce different limit spaces, or non-contracting iterated monodromy groups.

\begin{figure}
\centering
\includegraphics{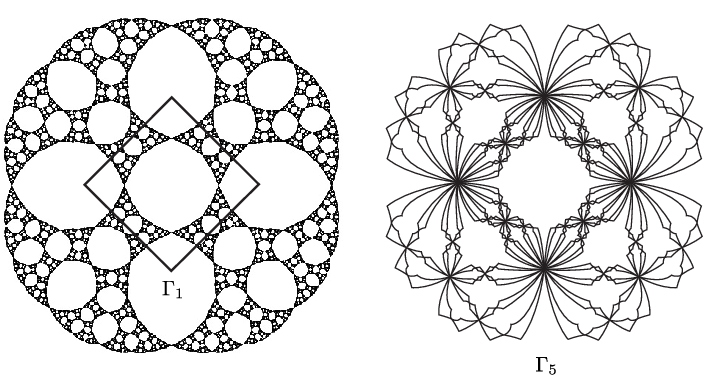}
\caption{The Julia set of $z^2+\frac{1}{16z^2}$ and the graphs $\G_1$ and $\G_5$} 
\label{fig:julia16}
\end{figure}

The group given by the wreath recursion~\eqref{eq:diamondimg} is the iterated monodromy group of the rational function $F(z)=z^2+\frac{1}{16z^2}$. The Julia set of this function, i.e., the limit space of the iterated monodromy group is shown on the left-hand side of Figure~\ref{fig:julia16} together with the square with vertices $\frac 12, -\frac 12, \frac i2, -\frac i2$. The preimage of the square under $F^n$ is isomorphic to $\G_{n+1}$. See, for example, the preimage of the square under $F^4$ on the right-hand side of the figure, and compare it to the graph $\G_5$ in Figure~\ref{fig:diamond}.

Another example of a uniformization of the last example is obtained by choosing $G_A$ to be $\img(z^2-1)$. It produces the iterated monodromy group with the recursion
\[a_0=(12)(34)(a_1, 1, b, 1),\quad a_1=(a_0, 1, 1, 1),\quad b=(14)(23).\]
It is the iterated monodromy group of the rational function $F(z)=z^2+\frac{c}{16 z^2}$, where $c\approx 0.08738$ is the real root of $x^3+5x^2+11x-1$. The Julia set of this function, i.e., the limit space of the iterated monodromy group, is not homeomorphic to the Julia set shown in Figure~\ref{fig:julia16}.

\begin{example}
A classical example of a polytopal substitution which is contracting on the faces, but not subhyperbolic is the \emph{barycentric subdivision rule}, see~\cite[Example~1.3.1]{cfp:subdivision} and~\cite[Example~4.6]{cfp:subdfromrat}.

Let $\G_0$ be the complex consisting of one triangular face (and its sides and vertices). Let $\G_1$ be its barycentric subdivision. Consider the substitution $f, \iota\colon\G_1\arr\G_0$, where $\iota$ is the identity map, and $f$ is the label-preserving branched covering shown in Figure~\ref{fig:barycentric}. Then $\G_n$ is the $n$th barycentric subdivision of $\G_0$. Since the degrees of vertices grow, this substitution is not subhyperbolic.

\begin{figure}
\centering
\includegraphics{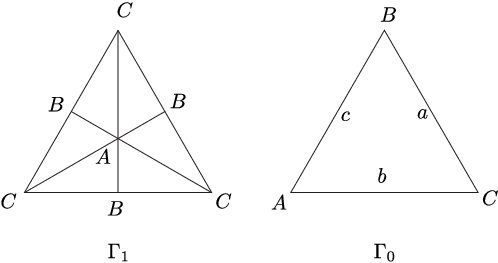}
\caption{Barycentric subdivision of the triangle}
\label{fig:barycentric}
\end{figure}

Same as in the previous example, we can uniformize $(f, \iota)$, if we allow infinite local groups. For example, we can choose the local groups of the sides of the triangle $\G_0$ are of order 2. We label the sides of the triangle $\G_0$ in Figure~\ref{fig:barycentric} by the generators $a, b, c$ of their local groups. The vertex groups will be generated by the adjacent edge groups with identity embeddings and trivial twisting elements. The wreath recursion for the iterated monodromy group is then
\begin{align*}
a &= (a, a, b, b, c, c),\\
b &=(16)(23)(45),\\
c &=(12)(34)(56).
\end{align*}
The vertex groups are dihedral groups $G_A=\langle b, c\rangle\cong  D_6$, $G_B=\langle a, c\rangle\cong D_2$, $G_C=\langle a, b\rangle\cong D_\infty$ of orders 12, 4, and infinity.

\begin{figure}
\centering
\includegraphics{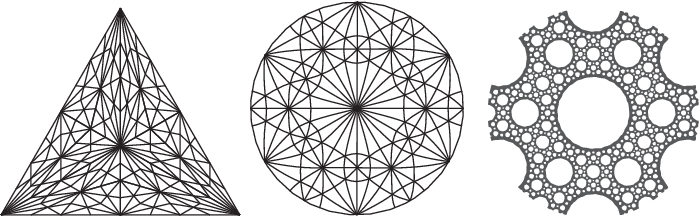}
\caption{Barycentric subdivision and the limit space}
\label{fig:barycentric2}
\end{figure}

This group $G$ is an index two extension of the iterated monodromy group of the rational function 
\[F(z)=\frac{z^6-2i\sqrt3\,z^3+1}{z^6+2i\sqrt3\,z^3+1},\]
and its limit orbispace $\lims$ is the quotient of the Julia set of $F$ by the reflection with respect to the unit circle. See~\cite[Example~4.6]{cfp:subdfromrat} where a conjugate map $F(z)=\frac{4(z^2-z+1)^3}{27z^2(z-1)^2}$ is given. See Figure~\ref{fig:barycentric2}, where the limit space (on the right-hand side of the picture), and the complex $\G_4$ are shown. The complex $\G_4$ are drawn in two ways: as the fourth barycentric subdivision of the triangle, and as the part inside the unit circle of the preimage of the unit circle under $F^4$.
\end{example}

\subsection{$G$-equivariant maps}
\label{ss:gequivariantmaps}

The language of complexes of groups is convenient for specifying virtual endomorphisms and computing their iterated monodromy groups, but is less suited to the proofs that follow. We therefore temporarily switch to the equivalent language of
$G$-complexes and $G$-equivariant maps, following~\cite{nek:models} and~\cite{nek:dimension}.

Let $(G, \bim)$ be a contracting self-similar group. Let us choose a basis $\alb\subset\bim$.
Let $\X$ be a locally compact space with a proper cocompact left action of $G$ on it. The space $\X\otimes\bim$ is the quotient of $\X\times\bim$, where $\bim$ is discrete, by the equivalence relation $x\otimes v=x\cdot g^{-1}\otimes g\cdot v$ for all $x\in\X, g\in G, v\in\bim$. 
The group $G$ naturally acts on $\X\otimes\bim$ by the action induced by the right action of $G$ on $\bim$.

A map $I\colon \X\otimes\bim\arr\X$ is $G$-equivariant if and only if the maps \[I_v=I(\cdot\otimes v)\colon \X\arr\X\] for $v\in\bim$ satisfy
\[I_v(\xi\cdot g)=I_{g\cdot v}(\xi),\qquad I_{v\cdot g}(\xi)=I_v(\xi)\cdot g\]
for all $\xi\in\X$, $v\in\bim$, and $g\in G$.

Consequently, an equivariant map $I\colon \X\otimes\bim\arr\X$ is uniquely determined by a collection of maps $I_x$ for $x\in\alb$ satisfying the condition
\[I_x(\xi\cdot g)=I_{g(x)}(\xi)\cdot g|_x\]
for all $\xi\in\X$, $x\in\alb$, and $g\in G$. Any collection of maps satisfying the above condition defines an equivariant map $I\colon \X\otimes\bim\arr\X$ by the rule $I(\xi\otimes x\cdot g)=I_x(\xi)\cdot g$.

Let $\bim_1, \bim_2$ be $G$-bisets, and let $\X_1, \X_2$ be proper cocompact left $G$-spaces.
If $I_1\colon \X_1\otimes\bim_1\arr\X_2$ and $I_2\colon \X_2\otimes\bim_2\arr\X_3$ are $G$-equivariant maps, then we denote by $I_1\otimes I_2\colon \X_1\otimes\bim_1\otimes\bim_2\arr\X_3$ the map given by
\[I_1\otimes I_2(\xi\otimes v_1\otimes v_2)=I_2(I_1(\xi\otimes v_1)\otimes v_2).\]

In particular, any equivariant map $I\colon \X\otimes\bim\arr\X$ naturally defines its \emph{iterations} $I^{\otimes n}\colon \X\otimes\bim^{\otimes n}\arr\X$, which we will denote just by $I^n$.

We also denote $I_n=I\otimes Id\colon \X\otimes\bim^{\otimes (n+1)}\arr\X\otimes\bim^{\otimes n}$, where $Id$ is the identity map on $\X\otimes\bim^{\otimes n}$. Then $I_n(\xi\otimes x_1x_2\ldots x_{n+1})=I(\xi\otimes x_1)\otimes x_2x_3\ldots x_{n+1}$, and $I^n=I_{n-1}\circ\cdots\circ I_0$.

\begin{defi}
Suppose that $\X$ is a metric space and $G$ acts on it by isometries.
We say that an equivariant map $I\colon \X\otimes\bim\arr\X$ is \emph{contracting} if there exists $\lambda<1$ such that
\begin{equation}
\label{eq:contractingI}
d(I(\xi_1\otimes x), I(\xi_2\otimes x))<\lambda d(\xi_1, \xi_2)
\end{equation}
for all $\xi_1, \xi_2\in\X$ and $x\in\bim$.

We say that it is \emph{asymptotically contracting} if $I^n\colon \X\otimes\bim^{\otimes n}\arr\X$ is contracting for some $n\ge 1$.
\end{defi}

Since $G$ acts by isometries on $\X$, it is enough to check the inequality~\eqref{eq:contractingI} for the maps $I_x$, $x\in\alb$.

Let $\limg$ be the limit $G$-space of $(G, \bim)$. Suppose that a point $\xi\in\limg$ is represented by $\ldots x_2x_1\cdot g\in\xmo\times G$, and let $x\in\alb$. Then we have a well-defined point $\xi\otimes x$ of $\limg$ represented by $\ldots x_2x_1\cdot g\cdot x=\ldots x_2x_1g(x)\cdot g|_x$. This defines a natural $G$-equivariant homeomorphism between $\X\otimes\bim$ and $\X$. 

The following is proved in~\cite[Theorem~2.11]{nek:dimension} and~\cite[Proposition~4.5.30]{nek:dyngroups}.

\begin{theorem}
Let $(G, \bim)$ be a finitely generated contracting group. Then there exists a $G$-invariant metric on $\limg$ such that the natural homeomorphism $\limg\otimes\bim\arr\limg$ is asymptotically contracting.
\end{theorem}

Replacing $G$ by the subgroup generated by the nucleus does not change $\limg$. Therefore, the condition that $G$ is finitely generated is not very restrictive.

The following theorem is proved in~\cite[Theorem~5.9]{nek:models}.

\begin{theorem}
\label{th:limgconvergence}
Let $(G, \bim)$ be a contracting group. Suppose that $I\colon \X\otimes\bim\arr\X$ is a contracting equivariant map. Then the inverse limit of the $G$-spaces $\X\otimes\bim^{\otimes n}$ with respect to the maps $I_n\colon \X\otimes\bim^{\otimes (n+1)}\arr\X\otimes\bim^{\otimes n}$ is equivariantly homeomorphic to the limit space $\limg$. 

If we denote by $I^\infty\colon \limg\arr\X$ the limit of the maps $I^n$, then there exists $\lambda\in (0, 1)$ such that for any compact subsets $T_1\subset\limg$ and $T_2\subset\X$ there exist $C>0$ such that 
\[d(I^{\infty}(\xi_1\otimes v), I^n(\xi_2\otimes v))<C\lambda^n\]
for all $\xi_1\in T_1, \xi_2\in T_2$, $v\in\bim^{\otimes n}$.
\end{theorem}

Let $(G, \bim)$ be a contracting group, and let $\Delta$ be a proper cocompact polytopal $G$-complex.  Then $\Delta\otimes\bim$ is also a proper cocompact polytopal $G$-complex. 

\begin{defi}
A \emph{contracting polytopal model} of $\limg$ is an equivariant morphism of polytopal complexes $J\colon\Delta\otimes\bim \arr\Delta$, where 
\begin{itemize}
\item $\Delta$ is a proper cocompact polytopal $G$-complex;
\item $J\colon A\mapsto\tilde J(A)$ is a contracting map for every face $A$ of $\Delta$;
\end{itemize}
\end{defi}

Let us describe a general method of constructing a contracting polytopal model directly from the recursive description of a self-similar group. 

Let $(G, \bim)$ be a self-similar contracting group, let $\alb\subset\bim$ be a basis, and let $\nuke$ be the associated nucleus. Denote by $\til$ the subset of the limit $G$-space $\limg$ equal to the image of the set of points representable by sequences $\ldots x_2x_1\cdot g$ for $g=1$. The sets $\til\otimes v$ for $v\in\bim^{\otimes n}$ are called the \emph{$n$th level tiles} of the limit space. In particular, the sets of the form $\til\cdot g$ are the $0$th level tiles.

The \emph{tiling nerve} $\Xi_\alb$ is the nerve of the cover of $\limg$ by the $0$th level tiles $\til\cdot g$. It is a simplicial complex with the set of vertices $G$ in which a subset $A\subset G$ is a simplex if and only if $\bigcap_{g\in A}\til\cdot g$ is non-empty. 

The group $G$ acts naturally on the tiling nerve by automorphisms, since the action of $G$ on itself by right multiplication is simplicial on the tiling nerve. Note that a simplex of $\Xi_\alb$ may be fixed under an element $g\in G$ with a non-trivial action of $g$ on it. Therefore, we have to replace $\Xi_\alb$ by its barycentric subdivision if we want our standing assumption for $G$-complexes to be satisfied. In particular, the complex of groups $\Xi_\alb/G$ is actually the complex of groups $\Br(\Xi_\alb)/G$.

Simplices of $\Xi_\alb$, i.e., subsets $A\subset G$ such that $\bigcap_{g\in A}\til\cdot g\ne\emptyset$, are called \emph{adjacency sets}. The barycentric subdivision $\Br(\Xi_\alb)$ is the simplicial complex associated with the poset of adjacency sets.

The following is proved in~\cite[Proposition~6.1, Proposition~6.5]{nek:models}.

\begin{proposition}
A set $A\subset G$ containing the identity is an adjacency set if and only if there exists a sequence $A_n$ of finite subsets of the nucleus and a sequence $x_n\in\alb$ such that 
$A_n|_{x_n}=A_{n-1}$
and $A_0=A$. In particular, edges of $\Xi$ are sets of the form $\{g_1, g_2\}$ such that $g_1g_2^{-1}$ belongs to the nucleus $\nuke$. 
\end{proposition}

We have a natural equivariant simplicial map $J\colon\Xi_\alb\otimes\bim\arr\Xi_\alb$ defined on the set of vertices by the condition
\[J(g\otimes x)=g|_x\]
for all $g\in G$ and $x\in\alb$. The corresponding iteration $J^n\colon\Xi_\alb\otimes\bim^{\otimes n}\arr\Xi_\alb$ is also given by
\[J(g\otimes v)=g|_v\]
for $g\in G$ and $v\in\alb^n$. 

The following theorem is proved in~\cite[Theorem~6.7]{nek:models}.

\begin{theorem}
\label{th:topnuclcontr}
Let $\Xi_\alb$ be the tiling nerve for $(G, \bim)$, where $\alb\subset\bim$ is a basis.
Let $J\colon\Xi_\alb\otimes\bim\arr\Xi_\alb$ be the natural equivariant map. Then for all $n$ large enough, the map $J^n\colon\Xi_\alb\otimes\bim^{\otimes n}\arr\Xi_\alb$ is equivariantly homotopic to a contracting map.
\end{theorem}

The tiling nerve is universal among proper cocompact $G$-complexes in the following sense.

\begin{proposition}
\label{pr:intoXi}
Let $\Delta$ be a simplicial complex with a proper cocompact right action of $G$ on it. Let $\Xi_\alb$ be the tiling nerve for $(G, \bim)$. Then there exists $n$ and a $G$-equivariant simplicial map $F\colon\Delta\otimes\bim^{\otimes n}\arr\Br(\Xi_\alb)$.
\end{proposition}

\begin{proof}
Let $T\subset\Delta$ be a finite sub-complex such that $\bigcup_{g\in G}T\cdot g=\Delta$. Let $S$ be a finite self-similar set (i.e., such that $g|_x\in S$ for all $g\in S$ and $x\in\alb$) containing all elements $g\in G$ such that $T\cdot g$ and $T$ have non-empty intersections with a common simplex of $\Delta$. It exists, since the action of $G$ on $\Delta$ is proper. Consider the corresponding Rips complex, i.e., the simplicial complex $\Xi_S$ with the set of vertices $G$ in which a set $A\subset G$ is a simplex if and only if $gh^{-1}\in S$ for all $g, h\in A$. Note that $\Xi_\alb\subset\Xi_S$, since $S$ must contain the nucleus.

We have a $G$-equivariant map $I\colon\Xi_S\otimes\bim\arr\Xi_S$ given by
\[I(g\otimes x\cdot h)=g|_x\cdot h\]
for $g, h\in G$ and $x\in\alb$. It is simplicial, since $S|_\alb\subset S$.

For every vertex $x$ of $\Delta$, define $A_x=\{g\in G\colon x\in T\cdot g\}$.
Suppose that $A=\{x_1, x_2, \ldots, x_n\}$ is a simplex of $\Delta$. Let $g_1\in A_{x_i}$ and $g_2\in A_{x_j}$ for some $i, j\in\{1, 2, \ldots, n\}$. Then $x_1\in T\cdot g_i$ and $x_2\in T\cdot g_j$, hence the simplex $A\cdot g_j^{-1}$ of $\Delta$ intersects both $T$ and $T\cdot g_ig_j^{-1}$. Consequently, $A_{x_i}$ are simplices of $\Delta$ contained in one simplex $\bigcup_{i=1}^nA_{x_i}$. Define the map from the set of vertices of $\Delta$ to the barycentric subdivision of $\Xi_S$ mapping a vertex $x$ to the barycenter of $A_x$. The map is $G$-equviariant, and extends to a simplicial $G$-equivariant map $F_0\colon\Delta\arr\Br(\Xi_S)$.

It follows from the description of the simplices of the tiling nerve that there exists $n$ such that the range of the map $I^n\colon \Xi_S\otimes\bim^{\otimes n}\arr\Xi_S$ is equal to $\Xi_\alb$.
We get then the $G$-equviariant map $F=F_0\otimes I^n\colon \Delta\otimes\bim^{\otimes n}\arr\Br(\Xi_\alb)$.
\end{proof}

\subsection{From $G$-equivariant maps to virtual endomorphisms and back}
Let us describe how to move back and forth between the formalisms of equivariant maps and virtual endomorphisms.

Let $I\colon\X\otimes\bim\arr\X$ be an equivariant map. Consider the quotient spaces $\M_1=(\X\otimes\bim)/G$ and $\M_0=\X/G$. The map $I\colon\X\otimes\bim\arr\X$ induces a continuous map $\iota\colon\M_1\arr\M_0$. 

Suppose that $\xi_1\otimes x_1$ and $\xi_2\otimes x_2\in\X\otimes\bim$ belong to the same $G$-orbit. Then there exist elements $g, h\in G$ such that $\xi_1\cdot g=\xi_2$ and $g^{-1}\cdot x_1\cdot h=x_2$. In particular, $\xi_1$ and $\xi_2$ belong to the same $G$-orbit. Consequently, the correspondence $\xi\otimes x\mapsto \xi$ induces a well defined map $f\colon\M_1\arr\M_0$. This map is a covering if the action of $G$ on $\Xi$ is free.

We get a topological correspondence $f, \iota\colon\M_1\arr\M_0$, which we call the \emph{correspondence defined by} $I\colon\X\otimes\bim\arr\X$.

\begin{proposition}
Let $f, \iota\colon\M_1\arr\M_0$ be the correspondence defined by $I\colon\X\otimes\bim\arr\X$. Let $\M_n$ and the maps $f_n, \iota_n\colon\M_{n+1}\arr\M_n$ be the maps obtained by iterating it.

Then there exist homeomorphisms $\M_n\arr (\X\otimes\bim^{\otimes n})/G$ under which $\iota_n$ and $f_n$ correspond, respectively, to the maps induced by $I_n\colon\X\otimes\bim^{\otimes n+1}\arr\X\otimes\bim^{\otimes n}$ and by $\xi\otimes vx\mapsto \xi\otimes v$ for $v\in\bim^{\otimes n}$ and $x\in\bim$.
\end{proposition}

\begin{proof}
    Let us define the homomorphisms. Consider a point of $(\X\otimes\bim^{\otimes n})/G$ equal to the orbit of $\xi\otimes v$. We may assume that $v=x_1x_2\ldots x_n\in\alb^n$ for some fixed basis $\alb\subset\bim$, since $\alb^n$ is a transversal of right $G$-orbits. Then the sequence $(t_1, t_2, \ldots, t_n)$ of points of $\M_1$
    \begin{align*}
        t_1&=\xi\otimes x_1G,\\
        t_2&=I(\xi\otimes x_1)\otimes x_2G,\\
        t_3&=I^2(\xi\otimes x_1x_2)\otimes x_3G,\\
        \vdots\\
        t_n&=I^{n-1}(\xi\otimes x_1x_2\ldots x_{n-1})\otimes x_nG.
    \end{align*}
is a point of $\M_n$.

We get a map $\tau\colon\xi\otimes v G\mapsto (t_1, t_2, \ldots, t_n)$ from $\M_n$ to the space of sequences satisfying the conditions of the proposition. It is easy to check that it is well-defined. Since the maps $I^n$ are continuous, $\tau(\xi\otimes v)$ depends continuously on $\xi$ for all $v\in\alb^n$, hence it is a continuous map from $\M_n$ to the space of orbits of the correspondence.

Let us show that $\tau$ is bijective. Since it is a continuous map between two compact Hausdorff spaces, this will imply that it is a homeomorphism.
Suppose that $(t_1, t_2, \ldots, t_n)\in\M_1^n$ satisfies the condition of the proposition. Let $\xi\in\X$ be such that $f(t_1)=\xi\cdot G$. Then $t_1=\xi\otimes x_1\cdot G$ for some $x_1\in\alb$. Since $f(t_2)=\iota(t_1)$, $t_2=I(\xi\otimes x_1)\otimes x_2\cdot G$ for some $x_2\in\alb$. Continuing this way, we will get a word $v=x_1x_2\ldots x_n$ such that $(t_1, t_2, \ldots, t_n)=\tau(\xi\otimes v)$.
\end{proof}

Let $\Delta$ be a proper cocompact polytopal $G$-complex. Let us show how an equivariant map $I\colon\Delta\otimes\bim\arr\Delta$ defines a virtual endomorphism of complexes of groups. 

Let $\G_0=\Delta/G$ be the complex of groups associated with the action of $G$ on $\Delta$, see Definition~\ref{def:complexofaction}.
Consider the direct product $\Delta\times(\bim/G)$ and the projection on its first coordinate $P\colon\Delta\times(\bim/G)\arr\Delta$. The complex of groups $\mathcal{CG}(\Delta/G)$ acts on $\Delta\times(\bim/G)$ by its natural action on the second coordinate:
\[(g, \varphi)(A, xG)=(B\cdot g^{-1}, g\cdot xG),\]
where $xG\in\bim/G$ for $x\in\bim$, $g\in G$, and $\varphi$ is an attaching map from a face $A=\be(g, \varphi)$ to its subface $B$, so that $\en(g, \varphi)=B\cdot g^{-1}$. It is easy to check that this is an action. We can replace $\bim/G$ by a basis $\alb\subset\bim$, so that the action becomes
\[(g, \varphi)(A, x)=(B\cdot g^{-1}, g(x)).\]

Let $\G_1$ be the covering complex associated with this action, and let $f\colon\G_1\arr\G_0$ be the corresponding covering. The elements of $\mathcal{CG}(\G_1)^{(1)}$ are then triples $(g, \varphi, x)$ such that $(g, \varphi)\in\mathcal{CG}(\G_0)^{(1)}$ and $x\in\alb$, composed by the rule
\[(g_2, \varphi_2, g_1(x))\circ (g_1, \varphi_1, x)=(g_2g_1, (\varphi_2\cdot g_1)\circ\varphi_1, x).\]

Let $(g, \varphi, x)\in\mathcal{CG}(\G_1)^{(1)}$. Then $\varphi$ is an attaching map from a face $A$ to its subface $B$, and we have $\be(g, \varphi)=A$ and $\en(g, \varphi)=B\cdot g^{-1}$. Denote by $\varphi\otimes x$ the attaching map from the face $A\otimes x$ to its subface $B\otimes x$ in the complex $\Delta\times\alb$. Then $\tilde I(\varphi\otimes x)$ is the attaching map from $\tilde I(A\otimes x)$ to its subface $\tilde I(B\otimes x)$. 

It is checked directly that the map
\[\tilde\iota(g, \varphi, x)=(g|_x, \tilde I(\varphi\otimes x)),\]
is a functor that together with the maps of faces \[\iota(\xi, x)=I(\xi\otimes x)\colon (A, x)\mapsto I(A\otimes x)\] defines a morphism of complexes of groups $\iota\colon\G_1\arr\G_0$.

The virtual endomorphism $f, \iota\colon\G_1\arr\G_0$ is \emph{defined by the equivariant map $I\colon\Delta\otimes\bim\arr\Delta$.}

We can also go in the opposite direction: from virtual endomorphisms to equivariant maps. 

Let us show how to do it in the case of topological virtual endomorphisms. 
Let $f, \iota\colon\M_1\arr\M_0$ be a topological virtual endomorphism, and suppose that $\M_0$ is a path-connected space. Suppose that there exists $\epsilon>0$ such that every loop of diameter less than $\epsilon$ acts as identity on the preimage tree of its starting point. This is true for virtual endomorphisms of polytopal complexes, for locally expanding covering maps, and in the case when $\M_0$ is semi-locally simply connected. 

Choose a basepoint $t$ in $\M_0$. Let $\widehat{\M_0}$ be the space of paths in $\M_0$ starting in $t$ modulo the equivalence relation identifying two paths $\gamma_1, \gamma_2$ if $\gamma_2^{-1}\gamma_1$ is a loop that has trivial image in the iterated monodromy group $\img{(f, \iota)}$. Since small loops are trivial in the iterated monodromy group, we can introduce a topology on $\widehat{\M_0}$ for which a neighborhood of a path $\gamma$ ending in $z$ is the set of its continuations of $\delta\gamma$, where $\delta$ is a path of diameter $<\epsilon$. The map from such a path to its endpoint is, by definition, a homeomorphism of the neighborhood of the class of $\gamma$ in $\widehat{\M_0}$ with a neighborhood of the end of $\gamma$.

The iterated monodromy group $\img{(f, \iota)}$ acts naturally on $\widehat{\M_0}$ by appending loops to the beginnings of the paths. The natural equivariant map $I\colon\widehat{\M_0}\otimes\bim_{f, \iota}\arr\widehat{\M_0}$ is given by
\[I([\alpha]\otimes[(z, \ell)])=[\iota(\alpha_z)\ell],\]
where $\alpha_z$ is the lift of $\alpha$ by $f$ to a path starting in $z$. One can show that the correspondence associated with the map $I$ can be naturally identified with $f, \iota\colon\M_1\arr\M_0$. 

The same construction works also for the case when $f, \iota\colon\G_1\arr\G_0$ is a virtual endomorphism of a complex of groups such that the local groups of $\G_0$ are faithfully represented in the iterated monodromy group.  
We define $\widehat{\G_0}$ as the set of paths in $\G_0$ (in the sense of Definition~\ref{def:pathsinG})  starting in a fixed face $A$ of $\G_0$ modulo the same identification as in the topological case. The complex of groups $\G_0$ will be equivalent to the complex of groups associated with the action of the iterated monodromy group on $\widehat{\G_0}$, since the isotropy groups of $\G_0$ are faithfully represented in the iterated monodromy group.

\section{Polytopal substitutions and conformal dimension }
\label{subsec:pcontract} 

\subsection{Conformal dimension of a contracting group}

We give a short overview of the results of~\cite{nek:confdim}.

Let $(G, \bim)$ be a contracting group, and let $\alb\subset\bim$ be a basis.
Let $\nuke$ be the associated nucleus, and consider the graph $\G(\alb^n, \nuke)$ with the set of vertices $\alb^n$, the set of edges $\nuke\times\alb^n$, where $(g, v)$ connects $v$ to $g(v)$. A path in $\G(\alb^n, \nuke)$ is a sequence $(g_n, v_n)\ldots(g_2, v_2)(g_1, v_1)$, where $g_i(v_i)=v_{i+1}$. Therefore, it connects $v_1$ to $g_n\cdots g_2g_1(v_1)$. We have
\[(g_n\cdots g_2g_1)|_{v_1}=g_n|_{v_n}\cdots g_2|_{v_2}g_1|_{v_1}.\]
We say that the path is \emph{essential} if $(g_n\cdots g_2g_1)|_{v_1}\notin\nuke$.
More generally, if $N\subset G$ is any finite set containing $\nuke$, then we say that the path is \emph{$N$-essential} if $(g_n\cdots g_2g_1)|_{v_1}\notin N$.

A function $\phi\colon\nuke\times\alb^n\arr [0, \infty)$ is said to be \emph{admissible} (resp.\ \emph{$N$-admissible}) if for every essential (resp.\ $N$-essential) path $(g_n, v_n)\ldots (g_2, v_2)(g_1, v_1)$ we have $\phi(g_n, v_n)+\cdots+\phi(g_2, v_2)+\phi(g_1, v_n)\ge 1$.

Denote by $\overline M_p(n)$ (resp.\ $\overline M_{p, N}(n)$) the infimum of $\sum_{(g, v)\in\nuke\times\alb^n}\phi(g, v)^p$ over all admissible (resp.\ $N$-admissible) functions $\phi$. We call it the \emph{combinatorial $p$-modulus} of the graph $\G(\alb^n, \nuke)$.

We denote by $\overline M_\infty(n)$ (resp.\ $\overline M_{\infty, N}(n)$) the infimum of  $\max\phi(g, v)$ over all admissible (resp.\ $N$-admissible) functions $\phi$.

\begin{defi}
The \emph{asymptotic $p$-modulus} of $(G, \bim)$ is 
\[\theta_p=\lim_{n\to\infty}(\overline M_p(n))^{\frac{1}{np}}\]
for $p\in [1, \infty)$, and
\[\theta_\infty=\lim_{n\to\infty}(\overline M_\infty(n))^{\frac 1n}.\]
\end{defi}

The asymptotic $p$-modulus does not depend on the choice of the basis $\alb$. It is shown in \cite{nek:confdim} that $\theta_p$ satisfies, for all $1\le p<q\le\infty$, the following inequalities:
\begin{equation}
\label{eq:pqinequalities}
\theta_p\ge\theta_q,\qquad \mathsf{D}^{-1/p}\theta_p\le \mathsf{D}^{-1/q}\theta_q,
\end{equation}
where $\mathsf{D}=|\alb|=|\bim/G|$.

It is also shown that the value $\theta_p$ does not change if we replace $\overline M_p(n)$ by $\overline M_{p, N}(n)$ in its definition.

\begin{defi}
The \emph{conformal dimension} of $(G, \bim)$ is 
\[\cdim(G, \bim)=\inf\{p\colon\theta_p<1\}.\]
\end{defi}

The conformal dimension is usually defined as the infimum of Hausdorff dimensions of a natural class of metrics on the limit space $\lims$ of the group. The above definition is then a theorem (see~\cite[Theorem~4.12]{nek:confdim}).

Denote by $C(\lims)$ the $C^*$-algebra of continuous complex  functions on $\lims$. We will naturally identify it with the algebra of $G$-invariant complex functions $\limg\arr\C$. We denote by $C_0(\limg)$ the $C^*$-algebra of continuous functions $\limg\arr\C$ vanishing at infinity.

\begin{defi} For $p\ge 1$, denote by $\mathcal{D}_p(\limg)$ the sub-algebra of $C_0(\limg)$ consisting of functions $\phi\colon\limg\arr\C$ such that
\[\sum_{v\in\bim^*}|\phi(\xi_1\otimes v)-\phi(\xi_2\otimes v)|^p<\infty\]
for any $\xi_1, \xi_2\in\limg$.

For $p\ge 1$, denote by $\mathcal{D}_p(\lims)$ the sub-algebra of $C(\lims)$ consisting of $G$-invariant functions $\phi\colon\limg\arr\C$ such that 
\[\sum_{v\in\xs}|\phi(\xi_1\otimes v)-\phi(\xi_2\otimes v)|^p<\infty\]
for any $\xi_1, \xi_2\in\limg$.

For $p\ge 1$ and $\lambda\in (0, 1)$, denote by $D_{p, \lambda}(\limg)$ the sub-algebra of $C_0(\limg)$ consisting of functions $\phi$ such that for every pair of points $\xi_1, \xi_2\in\limg$ there exists $C>0$ such that
\[\sum_{v\in\bim^{\otimes n}}|\phi(\xi_1\otimes v)-\phi(\xi_2\otimes v)|^p\le C\lambda^{pn}\]
for all $n\ge 1$.

For $p\ge 1$ and $\lambda\in (0, 1)$, denote by $D_{p, \lambda}(\lims)$ the sub-algebra of $C(\lims)$ consisting of $G$-invariant functions on $\limg$ such that for any pair of points $\xi_1, \xi_2\in\limg$ there exists $C>0$ such that
\[\sum_{v\in\alb^n}|\phi(\xi_1\otimes v)-\phi(\xi_2\otimes v)|^p\le C\lambda^{pn}\]
for all $n\ge 1$.
\end{defi}

We obviously have $\mathcal{D}_{p, \lambda}(\lims)\subset\mathcal{D}_p(\lims)$ and $\mathcal{D}_{p, \lambda}(\limg)\subset\mathcal{D}_p(\limg)$.

It is shown in~\cite[Proposition~6.3]{nek:confdim} that the algebra $\mathcal{D}_p(\lims)$ (resp.\ $\mathcal{D}_{p, \lambda}(\lims)$) is dense in $C(\lims)$ if and only if the algebra $\mathcal{D}_p(\limg)$ (resp.\ $\mathcal{D}_{p, \lambda}(\limg)$) is dense in $C_0(\limg)$. 

Note that the defined algebras are invariant under compositions with Lipschitz maps $\C\arr\C$.

The following is~\cite[Theorem~6.4]{nek:confdim}.

\begin{theorem}
\label{th:densealgebras}
The conformal dimension $\cdim(G, \bim)$ is equal to the infimum of the values $p$ for which $\mathcal{D}_p(\lims)$ is dense in $C(\lims)$.

Moreover, for every $p>\cdim(G, \bim)$ the infimum of values of $\lambda$ for which $\mathcal{D}_{p, \lambda}(\lims)$ is dense in $C(\lims)$ is equal to $\theta_p$.
\end{theorem}

\subsection{$p$-contracting virtual endomorphisms}
Let $(G, \bim)$ be a contracting self-similar group. Let $\alb\subset\bim$ be a basis.
Consider a right proper cocompact action of $G$ on a metric space $(\X, d)$ by isometries.
Let $I\colon\X\otimes\bim\arr\X$ be a $G$-equivariant map. 
We assume that the maps $I_x(\xi)=I(\xi\otimes x)$ are Lipschitz for every $x\in\bim$ (equivalently, for every $x\in\alb$). 

Denote
\[Q_p(I)=\sup_{\xi_1, \xi_2\in\X, \xi_1\ne\xi_2}\frac{\left(\sum_{x\in\alb}d(I(\xi_1\otimes x), I(\xi_2\otimes x))^p\right)^{1/p}}{d(\xi_1, \xi_2)},\]
for $p\ge 1$. We also denote by $Q_\infty(I)$ the infimum of the common Lipschitz constant for the maps $I_x$ for all $x\in\alb$. 

It follows from $G$-invariance of the metric that $Q_p(I)$ does not depend on the choice of the basis $\alb$. It is finite, since  $Q_p(I)\le |\alb|\cdot Q_\infty(I)$.

We have 
\begin{multline*}
\sum_{w\in\alb^{n+m}}d(I^{n+m}(\xi_1\otimes w), I^{n+m}(\xi_2\otimes w))^p=\\ \sum_{u\in\alb^m}\quad\sum_{v\in\alb^n}d(I^n(I^m(\xi_1\otimes u)\otimes v), I^n(I^m(\xi_2\otimes u)\otimes v))^p\le \\ \sum_{u\in\alb^m}Q_p(I^n)^p d(I^m(\xi_1\otimes u), I^m(\xi_2\otimes u))^p\le Q_p(I^n)^pQ_p(I^m)^p d(\xi_1, \xi_2)^p. 
\end{multline*}
Consequently,
\[Q_p(I^{n+m})\le Q_p(I^n)Q_p(I^m)\]
for all $n, m$. It is also easy to check that the inequality holds also for $p=\infty$.

It follows then from the Fekete's Lemma that the limit \[\mathcal{Q}_p(I)=\lim_{n\to\infty}\sqrt[n]{Q_p(I^n)}\] exists for all $p\in [1, \infty]$ and is equal to $\inf_n\sqrt[n]{Q_p(I^n)}$.

\begin{defi}
\label{def:pcontractinggeneral}
    We say that a Lipschitz $G$-equivariant map $I\colon \X\otimes\bim\arr\X$ is \emph{$p$-contracting} if $Q_p(I)<1$. We say that it is \emph{asymptotically $p$-contracting} if $\mathcal{Q}_p(I)<1$.
\end{defi}

\begin{proposition}
\label{pr:ifpcontracting} Let $(G, \bim)$ be a contracting group.
For every contracting (i.e., $\infty$-contracting) $G$-equivariant map $I\colon \X\otimes\bim\arr\X$ and every $p\in [1, \infty]$, we have \[\theta_p\le \mathcal{Q}_p(I)\le Q_p(I).\]
\end{proposition}

\begin{proof}
Since $\mathcal{Q}_p(I)=\inf_n Q_p(I^n)^{1/n}$, it is enough to prove that $\theta_p\le Q_p(I)$ for all self-similar contracting groups $(G, \bim)$ and all $G$-equivariant contracting maps $I\colon \X\otimes\bim\arr\X$.

Let $\lambda>Q_p(I)$. Let $T\subset\X$ be a compact set intersecting every $G$-orbit. Since $I$ is contracting, we may assume that $I(T\otimes x)\subset T$ for all $x\in\alb$. Let $N\subset G$ be an arbitrary finite set such that $T\cap T\cdot g=\emptyset$ for all $g\in G\setminus N$.  It exists, since the action of $G$ on $\X$ is proper. We may assume that $\nuke\subset N$.

Let $r$ be the minimal distance between any two disjoint sets of the form $T, T\cdot g$ for $g\in G$.
Pick $\xi\in T$, and define the map $\phi_n$ on the set of edges of $\G(\alb^n, \nuke)$ by $\phi(s, v)=r^{-1}d(I^n(\xi\otimes v), I^n(\xi\otimes s\cdot v))$.

Suppose that a path $(s_1, v_1), (s_2, v_2), \ldots, (s_m, v_m)$ in $\G(\alb^n, \nuke)$  is $N$-essential, i.e., such that \[(s_m\cdots s_2s_1)|_{v_1}\notin N.\]
We have $s_{i-1}\cdots s_2s_1\cdot v_1=v_i\cdot (s_{i-1}\cdots s_2s_1)|_{v_1}$, so $d(I^n(\xi\otimes v_i), I^n(\xi\otimes s_i\cdot v_i))=d(I^n(\xi\otimes s_{i-1}\cdots s_2s_1\cdot v_1), I^n(\xi\otimes s_i\cdots s_2s_1\cdot v_1))$ by $G$-equivariance.

Therefore, $\sum_{i=1}^n\phi(s_i, v_i)\ge r^{-1}d(I^n(\xi\otimes v_1), I^n(\xi\otimes v_1\cdot (s_m\cdots s_2s_1)|_{v_1})$. We have $\xi\otimes v_1\in\til$ and $\xi\otimes v_1\cdot (s_m\cdots s_2s_1)|_{v_1}\in\til\cdot g$ for $g=(s_m\cdots s_2s_1)|_{v_1}\notin N$. Consequently, the sum of the values of $\phi_n$ along the path is greater than 1, so the function $\phi_n$ is $N$-admissible.

We have $\sum_{s\in\nuke, v\in\alb^n}\phi_n(v)^p=\sum_{s\in\nuke, v\in\alb^n}r^{-p}d(I^n(\xi\otimes v), I^n(\xi\otimes s\cdot v))^p\le r^{-p}Q_p(I^n)^p\sum_{s\in\nuke}d(\xi, \xi\cdot s)^p$. Consequently,
\[\overline{M}_{p, N}(n)\le Q_p(I^n)^p\cdot r^{-p}\sum_{s\in\nuke}d(\xi, \xi\cdot s)^p,\]
which implies $\theta_p\le Q_p(I^n)^{1/n}\le Q_p(I)$.
\end{proof}

Suppose now that $\Delta$ is a proper cocompact connected polytopal $G$-complex. Let 
$I\colon \Xi\otimes\bim\arr\Xi$ be a $G$-equivariant smooth morphism. 

For every $x\in\bim$, every face $A$, and every point $\xi\in A$ (seen as a point in the corresponding ambient space $\R^k$ of the polytope $A$) we have the associated differential $DI_v$ from the tangent space $\R^k$ at $\xi$ to the tangent space to $I_v(\xi)\in\tilde I(A)$. Let $\|\cdot \|$ be the Euclidean norm on the tangent spaces.
Define then
\[DQ_p(I)=\sup
\frac{\left(\sum_{x\in\alb}\|DI_x(\vec a)\|^p\right)^{1/p}}{\|\vec a\|},\]
where the supremum is taken over all faces $A$ of $\Delta$, all points $\xi\in A$, and all non-zero tangent vectors $\vec a$.
We also define 
\[\mathcal{DQ}_p(I)=\lim_{n\to\infty}DQ_p(I^n)^{1/n}.\]

Replacing in the definition of $DQ_p(I)$ the Euclidean norm on $\R^k$ by any other norm will not change $\mathcal{DQ}_p(I)$, since $\Delta$ is finite-dimensional. In particular, it may be convenient to consider the $\ell^p$-norm on $\R^k$. Then $DQ_p(I)$ becomes the supremum over all points $\xi$ of all faces $A$ of $\Delta$ of the $\ell^p$-norm of the differential at $\xi$ of the map
\[A\to \prod_{v\in\alb^n}\tilde I(A\otimes v)\colon \xi\mapsto (I(\xi\otimes v))_{v\in\alb^n}.\]

Let us translate the definition of $DQ_p$ and $\mathcal{DQ}_p$ to the setting of virtual endomorphisms of complexes of groups.

Let $f, \iota\colon\G_1\arr\G_0$ be the virtual endomorphism of complexes of groups.
Let $A$ be a face of $\G_0$. Consider all of its $f$-preimages $B_1, B_2, \ldots, B_d$ in the sense of coverings of complexes of groups, i.e.,  each face $B$ such that $f(B)=A$ is repeated $[G_A:\tilde f(G_B)]$ times. The maps $f\colon B_i\arr A$ are isometries, so we can consider their inverses $f_i^{-1}\colon A\arr B_i$. Consider the map $F_A\colon A\arr\prod_{i=1}^d\tilde\iota(B_i)$ given by
\[F_A(\xi)=(\iota(f_1^{-1}(\xi), \iota(f_2^{-1}(\xi)), \ldots, \iota(f_d^{-1}(\xi))).\]
Denote by $DQ_p(f, \iota)$ the supremum of the $\ell^p$-norm of the differential of $F_A$ over all $A$ and all points $\xi\in |A|$. 

Next proposition follows directly from the definitions.

\begin{proposition}
Suppose that a virtual endomorphism $f, \iota\colon\G_1\arr\G_0$ is associated with an equivariant map $I\colon\Delta\otimes\bim\arr\Delta$. Then $DQ_p(I)=DQ_p(f, \iota)$.
\end{proposition}

Let us relate the two notions of $\ell^p$-contraction.

\begin{proposition}
\label{prop:DQQ}
Let $I\colon\Delta\otimes\bim\arr\Delta$, and suppose that 
$DQ_p(I)$ is defined using the Euclidean norm on the tangent spaces.

Then $Q_p(I)\le DQ_p(I)$.
\end{proposition}

\begin{proof}
Suppose that $\gamma\colon [0, L]\arr\Xi$ is a geodesically parametrized path (i.e., $\left\|\frac{d\gamma(t)}{dt}\right\|=1$ for a.e.\ $t\in [0, L]$. In particular, $l(\gamma)=L$.

Then for every $x\in\alb$ we have, by Jensen's inequality,
\begin{multline*}
l(I(\gamma\otimes x))^p=\left(\int_0^L\left\|\frac{d I(\gamma(t)\otimes x)}{dt}\right\|\;dt\right)^p\le\\
L^{p-1}\int_0^L\left\|\frac{d I(\gamma(t)\otimes x)}{dt}\right\|^p\;dt=
L^{p-1}\int_0^L\left\|I\left(\frac{d\gamma(t)}{dt}\otimes x\right)\right\|^p\;dt.
\end{multline*}
Consequently,
\[\sum_{x\in\alb}l(I(\gamma\otimes x))^p\le 
L^{p-1}\int_0^L\sum_{x\in\alb}\left\|I\left(\frac{d\gamma(t)}{dt}\otimes x\right)\right\|^p\;dt\le DQ_p(I)L^p,\]
which finishes the proof.
\end{proof}

We get the following corollary of Propositions~\ref{pr:ifpcontracting}, \ref{prop:DQQ}, and inequalities~\eqref{eq:pqinequalities}.

\begin{corollary}
Let $f, \iota\colon\G_1\arr\G_0$ be a virtual endomorphism of a complex of groups.
If  $DQ_p(f, \iota)<1$, for $p\ge 1$, then 
\[\cdim(G, \bim)\le \left(\frac 1p-\frac{\log DQ_p(I)}{n\log \mathsf{D}}\right)^{-1},\]
where $\mathsf{D}=|\bim/G|$.
\end{corollary}

\begin{example}
Consider a covering map $f\colon\R^n/\Z^n\arr\R^n/\Z^n$ defined by a matrix $A$ with integer entries. We consider it as a topological virtual endomorphism $f, \iota\colon\R^n/\Z^n\arr\R^n/\Z^n$, where $\iota$ is the identity map. After subdividing the cube $[0, 1]^n$, we can transform this correspondence into a polytopal substitution. The virtual endomorphism is contracting if and only if the matrix $A$ is expanding, i.e., all its eigenvalues are less than one in absolute value.
We have $DQ_p(\iota)=|\det A|\cdot \|A^{-1}\|^p$, and if $f^n, \iota^n\colon\R^n/\Z^n\arr\R^n/\Z^n$ is the $n$th iterate of the virtual endomorphism, then \[DQ_p(\iota^n)=|\det A|^n\cdot\|A^{-n}\|^p.\]
We get 
\[\mathcal{DQ}_p(\iota)=|\det A|^p\cdot \rho,\]
where $\rho$ is the spectral radius of  $A^{-1}$. In fact, we have $\theta_p=\mathcal{DQ}_p$ in this case, see~\cite{nek:confdim}.
\end{example}

\begin{example}
Consider the substitution from Example~\ref{ex:sierptriangle}.
Let use the substitution from Example~\ref{ex:sierptriangle} without iterating it, 
and minimize $DQ_p$ inside the homotopy class of $\iota$. Suppose that the vertices of $\G_1$ labeled by $O$ are mapped by $\iota$ to the points on distance $x$ from the vertex $O$ of $\G_0$, and the central vertices of $\G_1$ labeled by $A, B, C$ are mapped to the vertex $O$ of $\G_0$. Assuming that the lengths of the edges of $\G_0$ are equal to 1, we get
\[DQ_p=\left(2x^p+(1-x)^p\right)^{1/p}.\]
The minimal value at $x=\left(1+2^{\frac{1}{p-1}}\right)^{-1}$ gives the estimate.
\[\theta_p\le DQ_p=\left(1+2^{\frac{1}{1-p}}\right)^{\frac 1p-1},\]
which proves that $\cdim(G, \bim)=1$. One can get better estimates by considering higher iterates $\iota^n$.

It is shown in~\cite{HPS:penergy} that \[\theta_p\le\left(1+\frac{2^{2-1/(p-1)}}{1+\sqrt{1+2^{3-1/(p-1)}}}\right)^{\frac 1p-1}.\]

It is not hard to see that $\theta_\infty=1/2$. It is also known that $\theta_2=\sqrt{3/5}$, see~\cite[Example~3.1.5]{kigami:book} and~\cite[Lemma~1.1]{strichartz:taylor}.
\end{example}

\subsection{Existence of $p$-contracting virtual endomorphisms}

\begin{theorem}
\label{th:existpcontracting}
Let $(G, \bim)$ be a contracting group. Let $\Delta$ be a proper cocompact simplicial $G$-complex, and let $I\colon\Delta\otimes\bim\arr\Delta$ be a contracting equivariant map.
    
If $p>\cdim(G, \bim)$ and $\lambda>\theta_p$, then there exists $n\ge 1$ such that $I^n\colon\Delta\otimes\bim^{\otimes n}\arr\Delta$ is $G$-equivariantly homotopic to a map $J\colon\Delta\otimes\bim^{\otimes n}\arr\Delta$ such that $DQ_p(J)\le\lambda^n$.
\end{theorem}

\begin{proof}
We will construct the map $J$ in a formally more general setting than the conditions of the theorem.
Let $\Delta$ be a proper cocompact simplicial $G$-complex. Let $\limg$ be the limit $G$-space, and suppose that there exists an equivariant map $L\colon\limg\arr\Delta$.  This is true when there is a contracting equivariant map $\Delta\otimes\bim\arr\Delta$, by Theorem~\ref{th:limgconvergence}.

The \emph{star} $\Sigma_\zeta$ of a vertex $\zeta$ of $\Delta$ is the union of the simplices of $\Br(\Delta)$ to which $\zeta$ belongs. The set of vertices of the star is therefore the set of barycenters of simplices containing $\zeta$.

We say that a vertex $\zeta$ of $\Delta$ is a \emph{neighbor} of a point $\xi\in |\Delta|$ if $\xi$ belongs to $\Sigma_\zeta$.

\begin{lemma}
\label{lem:neighborvertex}
Let $\Delta$ be a proper cocompact $G$-complex. 
There exists $\epsilon>0$ such that if $T$ is a subset of $|\Delta|$ of diameter $\le\epsilon$, then the set of vertices that are neighbors of some points belonging to $T$ is a simplex of $\Delta$.
\end{lemma}

\begin{proof}
A set of vertices $A$ of $\Delta$ is a simplex if and only if the stars of the vertices of $A$ have a non-empty intersection. Therefore, it is enough to take $\epsilon$ to be smaller than the distance between any star $S$ and any non-empty intersection $B$ of stars of $\Delta$ such that $S\cap B=\emptyset$.
\end{proof}

Choose $\lambda_0\in (\theta_p, \lambda)$. Then the algebra $\mathcal{D}_{p, \lambda_0}(\limg)$ is dense in $C_0(\limg)$, see Theorem~\ref{th:densealgebras}. 
We can find a collection of open relatively compact sets $U_\zeta\subset|\Delta|$, labeled by vertices $\zeta$ of $\Delta$, such that $U_\zeta$ contains the star $\Sigma_\zeta$ of $\zeta$, and for any set $\{\zeta_1, \ldots, \zeta_k\}$ of vertices, the intersection $\bigcap_i\overline{U_{\zeta_i}}$ is non-empty if and only if the intersection $\bigcap_i\Sigma_{\zeta_i}$ is non-empty (which is equivalent to the condition that $\{\zeta_1, \zeta_2, \ldots, \zeta_k\}$ is a simplex of $\Delta$).

Moreover, since the action of $G$ on $\Delta$ is proper and cocompact, we can choose them so that we have $U_{\zeta\cdot g}=U_\zeta\cdot g$ for all $g\in G$ and all vertices $\zeta$.

Then $\{L^{-1}(U_\zeta)\}$ is a $G$-invariant open cover of $\limg$. For every $U_\zeta$, we can find a non-negative function $\phi_\zeta\colon\limg\arr\R$ supported in $L^{-1}(U_\zeta)$ and positive on $L^{-1}(\Sigma_\zeta)$. We can choose them in such a way that $\phi_{\zeta\cdot g}(\xi\cdot g)=\phi_{\zeta}(\xi)$ for every $\xi\in\limg$ and every vertex $\zeta$ of $\Delta$. 

After approximating $\phi_\zeta$ by an element of $\mathcal{D}_{p, \lambda_0}(\limg)$ and composing the approximation with a Lipschitz function $\R\arr\R$, we may assume that $\phi_\zeta\in\mathcal{D}_{p, \lambda_0}(\limg)$. We also can do it in a $G$-equivariant way.

The sum $\phi$ of the functions $\phi_\zeta$ for all vertices of $\Delta$ is $G$-invariant, belongs to $\mathcal{D}_{p, \lambda_0}(\lims)$, and is strictly positive. Since its set of values belongs to a compact sub-interval of $(0, \infty)$ and the function $t\mapsto t^{-1}$ is Lipschitz on the set of the values of $\phi$, the function $\frac{1}{\phi}$ belongs to $\mathcal{D}_{p, \lambda_0}(\lims)$. Multiplying the functions $\phi_\zeta$ by $\phi^{-1}$, we get a $G$-invariant partition of unity $\{\psi_\zeta\}$ subordinate to the cover $\{U_\zeta\}$, such that $\psi_\zeta\in\mathcal{D}_{p, \lambda_0}(\limg)$.

Let $\xi\in\limg$. Let $\{\zeta_1, \zeta_2, \ldots, \zeta_k\}$ be the set of all vertices of $\Delta$ such that $L(\xi)\in U_{\zeta_i}$ for every $i$. Then $\psi_{\zeta_1}(\xi)+\psi_{\zeta_2}(\xi)+\cdots\psi_{\zeta_k}(\xi)=1$, and the set $\{\zeta_1, \zeta_2, \ldots, \zeta_k\}$ is a simplex of $\Delta$. Denote by $\Psi(\xi)$ the point of the geometric realization of the simplex $(\zeta_1, \zeta_2, \ldots, \zeta_k)$ with the barycentric coordinates $(\psi_{\zeta_1}(\xi), \psi_{\zeta_2}(\xi), \ldots, \psi_{\zeta_k}(\xi))$. Then $\Psi\colon\limg\arr\Delta$ is a $G$-equivariant continuous map.

Choose an arbitrary $G$-equivariant map $\Lambda$ from the set of vertices of $\Delta$ to $\limg$. 

Let $m$ be large enough (to be chosen later). Consider the map $J_m$ from the set of vertices of $\Delta\otimes\bim^{\otimes m}$ to $|\Delta|$ defined by
\[J_m(\zeta\otimes v)=\Psi(\Lambda(\zeta)\otimes v)\]
for $v\in\bim^{\otimes m}$ and a vertex $\zeta$ of $\Delta$.

Denote by $\mathcal{U}$ the set of closures of the sets $U_\zeta$ for all vertices $\zeta$ of $\Delta$. Let $\epsilon>0$ be smaller than the distance from any element $V\in\mathcal{U}$ to any non-empty intersection $V_1\cap V_2\cap\cdots\cap V_k$ of elements of $\mathcal{U}$ disjoint from it. If $V\subset|\Delta|$ is a set of diameter less than $\epsilon$, and $\zeta_1, \zeta_2, \ldots, \zeta_k$ are all vertices of $\Delta$ such that $V\cap U_{\zeta_i}\ne\emptyset$, then $\bigcup_i\overline{U_{\zeta_i}}\ne\emptyset$, hence $\{\zeta_1, \zeta_2, \ldots, \zeta_k\}$ is a simplex of $\Delta$.

Since $L$ is equivariant, and the action of $G$ on $\Delta$ is proper and cocompact, there exists $D$ such that for every simplex $A$ of $\Delta$ the diameter of the set $\Lambda(A)$ is less than $D$.
Consequently, the diameter of $\Lambda(A)\otimes v\subset\limg$ is smaller than $C\theta^mD$ for some constants $C>0$ and $\theta\in (0, 1)$. It follows that if $m$ is large enough, then the diameter of $\Lambda(A)\otimes v$ is less than $\epsilon$. 

Then the set of points $\{\zeta_1, \ldots, \zeta_k\}$ such that $\Lambda(\zeta)\otimes v \cap \overline{U_{\zeta_i}}\ne\emptyset$ is a simplex of $\Delta$. Consequently, $\Psi(\Lambda(A)\otimes v)$ is contained in the simplex $\{\zeta_1, \ldots, \zeta_k\}$. It follows that we can extend the map $J_m$ from the set of vertices to an equivariant affine morphism of complexes $J_m\colon\Delta\otimes\bim^{\otimes m}\arr\Delta$.

Let us estimate $DQ_p(J_m)$. 
Let $A=\{\zeta_1, \zeta_2, \ldots, \zeta_k\}$ be a simplex of $\Delta$, and consider the map $J_m(\cdot\otimes v)\colon A\arr |\Delta|$ induced by it. It is the affine extension of the map 
\[\zeta_i\mapsto \Psi(\Lambda(\zeta_i)\otimes v).\]

Consider a tangent vector $\vec a=(a_1, a_2, \ldots, a_k)$ to $\zeta$ where all coordinates are equal to $0$ except for $a_r=1$ and $a_s=-1$ for some $r, s\in\{1, 2, \ldots, k\}$. Note that vectors of this form span the tangent space. The derivative of $\Psi(\cdot\otimes v)$ in the direction of this vector is equal to
\[\left(\psi_{\xi_i}(\Lambda(\zeta_r)\otimes v)-\psi_{\xi_i}(\Lambda(\zeta_s)\otimes v)\right)_{i=1}^m,\]
for some vertices $\xi_1, \xi_2, \ldots, \xi_m$ of $\Delta$.

Let us choose a finite set $\{A_1, A_2, \ldots, A_t\}$ of $G$-orbit representatives of simplices of $\Delta$. Then the union of all sets of the form $\Psi(\Lambda(A_i)\otimes v)$ for all $v\in\xs$ is contained in a compact set, so the set of simplices of $\Delta$ intersecting them is finite. Let us denote the union of their sets of vertices by $K$. Assuming that $A$ belongs to the chosen set of representatives (which we can do by equivariance), we get
\[\sum_{v\in\alb^m}\|D_{\vec a}J_m(\cdot\otimes v)\|^p_p=\sum_{\xi\in K}|\psi_\xi(\Lambda(\zeta_r)\otimes v)-\psi_\xi(\Lambda(\zeta_s)\otimes v)|^p\le C_1\lambda_0^{mp}\]
for some constant $C_1$, not depending on $\vec a$, $A$, and $m$, since $\psi_\xi\in \mathcal{D}_{p, \lambda_0}(\limg)$.

Since the vectors $\vec a$ we considered span the tangent space to $A$, there exists a constant $C_2$, not depending on $m$, such that we have \[\sum_{v\in\alb^n}\|D_{\vec a}J_m(\cdot\otimes v)\|_p^p\le C_2\lambda_0^{mp}\] for all tangent vectors $\vec a$. Consequently, choosing $m$ large enough, we get
\[DQ_p(J_m)\le\lambda^m,\]
as $\lambda_0<\lambda$.

It remains to show that $\xi\otimes v\mapsto J_m(\xi\otimes v)$ is homotopic to $I^m$. 
Let $A$ be a simplex of $\Delta$. Since the diameter of $A\cup L(\Lambda(A))$ is uniformly bounded, and $I$ is contracting, the diameter of the set $I^m(A\otimes v)\cup I^m(L(\Lambda(A))\otimes v)$ goes to $0$ as $m\to\infty$. In particular, we may assume that it is less than $\epsilon$. We have $I^m(L(\Lambda(A))\otimes v)=L(\Lambda(A)\otimes v)$. By definition of $\Psi$, every point $\xi$ of $J_m(A\otimes v)=\Psi(\Lambda(A)\otimes v)$ belongs to a simplex of $\Delta$ whose vertices $\zeta$ satisfy $L(\Lambda(\xi)\otimes v)\in U_\zeta$. It follows that $I^m(A\otimes v)$ and $J_m(A\otimes v)$ are contained in the geometric realization of a common simplex of $\Delta$. This implies that the maps $I^m$ and $J_m$ are equivariantly homotopic.

If $\Delta$ is not a simplicial complex, then the same arguments construct a map $J_m$ from $\Delta\otimes\bim^{\otimes m}$ to the second barycentric decomposition $\Br(\Br(\Delta))$. By the same contraction argument as above, we show that the image of every face of $\Delta$ is contained in a simplex of the barycentric decomposition, hence we get a morphism of $G$-complexes $J_m\colon\Delta\otimes\bim^{\otimes m}\arr\Delta$ equivariantly homotopic to $I^m$.
\end{proof}

As a direct corollary of Theorem~\ref{th:existpcontracting} and Proposition~\ref{pr:ifpcontracting}, we get the following characterization of asymptotic modulus and conformal dimension of the limit space.

\begin{theorem}
\label{th:asymptoticDQvirtend}
Let $f, \iota\colon\G_1\arr\G_0$ be a contracting virtual endomorphism of complexes of groups, and let $(G, \bim)$ be its iterated monodromy group. Let $\widehat{DQ}_p(\iota^n)$ be the infimum of $DQ_p(f^n, \widehat{\iota^n})$ over all morphisms $\widehat{\iota^n}\colon\G_n\arr\G_0$ homotopic to $\iota^n$. Then
\[\theta_p=\lim_{n\to\infty}\widehat{DQ}_p(n)^{1/n}.\]
\end{theorem}

\subsection{$p$-energies of graphs}

It is sometimes convenient to reformulate the definition of $DQ_p(f, \iota)$ in terms of lengths of images of edges instead of derivatives.

\begin{defi}
Let $f\colon\G_1\arr\G_0$ be a morphism of polytopal complexes of groups.
Then its \emph{$p$-energy} $E_p(f)$ is the sum of $p$th powers of distances between the images of the endpoints of edges (i.e., one-dimensional faces) of $\G_1$.

Denote by $\\widehat{E}_p(f)$ the infimum of the $p$-energies over all morphisms from the one-skeleton of $\G_1$ (seen as a graph of groups) to $\G_0$ homotopic to the restriction of $f$ to the one-skeleton.
\end{defi}

\begin{theorem}
Let $f, \iota\colon\G_1\arr\G_0$ be a contracting virtual endomorphism of a complex of groups, and let $(G, \bim)$ be its iterated monodromy group. Then
\[\lim_{n\to\infty}\widehat{E}_p(\iota^n)^{1/n}=\theta_p^p.\]
\end{theorem}

\begin{proof}
Similar arguments to the ones in the proof of Proposition~\ref{pr:ifpcontracting} show that \[\theta_p^p=\lim_{n\to\infty}\overline M_p(n)^{1/n}\le\liminf_{n\to\infty}\widehat{E}_p(\iota^n)^{1/n}.\]

On the other hand, Theorem~\ref{th:asymptoticDQvirtend} implies that $\limsup_{n\to\infty}\widehat{E}_p(\iota^n)^{1/np}\le\theta_p$.
\end{proof}

\begin{example}
Consider the substitution shown in Figure~\ref{fig:zsqponemodel}. The map $f$ is simplicial and maps a vertex labeled by a latter of the form $X_i$ of $\G_1$ to the vertex labeled by the corresponding label $X$.  The map $\iota\colon\G_1\arr\G_0$ acts by
\[B_0\mapsto A,\quad C_0\mapsto B, \quad B_1\mapsto C,\quad O_0\mapsto O,\]
maps the vertex $A$, and the edge $O_1C_1$ into the iterior of the edge $OC$. It maps the whole edge $O_1C_1$ to a single point of $OC$.

This substitution is subhyperbolic, and its limit dynamical system is topologically conjugate to the action of the polynomial $z^2+i$ on its Julia set.

\begin{figure}
\centering
\includegraphics{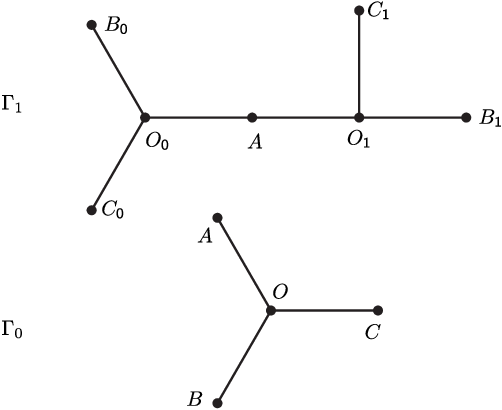}
\caption{The Hubbard tree for $z^2+i$}
\label{fig:zsqponemodel}
\end{figure}

Let $\G_1'$ be the subgraph of $\G_1$ consisting of all vertices and edges except for $C_1$ and $O_1C_1$, i.e., the convex hull of $B_0, B_1, C_0$. Then the restriction of $\iota$ is a homeomorphism $\iota\colon\G_1'\arr\G_0$. The complement of $\G_1$ is mapped by $\iota$ to a single point.

By iterating the substitution $f, \iota\colon\G_1'\arr\G_0$, we get a sequence of subgraphs $\G_n'$ of $\G_n$ such that $\iota^n$ maps $\G_n'$ homeomorphically to $\G_0$, and it maps each connected component of the complement to a singleton. Due to the latter condition, the minimal energy $\mathcal{E}_p(\iota^n)$ in the homotopy class of the morphism of complexes of groups $\iota^n\colon\G_n\arr\G_0$ coincides with the minimal energy in the homotopy class of the homeomorphism $\iota^n\colon\G_n'\arr\G_0$, relative the degree 1 vertices of the graph $\G_n'$.

The graph $\G_n'$ is a tripod of chains of edges, mapped by $f^n$ (``labeled'' in the terminology of polytopal substitutions) to the edges $OA, OB$, and $OC$ of $\G_0$.  According to the substitution, the labeled graph $\G_n'$ is obtained from $\G_{n-1}'$ by replacing the edges labeled by $OA$ and $OB$ by edges labeled by $OB$ and $OC$, respectively, and replacing every edges labeled by $OC$ by a chain of three edges labeled by $OA$, $OA$, and $OB$.

Consequently, if a leg of the tripod $\G_{n-1}'$ has $x_{n-1}, y_{n-1}, z_{n-1}$ edges labeled by $OA, OB, OC$, respectively, then the corresponding leg of $\G_n'$ has $2z_{n-1}$, $x_{n-1}+z_{n-1}$, $y_{n-1}$ edges labeled by $OA, OB, OC$, respectively. It follows that the vectors $\left(\begin{array}{c}x_n\\ y_n\\ z_n\end{array}\right)$ for the three legs of $\G_n'$ are the columns of the matrix
\[\left(\begin{array}{ccc} 0 & 0 & 2\\
1 & 0 & 1\\ 0 & 1 & 0\end{array}\right)^n,\]
so the lengths of the legs of $\G_n'$ are the coordinates of the row-vector
\[(a_n, b_n, c_n)=(1, 1, 1)\left(\begin{array}{ccc} 0 & 0 & 2\\
1 & 0 & 1\\ 0 & 1 & 0\end{array}\right)^n.\]

This implies 
\[\widehat{E}_p(\iota^n)=a_n^{1-p}+b_n^{1-p}+c_n^{1-p}\]
and
\[\theta_p=\lim_{n\to\infty}(a_n^{1-p}+b_n^{1-p}+c_n^{1-p})^{\frac 1{np}}=\lambda^{1-\frac 1p},\]
where $\lambda\approx 1.5214$ is the spectral radius of the matrix above, i.e., the real root of $x^3-x-2$.

Similar arguments can be used to show that, for the basilica substitution from Example~\ref{ex:basilicasubst}, we have $\theta_p=\sqrt{2}^{\frac 1p-1}$.
\end{example}

\begin{example}
Consider Example~\ref{ex:gromovsubst}. The graph $\G_{n+1}$ is obtained by taking four copies of $\G_n$ (equal to the preimages of the edges of $\G_0$ under $\iota^n$) and pasting them together in the same way as the edges are pasted together to construct $\G_1$ from $\G_0$. It follows that $\widehat{E}_p(\iota^n)$ is the minimum of 
\[(a_1^p+2a_2^p+a_3^p)\widehat{E}_p(\iota^{n-1})\]
over all triples $(a_1, a_2, a_3)$ of non-negative real numbers satisfying $a_1+a_2+a_3=1$. Consequently, \[\widehat{E}_p(\iota^n)=\left(2+2^{\frac{1}{1-p}}\right)^{n(1-p)},\]
which implies 
\[\theta_p=\left(2+2^{\frac{1}{1-p}}\right)^{\frac 1p-1}.\]
\end{example}

\def\cprime{$'$}\def\ocirc#1{\ifmmode\setbox0=\hbox{$#1$}\dimen0=\ht0
  \advance\dimen0 by1pt\rlap{\hbox to\wd0{\hss\raise\dimen0
  \hbox{\hskip.2em$\scriptscriptstyle\circ$}\hss}}#1\else {\accent"17 #1}\fi}
  \def\cprime{$'$}
\providecommand{\bysame}{\leavevmode\hbox to3em{\hrulefill}\thinspace}
\providecommand{\MR}{\relax\ifhmode\unskip\space\fi MR }
\providecommand{\MRhref}[2]{%
  \href{http://www.ams.org/mathscinet-getitem?mr=#1}{#2}
}
\providecommand{\href}[2]{#2}

\end{document}